%% file: main.tex
\documentclass[oneside, headsepline, 10pt, BCOR=1mm,DIV=18, titlepage=false]{scrartcl}

\usepackage{amsfonts}
\usepackage{amssymb}
\usepackage{amsthm}
\usepackage{amsmath, amsfonts, amsxtra}
\usepackage{esint}
\usepackage{verbatim} %

\input{environments}

\input{commands}

\numberwithin{equation}{section}%

\title{Interior and Boundary-Regularity for\\Fractional Harmonic Maps on Domains}
\author{Armin Schikorra}
\date{\today}

 \usepackage{fancyhdr}

\begin{document}
\maketitle

\thispagestyle{empty}
\begin{abstract}
\noindent%
\input{abstract}
\end{abstract}

\tableofcontents

 \pagestyle{fancy}                                 %Erstellen der Kolumnentitel

 \newcommand{\chaptermark}[1]{\markboth{\it \thechapter \ #1}{}}
  %Definition des Textes, der links oben erscheinen soll
 \renewcommand{\sectionmark}[1]{\markright{\it \thesection \ #1}}
  %Definition des Textes, der rechts oben erscheinen soll
 %\fancyheadoffset[RO]{30pt}
 \fancyhead{}                                %Kopfzeilen leeren
 %\fancyhead[RE]{} %	{A. Schikorra}
 %\fancyhead[LE]{\thepage}                        %links oben
 %\fancyhead[RO]{\thepage}                        %rechts oben
 %\fancyhead[LO]{blabla}                        %rechts oben
\fancyhead[C]{\thepage}
\fancyfoot{}                                %Fußzeilen leeren
\fancyfoot[C]{}
 %\fancyfoot[LE,RO]{\rule{17pt}{1.5pt} \\ \thepage}
  %Seitenzahlen jeweils links bzw. rechts unten

\input{intro}
\newpage

\input{facts}

\newpage

\input{lowerorder}

\newpage

\input{cutoff}

\newpage

\input{energy}

\newpage
% \input{sketch}
% \newpage

%\input{bel2}
%\newpage

\input{regproof}
\newpage

\input{boundary}

\newpage

\input{iteration}

\newpage

\bibliographystyle{alpha}%
\bibliography{bib}%
\vspace{2em}
\begin{tabbing}
\quad\=Armin Schikorra\\
\>RWTH Aachen University\\
\>Institut f\"ur Mathematik\\
\>Templergraben 55\\
\>52062 Aachen\\
\>Germany\\
\\
\>email: schikorra@instmath.rwth-aachen.de
\end{tabbing}
\end{document}

%% file: environments.tex
\newtheorem{theorem}{Theorem}[section]%
\newtheorem{lemma}[theorem]{Lemma}%
\newtheorem{proposition}[theorem]{Proposition}%
\newtheorem{corollary}[theorem]{Corollary}%
\newtheorem{definition}[theorem]{Definition}%
\newtheorem{remark}[theorem]{Remark}%

\makeatletter

\newenvironment{proofP}[1]{\par
  \pushQED{{\small\it Proposition {\rm #1}} {\leavevmode
  \hbox to.77778em{%
  \hfil\vrule
  \vbox to.675em{\hrule width.6em\vfil\hrule}%
  \vrule\hfil}}}%
  \normalfont \topsep6\p@\@plus6\p@\relax
  {\it \textbf{\proofname { of Proposition {\rm #1}\@addpunct{.}}}\newline\rm}
}{%
   \begin{flushright}\popQED\end{flushright}\@endpefalse
}
\makeatother

\makeatletter

\newenvironment{proofT}[1]{\par
  \pushQED{{\small\it Theorem {\rm #1}} {\leavevmode
  \hbox to.77778em{%
  \hfil\vrule
  \vbox to.675em{\hrule width.6em\vfil\hrule}%
  \vrule\hfil}}}%
  \normalfont \topsep6\p@\@plus6\p@\relax
  {\it \textbf{\proofname { of Theorem {\rm #1}\@addpunct{.}}}\newline\rm}
}{%
   \begin{flushright}\popQED\end{flushright}\@endpefalse
}
\makeatother

\makeatletter

\newenvironment{proofL}[1]{\par
  \pushQED{{\small\it Lemma {\rm #1}} {\leavevmode
  \hbox to.77778em{%
  \hfil\vrule
  \vbox to.675em{\hrule width.6em\vfil\hrule}%
  \vrule\hfil}}}%
  \normalfont \topsep6\p@\@plus6\p@\relax
  {\it \textbf{\proofname { of Lemma {\rm #1}\@addpunct{.}}}\newline\rm}
}{%
   \begin{flushright}\popQED\end{flushright}\@endpefalse
}
\makeatother



\usepackage{array} %\displaystyle%
\newenvironment{ma}{\begin{array}{>{\displaystyle}r >{\displaystyle}c >{\displaystyle}l}}{\end{array}}%

%% file: commands.tex
\newcommand{\aleq}{\prec}
\newcommand{\ageq}{\succ}
\newcommand{\aeq}{\approx}

\newcommand{\subsubset}{\subset\subset}

\newcommand{\N}{{\mathbb N}}
\newcommand{\R}{{\mathbb R}}

\newcommand{\im}{ \mathbf{i}}

\newcommand{\Sw}{{\mathcal{S}}}
\newcommand{\Sws}{{\Sw^{'}}}

\newcommand{\fracm}[1]{{\frac{1}{#1}}}

\newcommand{\supp}{\operatorname{supp}}%
\newcommand{\dist}{\operatorname{dist}}%

\newcommand{\lap}{\Delta}

\newcommand{\laps}[1]{\lap^{\frac{#1}{2}}}
\newcommand{\lapms}[1]{\lap^{-\frac{#1}{2}}}
\newcommand{\lapn}{\lap^{\frac{n}{4}}}
\newcommand{\lapmn}{\lap^{-\frac{n}{4}}}
\newcommand{\abs}[1]{{\left \vert #1 \right \vert}}%
\newcommand{\brac}[1]{{\left ( #1 \right )}}%
\newcommand{\Vrac}[1]{{\left \Vert #1 \right \Vert }}%

\newcommand{\ontop}[2]{{\genfrac{}{}{0pt}{}{#1}{#2}}}

\newcommand{\intl}{\int \limits}

\newcommand{\mvint}{\fint}
\def\XXint#1#2#3{{\setbox0=\hbox{$#1{#2#3}{\int}$}%
     \vcenter{\hbox{$#2#3$}}\kern-.5\wd0}}%

\newcommand{\sref}[2]{#1.\ref{#2}}

%% file: abstract.tex
%abstract
We prove continuity on domains up to the boundary for $n/2$-polyharmonic maps into manifolds. Technically, we show how to adapt H\'elein's direct approach to the fractional setting. This extends a remark by the author that this is possible in the setting of Rivi\`ere's famous regularity result for critical points of conformally invariant variational functionals. Moreover, pointwise behavior for the involved three-commutators is established. Continuity up to the boundary is then obtained via an adaption of Hildebrandt and Kaul's technique to the non-local setting.\\[0.5ex]
%\\[0.5ex]
{\bf Keywords:} Harmonic maps, nonlinear elliptic PDE, regularity of solutions, boundary regularity.\\
{\bf AMS Classification:} 58E20, 35B65, 35J60, 35S05, 35S15.

%% file: intro.tex
\section{Introduction}
In his seminal work \cite{Hel91}, H\'elein proved that harmonic maps from a two-dimensional surface $D$ into a compact manifold $\mathcal{M} \subset \R^N$ are smooth, by an optimal choice of frame $(Pe)_i$ which was obtained by minimizing a simple energy functional of the general form
\[
E(P) := \int_D \abs{P\nabla P^T + P \Omega P^T}^2, \quad \mbox{for $P-I \in W_0^{1,2}(D,\R^{N\times N})$, $P \in SO(N)$ a.e.},
\]
where $\Omega \in L^2(D,so(N)\otimes \R^2)$ is a tensor stemming from the right-hand side of the respective Euler-Lagrange System. In \cite{IchEnergie} the author remarked that this kind of minimizing approach might still be considered helpful in the general setting of Rivi\`ere's celebrated result in \cite{Riv06} where it was shown that in general critical points $u$ of conformally invariant variational functionals between $D$ and $\mathcal{M}$ satisfy an equation like
\[
 \lap u = \Omega \cdot \nabla u,
\]
and are -- because of the antisymmetry of $\Omega$ -- continuous. In fact, instead of constructing a Coulomb gauge adapting the powerful, yet indirect and involved techniques by Uhlenbeck \cite{Uhlenbeck82}, one can still minimize $E(\cdot)$ in order to construct the same gauge, see \cite{IchEnergie} for more details.\\
Nevertheless, there are several settings inspired by Rivi\`ere's result where adaptions of Uhlenbeck's method have seemed more viable in order to show regularity. One of these settings is the work by Da Lio and Rivi\`ere regarding fractional polyharmonic maps, \cite{DR1dMan}, \cite{DndMan} - cf. also \cite{DR1dSphere}, \cite{NHarmS10Arxiv}. Here, we'd like to show how to adapt H\'elein's moving frame approach -- in a similar fashion as in \cite{IchEnergie} -- to the following setting, which can be considered a fairly general model case for these fractional polyharmonic maps $u$ if $v \approx \lapn u$ -- as was shown in \cite{DR1dMan}, \cite{DndMan}:\\
Let $v \in L^2(\R^n)$ be a solution to
\begin{equation}\label{eq:pdelapnvomegav}
 \lapn v = \Omega v \quad \mbox{in $D \subset \R^n$}.
\end{equation}
We then can prove the following theorem, which for $D = \R^n$ was proven first in \cite{DndMan} -- but we will be using H\'elein's direct approach instead of Uhlenbeck's.
\begin{theorem}\label{th:regthm}
Let $v \in L^2(\R^n)$ be a solution of \eqref{eq:pdelapnvomegav}. Then, for any $\tilde{D} \subsubset D$ there exists an $\alpha > 0$, $R > 0$ such that
\[
 \sup_{\ontop{x \in \tilde{D}}{r \in (0,R)}} r^{-\alpha} \Vert v \Vert_{(2,\infty),B_r(x)} < \infty.
\]
In particular, (see \cite{DR1dMan}, \cite{DndMan}) we have $v \in L^p_{loc}(D)$ for any $p \in (1,\infty)$.
\end{theorem}
Moreover, by an extension of techniques by Hildebrandt and Kaul \cite{HildKaul72}, see also \cite{Strzelecki03}, we are able to show that solutions are continuous up to the boundary similar to the two-dimensional case as in \cite{MuellerS09}. More precisely, we have
\begin{theorem}\label{th:regthmGen}
Let $u \in L^2(\R^n)$, $v := \lapn u \in L^2(\R^n)$ be a solution of \eqref{eq:pdelapnvomegav}. Then for some $\alpha \in (0,1)$ we have $u \in C^{0,\alpha}(D)$. Moreover, if $D \subsubset \R^n$, $u \in C^0(\R^n \backslash D)$ and $\partial D \in C^\infty(\R^n)$ we have $u \in C^{0,\alpha}(D) \cap C^0(\R^n)$, in other words $u$ is continuous up to the boundary $\partial D$.
\end{theorem}
Let us sketch the new arguments involved (for necessary definitions we refer to Section \ref{s:basics}): Transforming equation \eqref{eq:pdelapnvomegav} as Da Lio and Rivi\`ere, we have (cf. \eqref{eq:pdew}) for $w := Pv$ and $P-I \in H^{\frac{n}{2}}_0(D)$, $P \in SO(N)$ almost everywhere, for any $\varphi \in C_0^\infty(\R^n)$
\[
\intl_{\R^n} w\ \lapn \varphi = \intl_{\R^n}  so\brac{\Omega_P}\ w\ \varphi +  \intl_{\R^n} \brac{\fracm{2}H(P-I,P^T-I)\ \varphi - H(\varphi,P-I) P^T}\ w.
\]
Here, we denote
\[
  \Omega_Q := Q \lapn (Q^T-I) + Q\Omega Q^T,
\]
\[
  so(A) := \fracm{2} \brac{A-A^T},
\]
and
\[
 H(a,b) := \lapn \brac{ab} - a \lapn b - b \lapn a.
\]
Again, similar in its spirit to \cite{IchEnergie}, instead of using the ingenious adaption of Uhlenbeck's approach by Da Lio and Rivi\`ere, we simply minimize
\[
 E(Q) := \Vert \Omega_Q \Vert_{2,\R^n}^2
\]
on a suitable class of $Q$, cf. Section \ref{s:energy}. Note, that the arguments in \cite{DR1dMan} suggest, that the minimal value should be attained for some $P$ with $E(P) = \Vert \lapn P \Vert_{2}^2$, although we were not able to prove that with this kind of direct method. Instead, we are able to prove that Euler-Lagrange equations of this functional imply that
\begin{equation}\label{eq:soomegapinl21loc}
 so\brac{\Omega_P} \in L^{2,1}_{loc}(D).
\end{equation}
Indeed, in Lemma~\ref{la:en:el} we prove that
\[
 \int so(\Omega_P)\ \lapn \varphi = \int so(H(\varphi,P-I) P^T \Omega_P) \quad \mbox{for all $\varphi \in C_0^\infty(D)$}.
\]
This and the following Lemma, whose localized version will be shown in Lemma~\ref{la:wlapneqhterm}, imply \eqref{eq:soomegapinl21loc}.
\begin{lemma}\label{la:wlapneqhterm:nonlocalized}
Assume that $f,g,h \in L^2(\R^n)$, and that for all $\varphi \in C_0^\infty(B_{10 r})$
\[
 \intl_{\R^n} f \lapn \varphi = \intl_{\R^n} g H(h,\varphi).
\]
Then,
\[
 \Vert f \Vert_{(2,1),B_r} \leq C\ \Vert g \Vert_{2,\R^n}\ \Vert \lapn h \Vert_{2,\R^n} + C\ \Vert f \Vert_{2,\R^n}.
\]
\end{lemma}
In order to show the ``gain in integrability''-effect of Lemma~\ref{la:wlapneqhterm:nonlocalized}, we need some results on the behavior of $H(\cdot,\cdot)$ similar to the one used in \cite{DR1dSphere}, although we prefer to view these, as in \cite{NHarmS10Subm}, in the form of lower order operators:\\
In \cite{NHarmS10Subm} the author remarked that by a fairly simple argument inspired by Tartar's approach to Wente's inequality \cite{Tartar85}, quantities like $H(\cdot,\cdot)$ behave like a product of lower order operators -- after taking the Fourier transform. As we deal here with spaces different from $L^2$, Tartar's argument (which for our purposes relies on Plancherel's theorem) does not apply that easily any more in order to get our needed estimates. One might try bilinear real interpolation on the fractional ``Leibniz rule'' originally due to Kato and Ponce \cite{KP88}, see also \cite{Hofmann98}. Another possibility is the following, and it is closer to the argument in \cite{NHarmS10Subm}: Using simple estimates on multipliers appearing in the representation as \emph{potential} of the involved operators rather than their representation as Fourier multiplier, one can be quite specific (even pointwise) about how $H(\cdot,\cdot)$ behaves like a product of lower order operators\footnote{For the sake of shortness of presentation, we will restrict the proof to cases where $n \geq 5$ and $n-1 \in 4\N$.}\footnote{It seems likely, that using the general potential representation of $\laps{s}$ for arbitrary $s \in \R$, cf. \cite{SKM93}, by arguments similar to the ones we use here, one might obtain a more precise estimate. Nevertheless, this is not needed for our argument.}:
\begin{lemma}\label{la:commlowerorder}
For some constants $L \in \N$, $s_k \in (0,\frac{n}{2})$, $t_k \in (0,s_k]$, $C > 0$, for zero-multiplier operators $M_{k,1}$, $M_{k,2}$, $M_{k,3}$, and for any $a$, $b \in \Sw(\R^n)$
\[
 \abs{H(\lapmn a,\lapmn b)}(x) \leq C\ \sum_{k=1}^L M_{k,1}\lapms{s_k - t_k} \brac{M_{k,2}\lapms{t_k} \abs{a}\ M_{k,3}\lap^{-\frac{n}{4}+\frac{s_k}{2}} \abs{b}}(x).
\]
\end{lemma}
With this, instead of dealing with paraproducts (although, of course, the underlying arguments are similar), Sobolev's inequality shows all the necessary ``integrability gain'' or ``compensation phenomena'' to be used (see Proposition~\ref{pr:lpestloop}).\\
Then, an argument similar to the one in \cite{DR1dMan} (though locally in $D$ instead of $\R^n$), implies the following Lemma, of which Theorem~\ref{th:regthm} is a consequence by an iteration result as in Lemma~\ref{la:it}
\begin{lemma}\label{la:localcontroldPreIt}
Let $v \in L^2(\R^n)$ be a solution of \eqref{eq:pdelapnvomegav}. Then there exists $\Lambda_0 > 0$, $\gamma > 0$, $C \equiv C_{v,\Omega} > 0$ such that for any $\Lambda > \Lambda_0$ there is an $R \in (0,1)$  such that if $B_{\Lambda r}(x) \subset D$, $r \in (0,R)$
\[
\Vert v \Vert_{(2,\infty),B_r(x)} \leq C\ \Lambda^{-\gamma}\ \Vert v \Vert_{(2,\infty),B_{\Lambda r}(x)} + C\ \Lambda^{-\gamma}\ \sum_{k=1}^\infty 2^{-\gamma k} \Vert v \Vert_{(2,\infty),B_{2^{k} \Lambda r}(x)\backslash B_{2^{k-1} \Lambda r}(x)}.
\]
\end{lemma}
${}$\\[2ex]
We will use fairly standard \emph{notation}, similar to \cite{NHarmS10Subm}:\\
As usual, we denote by $\Sw \equiv \Sw(\R^n)$ the Schwartz class of all smooth functions which at infinity tend faster to zero than any quotient of polynomials, and by $\Sws \equiv \Sws(\R^n)$ its dual. We say that $A \subsubset \R^n$ if $A$ is a bounded subset of $\R^n$. For a set $A \subset \R^n$ we will denote its $n$-dimensional Lebesgue measure by $\abs{A}$, and $rA$, $r > 0$, will be the set of all points $rx \in \R^n$ where $x \in A$. By $B_r(x) \subset \R^n$ we denote the open ball with radius $r$ and center $x \in \R^n$. If no confusion arises, we will abbreviate $B_r \equiv B_r(x)$. For a real number $p \geq 0$ we denote by $\lfloor p \rfloor$ the biggest integer below $p$ and by $\lceil p \rceil$ the smallest integer above $p$. If $p \in [1,\infty]$ we usually will denote by $p'$ the H\"older conjugate, that is $\frac{1}{p} + \frac{1}{p'} = 1$. By $f \ast g$ we denote the convolution of two functions $f$ and $g$. We set $f^\wedge$ to be the Fourier transform and $f^\vee$ to be the inverse Fourier transform, which on the Schwartz class $\Sw$ shall be defined as
\[
 f^\wedge(\xi) := \intl_{\R^n} f(x)\ e^{-2\pi \im\ x\cdot \xi}\ dx,\quad f^\vee(x) := \intl_{\R^n} f(\xi)\ e^{2\pi \im\ \xi\cdot x}\ d\xi.
\]
By $\im$ we denote here and henceforth the imaginary unit $\im^2 = -1$. We will speak of a zero-multiplier operator $M$, if there is a function $m \in C^\infty(\R^n \backslash \{0\})$ homogeneous of order $0$ and such that $(Mv)^\wedge(\xi) = m(\xi)\ v^\wedge(\xi)$ for all $\xi \in \R^n \backslash \{0\}$. For a measurable set $D \subset \R^n$, we denote the integral mean of an integrable function $v: D \to \R$ to be $(v)_D \equiv \mvint_D v \equiv \frac{1}{\abs{D}} \int_D v$. Lastly, our constants -- frequently denoted by $C$ or $c$ --  can possibly change from line to line and usually depend on the space dimensions involved, further dependencies will be denoted by a subscript, though we will make no effort to pin down the exact value of those constants. If we consider the constant factors to be irrelevant with respect to the mathematical argument, for the sake of simplicity we will omit them in the calculations, writing $\aleq{}$, $\ageq{}$, $\aeq{}$ instead of $\leq$, $\geq$ and $=$.\\
We will use the same cutoff-functions as in, e.g., \cite{DR1dSphere}, \cite{NHarmS10Subm}: $\eta^k_{r} \in C_0^\infty(A_{r,k})$ where
\[
 B_{r,k}(x) := B_{2^kr}(x)
\]
for $k \geq 1$,
\[
 A_{r,k}(x) := B_{r,k+1}(x) \backslash B_{r,k-1}(x),
\]
and for $k = 0$
\[
 A_{r,0}(x) := B_{r,0}(x).
\]
Moreover, $\sum_k \eta^k_r \equiv 1$ pointwise everywhere, and we assume that $\abs{\nabla^l \eta^k_r} \leq C_l\  \brac{2^k r}^{-l}$.

\vspace{2ex}\noindent
{\bf Acknowledgment.} The author would like to thank Francesca Da Lio and Tristan Rivi\`{e}re for introducing him to the topic.\\
Parts of this work were conducted while the author was a guest at FIM at ETH Z\"urich which was supported by FIM and DAAD grant D/10/50763.

%% file: facts.tex
\section{Preliminaries}\label{s:basics}
\subsection{Lorentz Spaces}
In this section, we recall the definition of Lorentz spaces, which are a refinement of the standard $L^p$-spaces. For more on Lorentz spaces, the interested reader might consider \cite{Hunt66}, \cite{Ziemer89}, \cite[Section~1.4]{GrafC08}, and also \cite{Tartar07}.
\begin{definition}[Lorentz Space]
Let $f: \R^n \to \R$ be a Lebesgue-measurable function. We denote
\[
d_f(\lambda) := \abs{\{ x \in \R^n\ :\ \abs{f(x)} > \lambda\}}.
\]
The decreasing rearrangement of $f$ is the function $f^\ast$ defined on $[0,\infty)$ by
\[
f^\ast(t) := \inf \{s > 0:\ d_f(s) \leq t\}.
\]
For $1 \leq p \leq \infty$, $1 \leq q \leq \infty$, the Lorentz space $L^{p,q} \equiv L^{p,q}(\R^n)$, is the set of measurable functions $f: \R^n \to \R$ such that $\Vert f \Vert_{L^{p,q}} < \infty$, where
\[
\Vert f \Vert_{(p,q),\R^n} \equiv \Vert f \Vert_{(p,q)} \equiv \Vert f \Vert_{L^{p,q}(\R^n)} := \begin{cases}
													\left (\intl_0^\infty \left (t^{\frac{1}{p}} f^\ast (t) \right )^q \frac{dt}{t} \right )^\frac{1}{q}, \quad &\mbox{if $q < \infty$,}\\
													\sup_{t > 0} t^{\frac{1}{p}} f^\ast (t),\quad &\mbox{if $q = \infty$, $p < \infty$},\\
													\Vert f \Vert_{L^\infty(\R^n)},\quad &\mbox{if $q = \infty$, $p = \infty$}.
														\end{cases}
\]
As usual, if $A \subset \R^n$ and $\chi_A$ denotes its characteristic function, we define
\[
 \Vert f \Vert_{(p,q),A} := \Vert \chi_{A} f \Vert_{(p,q),\R^n}.
\]
If $p = q$, and as a consequence $L^{p,q} = L^p$, we write instead of $\Vert \cdot \Vert_{(p,p),A}$ mostly $\Vert \cdot \Vert_{p,A}$.
\end{definition}

\begin{remark}
Observe that $\Vert \cdot \Vert_{\brac{p,q}}$ as defined here does \emph{not} satisfy the triangle inequality. Nevertheless, if one defines $\Vert \cdot \Vert_{\brac{p,q}}$ replacing $f^\ast$ by the averaged version $f^{\ast\ast}$,
\[
 f^{\ast\ast}(t) := \fracm{t} \int_0^t f^\ast(s)\ ds,
\]
one obtains an equivalent quantity for $1 < p < \infty$, $q \in [1,\infty]$, and this is in fact a norm, see \cite{Hunt66} or \cite[Ex. 1.4.3]{GrafC08}. We will switch between the two definitions without mentioning it again.
\end{remark}

\begin{proposition}[Some facts about Lorentz spaces]\label{pr:lorentzfacts}
\begin{itemize}
 \item[(i)] $\brac{L^{p,q}}^\ast = L^{p',q'}$ for $1 < p,q < \infty$.
 \item[(ii)] $\brac{L^{p,1}}^\ast = L^{p',\infty}$ for $1 < p < \infty$.
\item[(iii)] Simple functions are dense in $L^{p,q}$ for $p \in (1,\infty)$, $q \in [1,\infty)$.
\item[(iv)] Simple functions are not dense in $L^{p,\infty}$, $p \in [1,\infty]$.
\item[(v)] $C_0^\infty(\R^n)$ is dense in $L^{p,q}(\R^n)$ and for smoothly bounded $A \subsubset \R^n$ also $C_0^\infty(A)$ is dense in $L^{p,q}(A)$ for $p \in (1,\infty)$, $q \in [1,\infty)$.
 \item[(vi)] For $p \in (1,\infty)$, $q \in [1,\infty]$, $A = \R^n$ or $A$ a smoothly bounded domain
\[
 \Vert f \Vert_{(p,q),A} \aeq \sup_{\ontop{g \in C_0^\infty(A)}{\Vert g \Vert_{(p',q'),A} \leq 1}} \int fg.
\]
\end{itemize}
\end{proposition}
\begin{proofP}{\ref{pr:lorentzfacts}}
Facts (i), (ii) are stated and proven in \cite[Theorem 1.4.17]{GrafC08}. Claim (iii) can be found in \cite[Theorem 1.4.13]{GrafC08}. And claim (iv) is explained in \cite[Ex. 1.4.4]{GrafC08}. Claim (v) and (vi) then can be proven by an approximation scheme.
\end{proofP}

\begin{remark}\label{rem:nowidman}
Note, however, that there is no reason for $q \neq p$ that
\[
\Vert f \Vert_{(p,q),A}^p + \Vert f \Vert_{(p,q),B}^p = \Vert f \Vert_{(p,q),A \cup B}^p
\]
even if $A$ and $B$ are disjoint sets.
\end{remark}

\subsection{Some Facts about the Fractional Laplacian and its Inverse}
\begin{definition}
(cf. \cite[Chapter 5, \textsection 26]{SKM93})\\
For $f \in \Sw(\R^n)$, $s \geq 0$, we define the operator $\laps{s} f$ via
\[
 \brac{\laps{s} f}^\wedge(\xi) := \brac{\omega_{n}}^s \abs{\xi}^s f^\wedge(\xi),
\]
where $\omega_{n} \in \R$ is chosen such that the classic differential operator $\lap = \sum_i \partial_i\partial_i$ suffices
\[
 \lap f \equiv \laps{2} f \quad \mbox{for any $f \in \Sw(\R^n)$}.
\]
For $s \in (0,2)$ we set for a function $f \in L^1(\R^n) + L^\infty(\R^n)$
\[
  \laps{s} f (x):= c_s\ \mbox{P.V.-}\intl_{\R^n} \frac{f(x-y)+f(x+y)-2f(x)}{\abs{y}^{n+s}}\ dy,
\]
whenever this integral is defined. One can show (see, e.g., \cite{NHarmS10Arxiv}), that these operators coincide on $\Sw(\R^n)$ if $c_s$ is chosen appropriately.\\
Finally, for $s \in (0,-n)$ we set for $f \in L^1(\R^n) + L^\infty(\R^n)$
\[
 \lapms{s} f(x) := c_s \mbox{P.V.-}\intl_{\R^n} \abs{x-y}^{-n+\abs{s}}\ f(y),
\]
whenever this integral is well-defined. One checks, that this is the case if $f \in L^p(\R^n)$ for any $p \in [1,\frac{n}{s})$. If $f \in \Sw(\R^n)$ then for any $\varphi \in C_0^\infty(\R^n)$
\[
 \intl_{\R^n} \lapms{s} f\ \varphi = \intl_{\R^n} \abs{\cdot}^{-s} f^\wedge(\cdot)\ \varphi^\vee(\cdot),
\]
that is $\brac{\lapms{s} f}^\wedge = c \abs{\cdot}^s f^\wedge$ in the sense of distributions, and by the same argument $\laps{s} \lapms{s} f = f = \lapms{s} \laps{s} f$ in the sense of distributions for any $f \in \Sw(\R^n)$ and this is true also pointwise, see for example \cite[Chapter 5, \textsection 26.3]{SKM93}.
\end{definition}

\begin{definition}[Fractional Sobolev Spaces]\label{def:Hs}
For $s \geq 0$ we set
\[
 H^s \equiv H^{s}(\R^n) := \left \{ f \in L^2(\R^n), \quad \laps{s} f \in L^2(\R^n)\right \}.
\]
Moreover, for a domain $D \subset \R^n$ we denote
\[
 H^s_0(D) := \left \{ f \in H^s(\R^n), \quad \supp f \subset \overline{D}\right \},
\]
and more generally for measurable $\varphi: \R^n \to \R^n$ we define
\[
  H^s_\varphi(D) := \left \{ f \in L^1(\R^n) + L^\infty(\R^n), f-\varphi \in H^s_0(D) \right \}.
\]
As usual for any finitely dimensional vectorspace $V$ and any subset $A \subset V$ we mean by $H^s(\R^n,A)$ all these vectorvalued functions $f: \R^n \to V$ such that $f \in A$ almost everywhere and $\langle f, v \rangle \in H^s(\R^n)$. Similar definitions are used for $H^s_0(D,A)$ and $H^s_\varphi(D,A)$.
\end{definition}

\begin{proposition}\label{pr:finswlapsinlp}
Let $f \in \Sw(\R^n)$ then $\laps{s} f \in L^{p,q}(\R^n)$ for any $s \geq 0$ and arbitrary $p \in (1,\infty)$, $q \in [1,\infty]$, as well as $(p,q) = (1,1)$ and $(p,q) = (\infty,\infty)$. In particular,
\[
\Vert \laps{s} \eta_r^k \Vert_{(p,q),\R^n} \leq C_{p,q}\ \brac{2^k r}^{\frac{n}{p}-s}.
\]
\end{proposition}
\begin{proofP}{\ref{pr:finswlapsinlp}}
It suffices to show the claim for $s \in [0,2)$, as $\lap^k f \in \Sw(\R^n)$ for any $k \in \N$. The operator $\laps{s} f$ is defined by the Fourier-definition, so for any $p \in [2,\infty]$
\[
 \Vert \laps{s} f \Vert_{p,\R^n} \aleq{} \Vert \abs{\cdot}^s\ f^\wedge \Vert_{p',\R^n} \leq C,
\]
because $f^\wedge \in \Sw(\R^n)$. By interpolation arguments, it now suffices to show that the claim holds also for $p = 1$. In this case,
\[
 \intl_{\abs{y} \geq 2} \fracm{\abs{y}^{n+s}} \intl_{\R^n} \abs{f(x+y)}+\abs{f(x-y)}+2\abs{f(x)}\ dx\ dy \aleq{} \Vert f \Vert_{1,\R^n}.
\]
And
\[
\begin{ma}
 &&\intl_{\abs{y} \leq 2} \fracm{\abs{y}^{n+s}} \intl_{\R^n} \abs{f(x+y) + f(x-y) - 2f(x)}\ dx\ dy\\
&\aleq& \intl_{0}^1 \intl_0^1 \intl_{\abs{y} \leq 2} \fracm{\abs{y}^{n+s-2}} \intl_{\R^n} \abs{\nabla^2 f(x+sty)}\ dx\ dy\ dt\ ds\\
&\aleq&  \Vert \nabla^2 f \Vert_{1,\R^n}.
\end{ma}
\]
\end{proofP}

% \begin{lemma}
% We have for $s \in (0,n)$, $f \in \Sw(\R^n)$,
% \[
%  g := \lapms{s} f \in L^{p}(\R^n), \quad \mbox{for any $p > \frac{n}{n-s}$},
% \]
% and $\laps{s} g$ is well defined as above and
% \[
%  \laps{s} g \equiv f \quad \mbox{almost everywhere}.
% \]
% \end{lemma}
% \begin{proof}
% \cite{Samko98} \ToDo
% \end{proof}
%
% \begin{lemma}
% Let $f \in L^1(\R^n) + L^\infty(\R^n)$ and $h \in C^\infty(\R^n) \cap L^\infty(\R^n)$, then if $\laps{s} f \in L^1(\R^n) + L^\infty(\R^n)$ for some $s$, also $\laps{s} (hf) \in L^1(\R^n) + L^\infty(\R^n)$.
% \end{lemma}
% \begin{proof}
% \ToDo.
% \end{proof}
%
%

\subsection{Sobolev And Poincar\'e Inequalities}
As we stated in the introduction, once we show that certain quantities behave like products of lower order operators, all we need are Sobolev inequality, some versions of which we are going to state in this section.
\begin{proposition}[Sobolev inequality]\label{pr:sobineq}
For $s \in (0,n)$, $p_1 \in (1,\frac{n}{s})$, $p \in (\frac{n}{n-s},\infty)$ such that
\[
 \frac{1}{p_1} - \frac{s}{n} = \frac{1}{p},
\]
and for any zero-multiplier operator $M$ there is a constant $C_{M,p,s}$ such that for any $q \in [1,\infty]$
\[
\Vert M\lapms{s} a \Vert_{(p,q),\R^n} \leq C_{M,p,s}\ \Vert a \Vert_{(p_1,q),\R^n}.
\]
\end{proposition}

\begin{lemma}[Convolution, one limit case]\label{la:lorentzinftyconv}
(\cite[Lemma 4.7, p. 25]{Hunt66})\\
For any $p \in (1,\infty)$, $q \in [1,\infty]$ we have a constant $C_p > 0$ such that
\[
 \Vert f \ast g \Vert_{\infty} \leq C_p\ \Vert f \Vert_{(p,q)}\ \Vert g \Vert_{(p',q')}.
\]
\end{lemma}
\begin{proofL}{\ref{la:lorentzinftyconv}}
We have for any $x \in \R^n$ (see, e.g., \cite[Exercise 1.4.1(b)]{GrafC08})
\[
 \intl_{\R^n} f(y)\ g(x-y)\ dy \leq \intl_0^\infty f^\ast(t)\ \brac{g(x-\cdot)}^\ast(t)\ dt.
\]
As
\[
 \brac{g(x-\cdot)}^\ast(t) = g^\ast(t),
\]
this implies
\[
\begin{ma}
 \intl_{\R^n} f(y)\ g(x-y)\ dy &\leq& \intl_0^\infty f^\ast(t)\ g^\ast(t)\ dt\\
&=& \intl_0^\infty t^{\fracm{p}-\fracm{q}}f^\ast(t)\ t^{\fracm{p'}-\fracm{q'}}g^\ast(t)\ dt\\
&\aeq& \Vert f \Vert_{(p,q),\R^n}\ \Vert g \Vert_{(p',q'),\R^n}.
\end{ma}
\]
\end{proofL}

\begin{lemma}[Sobolev inequality for $\infty$]\label{la:soblinfty}
For any $s \in (0,n)$ and any zero-multiplier operator $M$ there is a constant $C_{M,s} > 0$ such that for any $g \in \Sw(\R^n)$,
\[
 \Vert M\lapms{s} g \Vert_{\infty,\R^n} \leq\ C_{M,s} \Vert g \Vert_{(\frac{n}{s},1),\R^n}
\]
\end{lemma}
\begin{proofL}{\ref{la:soblinfty}}
Consequence of Lemma~\ref{la:lorentzinftyconv} and of the fact, that $m(\cdot)\abs{\cdot}^{-n+s} \in L^{\frac{n}{n-s},\infty}(\R^n)$, where $m(\cdot)$ is the zero-homogeneous symbol of $M$.
\end{proofL}

\begin{proposition}[Lorentz Space Estimates for Products of Lower Order Operators]\label{pr:lpestloop}
For any $\beta \in (0,\frac{n}{2})$ and $\alpha \in (0,\beta]$ and zero-multiplier-operators $M_1$, $M_2$, $M_3$ the following holds: For any $p \in \brac{\frac{n}{n-(\beta-\alpha)},\frac{n}{\alpha}}$ (and in particular for $p = 2$) and any $q,\ q_1,\ q_2 \in [1,\infty]$ such that
\[
 \fracm{q} = \fracm{q_1}+\fracm{q_2}
\]
there exists a constant $C_{\alpha,\beta,p,q,q_1,M}$ such that for any $a,b \in \Sw(\R^n)$
\[
 \Vert M_1 \lapms{\beta - \alpha} \brac{M_2 \lapms{\alpha} \abs{a}\ M_3 \lap^{-\frac{n}{4}+\frac{\beta}{2}} \abs{b}} \Vert_{(p,q),\R^n} \leq C_{\alpha,\beta,p,q,q_1}\ \Vert a \Vert_{(p,q_1),\R^n}\ \Vert b \Vert_{(2,q_2),\R^n}.
\]
\end{proposition}
\begin{remark}
Note that for, say, $\alpha = \frac{n}{2}$ and $\beta = 0$ (or vice versa), the estimate
\[
 \Vert \brac{\lapmn a}\ b \Vert_{(2,q),\R^n} \aleq \Vert a \Vert_{(2,q_1),\R^n}\ \Vert b \Vert_{(2,q_2),\R^n}.
\]
does not need to hold.
\end{remark}
\begin{proofP}{\ref{pr:lpestloop}}
First of all, by Proposition~\ref{pr:sobineq} we have
\[
 \Vert M_1 \lapms{\beta - \alpha} \brac{M_2 \lapms{\alpha} \abs{a}\ M_3 \lap^{-\frac{n}{4}+\frac{\beta}{2}} \abs{b}} \Vert_{(p,q),\R^n} \aleq{}
\Vert M_2\lapms{\alpha} \abs{a}\ M_3 \lap^{-\frac{n}{4}+\frac{\beta}{2}} \abs{b} \Vert_{(p_1,q)},
\]
for
\[
  \frac{1}{p_1} - \frac{\beta-\alpha}{n} = \frac{1}{p},
\]
if $p \in (\frac{n}{n-(\beta-\alpha)},\infty)$ (of course, this holds also if $\beta-\alpha = 0$). Moreover,
\[
 \Vert M_2 \lapms{\alpha} \abs{a} \Vert_{(p_2,q_1),\R^n} \aleq{} \Vert a \Vert_{(p,q_1),\R^n},
\]
for
\[
  \frac{1}{p} - \frac{\alpha}{n} = \frac{1}{p_2},
\]
if $p \in (1,\frac{n}{\alpha})$. Next we define $p_3 \in (2,\infty)$ via
\[
 \fracm{p_3} = \fracm{p_1} - \fracm{p_2} =  \frac{1}{p} +  \frac{\beta-\alpha}{n} - \frac{1}{p} + \frac{\alpha}{n} = \frac{\beta}{n}.
\]
Then
\[
 \Vert M_3 \lap^{-\frac{n}{4}+\frac{\beta}{2}} \abs{b} \Vert_{(p_3,q_2)} \aleq{} \Vert b \Vert_{(p_4,q_2)},
\]
(the condition $p_3 \in (\frac{n}{\frac{n}{2} + \beta},\infty)$ is satisfied because $p_3 > 2$) for $p_4 \in (1,\infty)$ such that
\[
 \fracm{p_3} = \fracm{p_4} - \frac{\frac{n}{2}-\beta}{n},
\]
i.e.
\[
 \fracm{p_4} = \fracm{p_3} + \frac{\frac{n}{2}-\beta}{n} = \frac{1}{2},
\]
that is, $p_4 = 2$. Together we have
\[
 \begin{ma}
&& \Vert M_1 \lapms{\beta - \alpha} M_2 \brac{\lapms{\alpha} \abs{a}\ M_3\lap^{-\frac{n}{4}+\frac{\beta}{2}} \abs{b}} \Vert_{(p,q),\R^n}\\
&\aleq{}& \Vert M_2 \lapms{\alpha} \abs{a}\ M_3 \lap^{-\frac{n}{4}+\frac{\beta}{2}} \abs{b} \Vert_{(p_1,q)} \\
&\aleq{}& \Vert M_1 \lapms{\alpha} \abs{a} \Vert_{(p_2,q_1)}\ \Vert M_3 \lap^{-\frac{n}{4}+\frac{\beta}{2}} \abs{b} \Vert_{(p_3,q_2)} \\
&\aleq{}& \Vert a \Vert_{(p,q_1)}\ \Vert b \Vert_{(2,q_2)}.
\end{ma}
\]
\end{proofP}

\begin{lemma}[Localized Sobolev Inequality]\label{la:locsobineq}
Let $p_1,p_2 \in (1,\infty)$, $q \in [1,\infty]$ and $s \in [0,n)$ so that
\[
 \fracm{p_2} = \fracm{p_1}+\frac{s}{n}.
\]
For any zero-multiplier operator $M$ there is a constant $C_{p_1,q,s,M} > 0$ and $\gamma \equiv \gamma_{p_1,s} > 0$ such that for any $a \in \Sw(\R^n)$, $\Lambda \geq 1$ and $B_r \subset \R^n$
\[
 \Vert M\lapms{s} a \Vert_{(p_1,q),B_r} \leq C_{p_1,s} \brac{\Vert a \Vert_{(p_2,q),B_{2\Lambda r}} +  \Lambda^{-\gamma} \Vert a \Vert_{(p_2,q),\R^n}}.
\]
\end{lemma}
\begin{proofL}{\ref{la:locsobineq}}
W.l.o.g. we assume $\Lambda \geq 4$; For all smaller $\Lambda$, the claim is just Sobolev's inequality. We have
\[
 \abs{\eta_r M\lapms{s} a} \leq \abs{M\lapms{s} \brac{\eta_{\Lambda r} a}} + \abs{\eta_r M\lapms{s} \brac{(1-\eta_{\Lambda r}) a}},
\]
so
\[
  \Vert M\lapms{s} a \Vert_{(p_1,q),B_r} \aleq{} \Vert \eta_{\Lambda r} a \Vert_{(p_2,q),\R^n} + \Vert \eta_r M\lapms{s} \brac{ \brac{1-\eta_{\Lambda r}} a} \Vert_{(p_1,q),\R^n}.
\]
It remains to estimate the second term, which can be done by duality, cf. Proposition~\ref{pr:lorentzfacts}: For some $\psi \in C_0^\infty(B_{2r})$, $\Vert \psi \Vert_{(p_1',q'),\R^n} \leq 1$,
\[
\begin{ma}
 &&\Vert \eta_r M\lapms{s} (1-\eta_{\Lambda r}) a \Vert_{p_1,B_r}\\
&\aleq& \int \psi\ M\lapms{s} (1-\eta_{\Lambda r}) a\\
&\aleq& \int \abs{\xi}^{-n+s}\chi_{\abs{\xi} \geq \frac{\Lambda}{2}r}\ \brac{\abs{\psi} \ast \abs{a}}(\xi)\ d\xi.
\end{ma}
\]
Now let $p_3 \in (1,p_1')$ be close enough to $p_1'$ such that still
\[
 \frac{1}{p_4} := \fracm{p_3} - \fracm{p_1'}+\frac{s}{n} \in \brac{\frac{s}{n},1}.
\]
Then, for the $\psi$ from above
\[
\begin{ma}
 && \Vert \eta_r M\lapms{s} (1-\eta_{\Lambda r}) a \Vert_{(p_1,q),B_r}\\
&\aleq& \Vert \abs{\cdot}^{-n+s}\chi_{\abs{\cdot} \geq \frac{\Lambda}{2}r} \Vert_{p_4',\R^n}\ \Vert \abs{\psi} \ast \abs{a} \Vert_{p_4,\R^n}\\
&\aleq& \Vert \abs{\cdot}^{-n+s}\chi_{\abs{\cdot} \geq \frac{\Lambda}{2}r} \Vert_{p_4',\R^n}\ \Vert \psi \Vert_{(p_3,q'),\R^n}\ \Vert a \Vert_{(p_2,q),\R^n}\\
&\aleq& \Vert \abs{\cdot}^{-n+s}\chi_{\abs{\cdot} \geq \frac{\Lambda}{2}r} \Vert_{p_4',\R^n}\ r^{\frac{n}{p_3}-\frac{n}{p_1'}}\ \Vert \psi \Vert_{(p_1',q),\R^n}\ \Vert a \Vert_{(p_2,q),\R^n}\\
&\leq& \Vert \abs{\cdot}^{-n+s}\chi_{\abs{\cdot} \geq \frac{\Lambda}{2}r} \Vert_{p_4',\R^n}\ r^{\frac{n}{p_3}-\frac{n}{p_1'}}\ \Vert a \Vert_{(p_2,q),\R^n}.
\end{ma}
\]
Because $p_4 < \frac{n}{s}$, we have $\fracm{p_4'} < 1 - \frac{s}{n}$ and thus $(-n+s)p_4' < -n$. Hence,
\[
\begin{ma}
 \Vert \abs{\cdot}^{-n+s}\chi_{\abs{\cdot} \geq \frac{\Lambda}{2}r} \Vert_{p_4',\R^n} &\aeq& \brac{\intl_{\frac{\Lambda}{2} r}^\infty t^{(-n+s)p_4'+n-1}\ dt }^{\fracm{p_4'}}\\
&\aeq& \brac{\brac{\Lambda r}^{(-n+s)p_4'+n} }^{\fracm{p_4'}}\\
&=& \brac{\Lambda r}^{(-n+s)+n-\frac{n}{p_4}}\\
&=& \brac{\Lambda r}^{s-\frac{n}{p_4}}\\
&=& \brac{\Lambda r}^{\frac{n}{p_1'} - \frac{n}{p_3}}.
\end{ma}
\]
We conclude that
\[
\Vert \eta_r \lapms{s} (1-\eta_{\Lambda r}) a \Vert_{(p_1,q),\R^n} \aleq \Lambda^{\frac{n}{p_1'} - \frac{n}{p_3}}\ \Vert a \Vert_{(p_2,q),\R^n}.
\]
\end{proofL}

The following lemma is the Poincar\'e inequality. There are several ways to prove it, but in the general setting of Lorentz spaces, we preferred the following.
\begin{lemma}[Poincar\'{e} Inequality]\label{la:poinc}
For any $B_r \subset \R^n$, $p \in (1,\infty)$, $q \in [1,\infty]$ there exists a constant $C_{p,q} > 0$ such that
\begin{equation}\label{eq:poinc}
\Vert f \Vert_{(p,q),\R^n} \leq C_{s,p,q}\ r^{\frac{n}{2}}\ \Vert \lapn f \Vert_{(p,\infty),\R^n}, \quad \mbox{for all $f \in C_0^\infty(B_r)$}.
\end{equation}
\end{lemma}
\begin{proofL}{\ref{la:poinc}}
Let
\[
p_1 :=  \begin{cases}
  p \quad & \mbox{if $p > 2$,}\\
  4 \quad & \mbox{if $p = 2$,}\\
  \frac{2p}{2-p} \quad & \mbox{if $p < 2$},
 \end{cases}
\]
\[
q_1 :=  \begin{cases}
  q \quad & \mbox{if $p \geq 2$,}\\
  \infty \quad & \mbox{if $p < 2$},
 \end{cases}
\]
and
\[
p_2 :=  \begin{cases}
  \frac{2p}{2+p} \quad & \mbox{if $p > 2$,}\\
  \frac{4}{3} \quad & \mbox{if $p = 2$,}\\
  p \quad & \mbox{if $p < 2$}.
 \end{cases}
\]
Then,
\[
 \begin{ma}
 \Vert f  \Vert_{(p,q),\R^n} &=& \Vert f  \Vert_{(p,q),B_r} \\
&\aleq& r^{\frac{n}{p}-\frac{n}{p_1}}\ \Vert \lapmn \lapn f  \Vert_{(p_1,q_1),B_r}\\
&\aleq& r^{\frac{n}{p}-\frac{n}{p_1}}\ \Vert \lapn f  \Vert_{(p_2,q_1),\R^n}
\end{ma}
\]
If $p_2 = p$, i.e. if $p < 2$, this proves the claim. If $p \geq 2$, so $q_1 = q$, let
\[
 p_3 \in (1,p_2),
\]
and $p_4 \in (1,\infty)$ such that
\[
 1-\fracm{p_4} = \fracm{p_3} - \fracm{p_2}.
\]
Then, for some $\Lambda > 4$ to be chosen later,
\[
\begin{ma}
\Vert f  \Vert_{(p,q),\R^n}
&\aleq& r^{\frac{n}{p}-\frac{n}{p_1}}\ \Vert \eta_{\Lambda r} \lapn f  \Vert_{(p_2,q),\R^n} + r^{\frac{n}{p}-\frac{n}{p_1}}\ \Vert \brac{1-\eta_{\Lambda r}} \lapn f  \Vert_{(p_2,q),\R^n}\\

&\aleq& r^{\frac{n}{p}-\frac{n}{p_1}}\ \brac{\Lambda r}^{\frac{n}{p_2}-\frac{n}{p}}\ \Vert \lapn f  \Vert_{(p,\infty),\R^n} + r^{\frac{n}{p}-\frac{n}{p_1}}\ \Vert \brac{1-\eta_{\Lambda r}} \lapn f  \Vert_{(p_2,q),\R^n}\\

&=& \Lambda^{\frac{n}{p_2}-\frac{n}{p}}\ r^{\frac{n}{2}}\Vert \lapn f  \Vert_{(p,q),\R^n} + r^{\frac{n}{p}-\frac{n}{p_1}}\ \Vert \brac{1-\eta_{\Lambda r}} \lapn f  \Vert_{(p_2,q),\R^n}.
\end{ma}
\]
For some $\psi \in C_0^\infty(\R^n \backslash B_{\Lambda r})$, $\Vert \psi \Vert_{(p_2',q')} \leq 1$ we have
\[
\begin{ma}
  \Vert \brac{1-\eta_{\Lambda r}} \lapn f  \Vert_{(p_2,q),\R^n} &\aleq& \int \abs{\cdot}^{-\frac{3}{2}n} \chi_{\abs{\cdot}> \frac{\Lambda}{2}r}\ \abs{f}\ast \abs{\psi}\\
&\aleq& \brac{\Lambda r}^{-\frac{3}{2}n + \frac{n}{p_4}}\ \Vert f \Vert_{\brac{p_3,q},B_r} \Vert \psi \Vert_{(p_2',q')}\\
&\aleq& \brac{\Lambda r}^{-\frac{3}{2}n + \frac{n}{p_4}}\ r^{\frac{n}{p_3}-\frac{n}{p}}\ \Vert f \Vert_{\brac{p,q},B_r}.
\end{ma}
\]
Hence,
\[
 \begin{ma}
 \Vert f  \Vert_{(p,q),\R^n}
&\aleq& \Lambda^{\frac{n}{p_2}-\frac{n}{p}}\ r^{\frac{n}{2}}\ \Vert \lapn f  \Vert_{(p,\infty),\R^n} +
\Lambda^{-\frac{3}{2}n + \frac{n}{p_4}}\ r^{\frac{n}{p}-\frac{n}{p_1}-\frac{3}{2}n + \frac{n}{p_4}+\frac{n}{p_3}-\frac{n}{p}}\ \Vert f \Vert_{\brac{p,q},B_r}\\
&\aeq& \Lambda^{\frac{n}{p_2}-\frac{n}{p}}\ r^{\frac{n}{2}}\ \Vert \lapn f  \Vert_{(p,\infty),\R^n} +
\Lambda^{-\frac{3}{2}n + \frac{n}{p_4}} \Vert f \Vert_{\brac{p,q},\R^n}.
\end{ma}
\]
As $-\frac{3}{2}n + \frac{n}{p_4} < 0$ we can pick some $\Lambda > 1$ large enough, so that
\[
 \Vert f  \Vert_{(p,q),\R^n} \leq C_\Lambda\ r^{\frac{n}{2}}\ \Vert \lapn f  \Vert_{(p,\infty),\R^n} +
\fracm{2} \Vert f \Vert_{\brac{p,q},\R^n}
\]
which implies the claim.
\end{proofL}

\subsection{Compactness}
The proof of the following lemma can be found in, e.g., \cite{NHarmS10Arxiv}.
\begin{lemma}[Compactness]\label{la:bs:compact}
Let $D \subset \R^n$ be a smoothly bounded domain, $s > 0$. Assume that there is a constant $C > 0$ and $f_k \in H^s(\R^n)$, $k \in \N$, such that for any $k \in \N$ the conditions $\supp f_k \subset \bar{D}$ and $\Vert f_k \Vert_{H^s} \leq C$ hold. Then there exists a subsequence $f_{k_i}$, such that $f_{k_i} \xrightarrow{i \to \infty} f \in H^s$ weakly in $H^s$, strongly in $L^2(\R^n)$, and pointwise almost everywhere. Moreover, $\supp f \subset \bar{D}$.
\end{lemma}

%% file: lowerorder.tex
\section{Lower Order Products: Proof of Lemma~\ref{la:commlowerorder}}
In this section, we prove that
\[
 H(u,v) = \lapn (uv) - u\lapn v - v \lapn u
\]
behaves in some sense like a product of lower order operators, as one can see immediately if $n \in 4\N$. As in \cite{NHarmS10Subm} we use elementary multiplier estimates derived in Section \ref{ss:multest} in order to give in Section \ref{ss:proofoflorder} the proof of Lemma~\ref{la:commlowerorder}. For the latter we restrict our attention to dimensions $n \in 2\N+1$ for the sake of decent size of presentation.

\subsection{Multiplier Estimates}\label{ss:multest}
Similar to the multiplier estimates in \cite{NHarmS10Subm}, we will need the following estimates which are again basically just consequences of the mean value theorem:
\begin{proposition}[Yet Another Silly Estimate]\label{pr:sillyest2}
Let $\lambda := \min\{\abs{a},\abs{b}\}> 0$. Then for any $s \geq 0$ there is a constant $C_s > 0$ such that.
\[
 \abs{\abs{a}^{-s} - \abs{b}^{-s}} \leq C_s\ \lambda^{-s-1}\ \abs{a-b}.
\]
\end{proposition}
% \begin{proofP}{\ref{pr:sillyest2}}
% Set $\Lambda := \max\{\abs{a},\abs{b}\}$. First of all, for any $t \geq 1$ we have by the mean value theorem
% \[
%  \frac{\abs{\abs{a}^t - \abs{b}^t}}{\abs{a-b}} \leq \sup_{\xi \in [a,b]} t \abs{\xi}^{t-1} \aleq \Lambda^{t-1}.
% \]
% If on the other hand $t \in (0,1)$, there exists a $K \in \N$ such that $2^K t \geq 1$. Then
% \[
%  \begin{ma}
%  &&\abs{\abs{a}^{t} - \abs{b}^{t}} \\
% &\leq& \abs{\abs{a}^{2t} - \abs{b}^{2t}}\ \Lambda^{-t}\\
% &\leq& \abs{\abs{a}^{4t} - \abs{b}^{4t}}\ \Lambda^{-t-2t}\\
% &\leq& \ldots\\
% &\leq& \abs{\abs{a}^{2^K t} - \abs{b}^{2^K t}}\ \Lambda^{-\sum_{k=0}^{K-1}2^k t}\\
% &\aleq{}& \abs{a-b}\ \Lambda^{2^K t-1} \Lambda^{-\sum_{k=0}^{K-1}2^k t}\\
% &=& \abs{a-b}\ \Lambda^{t-1}.
% \end{ma}
% \]
% This implies for any $s >0$
% \[
%  \abs{\abs{a}^{-s} - \abs{b}^{-s}} \aleq{} \abs{a-b}\ \abs{a}^{-s} \abs{b}^{-s}\ \Lambda^{s-1} \leq \abs{a-b}\ \lambda^{-s-1}.
% \]
% \end{proofP}

If one replaces in Proposition~\ref{pr:sillyest2} $\abs{\cdot}^{-s}$ by $m(\cdot)\abs{\cdot}^{-s}$ for some zero-multiplier $m$, the same result is to be expected. In order to prove this, first we have the following

\begin{proposition}[Silly Zero-Multiplier Estimates]\label{pr:sillyzerom}
Let $m(\cdot) \in C^\infty(\R^n \backslash \{0\})$ be a zero-homogeneous function. Then there is a constant $C_m$ such that for any $a \neq b \in \R^n \backslash \{0\}$, denoting $\Lambda := \max\{\abs{a},\abs{b}\}$
\[
  \abs{m(a)-m(b)} \leq C_m\ \frac{\abs{a-b}}{\Lambda}.
\]
\end{proposition}
\begin{proposition}[Multiplier-Estimates]\label{pr:multest2}
For any zero-multiplier $m \in C^\infty(\R^n \backslash \{0\})$ and any $s \geq 0$ there exists a constant $C_{m,s} > 0$ such that for any $a \neq b \in \R^n \backslash \{0\}$ we have for $\lambda := \min\{\abs{a},\abs{b}\}$
\begin{equation}\label{eq:multest:bfbest}
 \abs{\abs{a}^{-s} m(a)- \abs{b}^{-s} m(b)} \leq C_{m,s}\ \lambda^{-s-1}\ \abs{a-b}.
\end{equation}
If $\abs{a} \aeq \abs{b} \aeq \lambda$, for any $\varepsilon \in (0,1)$, $s \in \R$
\begin{equation}\label{eq:multest:bfbest2}
  \abs{\abs{a}^{s} m(a)- \abs{b}^{s} m(b)} \leq C_{m,s}\ \lambda^{s-\varepsilon}\ \abs{a-b}^\varepsilon.
\end{equation}
\end{proposition}
% \begin{proofP}{\ref{pr:multest2}}
% In order to prove \eqref{eq:multest:bfbest} we use that
% \[
% \begin{ma}
%  &&\abs{\abs{a}^{-s} m(a)- \abs{b}^{-s} m(b)}\\
% &\leq& \abs{m(a)-m(b)}\ \abs{a}^{-s} + \abs{m(a)}\ \abs{\abs{a}^{-s}-\abs{b}^{-s}}\\
% &\overset{\ontop{\sref{P}{pr:sillyzerom}}{\sref{P}{pr:sillyest2}}}{\aleq}& \frac{\abs{a-b}}{\lambda}\abs{a}^{-s} + \abs{a-b}\ \lambda^{-s-1}\\
% &\leq& \abs{a-b}\ \lambda^{-s-1}.
% \end{ma}
% \]
% For $s \leq 0$ it is obvious that \eqref{eq:multest:bfbest2} follows from \eqref{eq:multest:bfbest}. If $s > 0$, this follows from the same estimates as above, using that
% \[
%  \abs{\abs{a}^{s} - \abs{b}^{s} } \aeq{} \lambda^{2s}\ \abs{\abs{a}^{-s} - \abs{b}^{-s}}.
% \]
% \end{proofP}

\subsection{Proving lower order behavior}\label{ss:proofoflorder}
In order to give the proof of Lemma~\ref{la:commlowerorder}, we need the following intermediate result.
\begin{lemma}\label{la:lapmsblap12lapmta}
Let $s \in (0,\frac{n}{2})$, $M$ and $N$ zero-multiplier operators and $a,b \in \Sw(\R^n)$. Then,
\[
\begin{ma}
&&\int M \lapms{s} \abs{b}(y)\ \frac{\abs{N \lap^{-\frac{n+1}{4}+\frac{s}{2}}a(x)-N\lap^{-\frac{n+1}{4}+\frac{s}{2}}a(y)}}{\abs{x-y}^{n+\frac{1}{2}}}\ dy \\
&\aleq& \abs{\abs{M} \lapms{s+\delta} \abs{b}}(x)\ \lap^{-\frac{n}{4}+\frac{s}{2}+\frac{\delta}{2}}\abs{a}(x)
+\abs{\lap^{-\frac{\delta}{2}} \brac{\abs{\abs{M} \lapms{s} \abs{b}}\ \lap^{-\frac{n}{4}+\frac{s}{2}+\frac{\delta}{2}}\abs{a}}(x)},
\end{ma}
\]
for some $\delta \in \brac{0,\frac{n}{2}-s}$.
\end{lemma}
\begin{proofL}{\ref{la:lapmsblap12lapmta}}
We have
\[
\begin{ma}
 &&\int M \lapms{s} \abs{b}(y)\ \frac{\abs{N \lap^{-\frac{n+1}{4}+\frac{s}{2}}a(x)-N\lap^{-\frac{n+1}{4}+\frac{s}{2}}a(y)}}{\abs{x-y}^{n+\frac{1}{2}}}\ dy\\
&\aleq& \int \int M \lapms{s} \abs{b}(y)\ \abs{a(\xi)} \frac{\abs{n(x-\xi)\abs{x-\xi}^{-\frac{n}{2}+\frac{1}{2}-s}-n(y-\xi)\abs{y-\xi}^{-\frac{n}{2}+\frac{1}{2}-s}}}{\abs{x-y}^{n+\frac{1}{2}}}\ dy\ d\xi.
\end{ma}
\]
We set
\[
 k(x,y,\xi) := \frac{\abs{n(x-\xi)\abs{x-\xi}^{-\frac{n}{2}+\frac{1}{2}-s}-n(y-\xi)\abs{y-\xi}^{-\frac{n}{2}+\frac{1}{2}-s}}}{\abs{x-y}^{n+\frac{1}{2}}}.
\]
We decompose the space $(x,y,\xi) \in \R^{3n}$ into several subspaces depending on the relations of $\abs{y-\xi}$, $\abs{x-y}$, $\abs{x-\xi}$:
\[
 1 \leq \chi_1(x,y,\xi) + \chi_2(x,y,\xi) + \chi_3(x,y,\xi) + \chi_4(x,y,\xi) \quad \mbox{for $x,y,\xi \in \R^n$},
\]
where
\[
 \chi_1 := \chi_{\abs{x-y} \leq  2\abs{y-\xi}}\ \chi_{\abs{x-y} \leq 2 \abs{x-\xi}},
\]
\[
 \chi_2 := \chi_{\abs{x-y} \leq  2\abs{y-\xi}}\ \chi_{\abs{x-y} > 2 \abs{x-\xi}},
\]
\[
 \chi_3 := \chi_{\abs{x-y} >  2\abs{y-\xi}}\ \chi_{\abs{x-y} \leq 2 \abs{x-\xi}}\ \chi_{\abs{x-\xi} \leq 2 \abs{y-\xi}}.
\]
\[
 \chi_4 := \chi_{\abs{x-y} >  2\abs{y-\xi}}\ \chi_{\abs{x-y} \leq 2 \abs{x-\xi}}\ \chi_{\abs{x-\xi} > 2 \abs{y-\xi}}.
\]
In fact, if we assumed $\abs{x-y} > 2 \abs{x-\xi}$ and $\abs{x-y} > 2 \abs{y-\xi}$, then
\[
 \abs{x-y} \leq \abs{x-\xi} + \abs{y-\xi} < \fracm{2} \abs{x-y} + \fracm{2} \abs{x-y} = \abs{x-y},
\]
which is clearly impossible. Thus,
\[
\begin{ma}
 &&k(x,y,\xi)\\
&\leq& \chi_1(x,y,\xi)\ k(x,y,\xi) + \chi_2(x,y,\xi)\ k(x,y,\xi) + \chi_3(x,y,\xi)\ k(x,y,\xi) + \chi_4(x,y,\xi)\ k(x,y,\xi)\\
&=:& k_1(x,y,\xi) + k_2(x,y,\xi) + k_3(x,y,\xi) + k_4(x,y,\xi).
\end{ma}
\]
\underline{As for $k_1$}, note that
\[
 \abs{x-\xi}\chi_1 \leq  \abs{x-y} \chi_1+\abs{y-\xi} \chi_1 \leq 3 \abs{y-\xi} \chi_1 \leq \ldots \leq 9 \abs{x-\xi} \chi_1,
\]
that is $\abs{x-\xi} \chi_1 \aeq \abs{y-\xi} \chi_1$ for some uniform constants. Then, by Proposition~\ref{pr:multest2} for $\varepsilon := \frac{1}{2}+\delta$, $\delta \in (0,\frac{1}{2})$
\[
k_1(x,y,\xi) \aleq{} C_\delta \abs{x-y}^{-n+\delta} \abs{x-\xi}^{-\frac{n}{2}-s-\delta},
\]
and we choose $\delta < \frac{n}{2}-s$, i.e. small enough so that $-\frac{n}{2}-s-\delta > -n$.
\underline{As for $k_2$}, we have that
\[
\begin{ma}
 k_2(x,y,\xi) &\leq& \frac{\abs{x-\xi}^{-\frac{n}{2}-s+\frac{1}{2}} +\abs{y-\xi}^{-\frac{n}{2}-s+\frac{1}{2}}}{\abs{x-y}^{n+\frac{1}{2}}}\chi_2\\
&\aleq& \frac{\abs{x-\xi}^{-\frac{n}{2}-s+\frac{1}{2}}}{\abs{x-y}^{n+\frac{1}{2}}}\chi_2\\
&\aleq& \frac{\abs{x-\xi}^{-\frac{n}{2}-s-\delta}}{\abs{x-y}^{n-\delta}}\chi_2.
\end{ma}
\]
\underline{As for $k_3$}, we argue in the same way as for $k_2$
\[
\begin{ma}
 k_3(x,y,\xi) &\leq& \frac{\abs{x-\xi}^{-\frac{n}{2}-s+\frac{1}{2}} +\abs{y-\xi}^{-\frac{n}{2}-s+\frac{1}{2}}}{\abs{x-y}^{n+\frac{1}{2}}}\chi_3\\
&\aleq& \frac{\abs{x-\xi}^{-\frac{n}{2}-s+\frac{1}{2}}}{\abs{x-y}^{n+\frac{1}{2}}}\chi_3\\
&\aleq& \frac{\abs{x-\xi}^{-\frac{n}{2}-s-\delta}}{\abs{x-y}^{n-\delta}}\chi_3.
\end{ma}
\]
Consequently,
\[
 k_1(x,y,\xi) + k_2(x,y,\xi) + k_3(x,y,\xi) \aleq \frac{\abs{x-\xi}^{-\frac{n}{2}-s-\delta}}{\abs{x-y}^{n-\delta}}.
\]
Thus, for $i = 1,2,3$
\[
 \begin{ma}
&& \int \int \abs{M \lapms{s} b(y)}\ \abs{a(\xi)}\ k_i(x,y,\xi)\ dy\ d\xi\\
&\aleq& \int \int \abs{M} \lapms{s} \abs{b}(y)\ \abs{a(\xi)} \frac{\abs{x-\xi}^{-\frac{n}{2}-s-\delta}}{\abs{x-y}^{n-\delta}}\ dy\ d\xi\\
&\aeq& \lap^{-\frac{n}{4}+\frac{s}{2}+\frac{\delta}{2}} \abs{a}(x)\ \lapms{s+\delta}\abs{M}\abs{b}(x). \\
\end{ma}
\]
It remains to consider the \underline{case of $\chi_4 \neq 0$}: We have
\[
\begin{ma}
  k_4(x,y,\xi) &\leq& \abs{y-\xi}^{-\frac{n}{2}-s+\frac{1}{2}}\ \abs{x-y}^{-n-\frac{1}{2}} \chi_4 \leq \abs{y-\xi}^{-\frac{n}{2}-s-\delta}\ \abs{x-y}^{-n+\delta}.
\end{ma}
\]
Thus,
\[
 \begin{ma}
&& \int \int \abs{M \lapms{s} b(y)}\ \abs{a(\xi)}\ k_4(x,y,\xi)\ dy\ d\xi\\
&\aleq& \int \int \abs{M \lapms{s} b(y)}\ \abs{a(\xi)}\ \abs{y-\xi}^{-\frac{n}{2}-s-\delta}\ \abs{x-y}^{-n+\delta}\ dy\ d\xi\\
&\aeq&  \int \abs{M \lapms{s} b(y)}\ \lap^{-\frac{n}{4}+\frac{s}{2}+\frac{\delta}{2}}\abs{a}(y)\ \abs{x-y}^{-n+\delta}\ dy\\
&\aeq&  \lap^{-\frac{\delta}{2}} \brac{\abs{M \lapms{s} b}\ \lap^{-\frac{n}{4}+\frac{s}{2}+\frac{\delta}{2}}\abs{a}}(x)\\
\end{ma}
\]
\end{proofL}

Now we are able to give the
\begin{proofL}{\ref{la:commlowerorder}}
We prove only the case where $n-1$ is divisible by $4$ and $n \geq 5$. Generally,
\begin{equation}\label{eq:14prodrule}
 \lap^{\frac{1}{4}} (uv)(x) = v(x) \lap^{\frac{1}{4}} u(x) + c_n \int u(y)\ \frac{v(x)-v(y)}{\abs{x-y}^{n+\frac{1}{2}}}\ dy.
\end{equation}
Set $K = \lfloor \frac{n}{4} \rfloor = \frac{n-1}{4} \in \N$. Let $\gamma$ denote multi-indices $\gamma \in \brac{\N_0}^n$, and let $M_\gamma$, $N_\gamma$ be certain zero-multiplier operators such that
\[
\begin{ma}
 \lapn \brac{\lapmn a\ \lapmn v}
&=& \lap^{\frac{1}{4}} \sum_{\abs{\gamma} = 1}^{2K-1} M_\gamma \lap^{\frac{2\abs{\gamma}-n}{4}} a\ N_\gamma \lap^{\frac{4K-2\abs{\gamma}-n}{4}} b\\
&& + \lap^{\frac{1}{4}} \brac{\lapmn a\ \lap^{\frac{4K-n}{4}} b} \\
&& + \lap^{\frac{1}{4}} \brac{\lap^{\frac{4K-n}{4}} a\ \lapmn b} \\
&\overset{\eqref{eq:14prodrule}}{=}& \lap^{\frac{1}{4}} \sum_{\abs{\gamma} = 1}^{2K-1} M_\gamma \lap^{\frac{2\abs{\gamma}-n}{4}} a\ N_\gamma \lap^{\frac{4K-2\abs{\gamma}-n}{4}} b\\
&& + \lapmn  a\ b + c_n \int \lap^{\frac{4K-n}{4}} b(y)\ \frac{\lapmn  a(x)-\lapmn  a(y)}{\abs{x-y}^{n+\frac{1}{2}}}\ dy\\
&& + a\ \lapmn  b + c_n \int \lap^{\frac{4K-n}{4}} a(y)\ \frac{\lapmn  b(x)- \lapmn  b(y)}{\abs{x-y}^{n+\frac{1}{2}}}\ dy.
\end{ma}
\]
Thus,
\[
\begin{ma}
 \lapn \brac{\lapmn a\ \lapmn b} - a\ \lapmn b - b\ \lapmn a
&=& \lap^{\frac{1}{4}} \sum_{\abs{\gamma} = 1}^{2K-1} M_\gamma \lap^{\frac{2\abs{\gamma}-n}{4}} a\ N_\gamma \lap^{\frac{4K-2\abs{\gamma}-n}{4}} b\\
&& + c_n \int \lap^{\frac{4K-n}{4}} b(y)\ \frac{\lapmn  a(x)-\lapmn  a(y)}{\abs{x-y}^{n+\frac{1}{2}}}\ dy\\
&& + c_n \int \lap^{\frac{4K-n}{4}} a(y)\ \frac{\lapmn  b(x)- \lapmn  b(y)}{\abs{x-y}^{n+\frac{1}{2}}}\ dy\\
&=& \sum_{\abs{\gamma} = 1}^{\frac{n-3}{2}} \lap^{\frac{1}{4}} \brac{M_\gamma \lap^{\frac{2\abs{\gamma}-n}{4}} a\ N_\gamma \lap^{\frac{-1-2\abs{\gamma}}{4}} b}\\
&& + c_n \int \lap^{-\frac{1}{4}} b(y)\ \frac{\lapmn a(x)-\lapmn a(y)}{\abs{x-y}^{n+\frac{1}{2}}}\ dy\\
&& + c_n \int \lap^{-\frac{1}{4}} a(y)\ \frac{\lapmn b(x)- \lapmn b(y)}{\abs{x-y}^{n+\frac{1}{2}}}\ dy\\
&=:& \sum_{\abs{\gamma} = 1}^{\frac{n-3}{2}} I_\gamma + II + III.
\end{ma}
\]
Note that
\[
\begin{ma}
 &&\abs{I_\gamma(x)}\\
&\overset{\eqref{eq:14prodrule}}{\aleq}& \abs{M_\gamma \lap^{\frac{2\abs{\gamma}+1-n}{4}} a(x)}\ \abs{N_\gamma \lap^{\frac{-1-2\abs{\gamma}}{4}} b(x)} \\
&& + \abs{\int M_\gamma \lap^{\frac{2\abs{\gamma}-n}{4}} a(y)\ \frac{N_\gamma \lap^{\frac{-1-2\abs{\gamma}}{4}} b(x)-N_\gamma \lap^{\frac{-1-2\abs{\gamma}}{4}} b(y)}{\abs{x-y}^{n+\frac{1}{2}}}\ dy}.
\end{ma}
\]
Because of $1 \leq \abs{\gamma} \leq \frac{n-3}{2}$ we have altogether for some constants $L \in \N$, $s_k \in (0,\frac{n}{2})$
\[
 \begin{ma}
 &&H(\lapmn a,\lapmn b)(x) \\
&\aleq{}& \sum_{k = 1}^{L_1} M_k \lapms{s_k} \abs{a}(x)\ N_k\lap^{-\frac{n}{4}+\frac{s_k}{2}} \abs{b}(x)\\
&& + \int M_k \lapms{s_k} \abs{a}(y)\ \frac{\abs{N_k \lap^{-\frac{n+1}{4}+\frac{s_k}{2}}b(x)-N_k\lap^{-\frac{n+1}{4}+\frac{s_k}{2}}b(y)}}{\abs{x-y}^{n+\frac{1}{2}}}\ dy \\
&& + \int N_k \lapms{s_k} \abs{b}(y)\ \frac{\abs{M_k \lap^{-\frac{n+1}{4}+\frac{s_k}{2}}a(x)- M_k\lap^{-\frac{n+1}{4}+\frac{s_k}{2}}a(y)}}{\abs{x-y}^{n+\frac{1}{2}}}\ dy.
\end{ma}
\]
Now one estimates the integral terms with Lemma~\ref{la:lapmsblap12lapmta}, and concludes.
\end{proofL}
\begin{remark}\label{rem:lowordest}
In particular, Proposition~\ref{pr:lpestloop} is applicable, and we have for $a := \lapn u$, $b := \lapn v$,
\begin{equation}\label{eq:HuvinL21}
 \Vert H(u,v) \Vert_{(2,1),\R^n} \aleq \Vert \lapn u \Vert_{2,\R^n}\ \Vert \lapn v \Vert_{2,\R^n},
\end{equation}
for $\fracm{q} = \fracm{q_1} + \fracm{q_2}$. In fact, this holds whenever $u$, $v$, $\lapn u$, $\lapn v \in L^2(\R^n)$, via approximation $u_k \xrightarrow{k \to \infty} u$, $v_k \xrightarrow{k \to \infty} v$ in $H^{\frac{n}{2}}(\R^n)$, because then
\[
 \Vert H(u_k,v_k) - H(u_l,v_l) \Vert_{(2,1),\R^n} \xrightarrow{ k,l \to \infty} 0,
\]
because of bilinearity of $H(\cdot,\cdot)$, and as the pointwise limit of $H(u_k,v_k) = H(u,v)$, we have that $H(u,v) \in L^{2,1}(\R^n)$ and \eqref{eq:HuvinL21} holds.
In the same way one can show that
\[
\Vert H(u,v) \Vert_{2,\R^n} \aleq \Vert \lapn u \Vert_{(2,\infty),\R^n}\ \Vert \lapn v \Vert_{2,\R^n}.
\]
\end{remark}

%% file: cutoff.tex
\section{General estimates and inequalities}

\subsection{Estimates on commutator-operators with cut-off functions}
\begin{proposition}\label{pr:lapmsabslapsteta}
For any $s,t \in (0,n)$, and a zero-multiplier operator $M$
\[
 \left \Vert M\lapms{s} \abs{\laps{t}(1-\eta_r)} \right \Vert_{\infty,\R^n} \leq C_{t,s}\ r^{-t+s},
\]
\end{proposition}
\begin{proofP}{\ref{pr:lapmsabslapsteta}}
This follows by scaling, once we prove that
\[
 \left \Vert M\lapms{s} \abs{\laps{t}(1-\eta_1)} \right \Vert_{\infty,\R^n} \overset{t > 0}{=} \left \Vert M\lapms{s} \abs{\laps{t}\eta_1} \right \Vert_{\infty,\R^n} \aleq C_{s,t,M}.
\]
But this again follows from the fact that $\abs{\laps{t}\eta_1} \in L^(p,q)(\R^n)$ for any $p \in (1,\infty)$, $q \in (1,\infty)$, as shown in Proposition~\ref{pr:finswlapsinlp}: In particular, $\abs{\laps{t}\eta_1} \in L^{\frac{n}{s},1}(\R^n)$, and one concludes via Lemma~\ref{la:soblinfty}.
\end{proofP}

\begin{proposition}\label{pr:estvlapn1metalampnvp}
For all $p \in (1,\infty)$, $q \in [1,\infty]$, there is a constant $C_{p,q} > 0$ and $\gamma \equiv \gamma_p > 0$ such that
\[
 \Vert \brac{\lapn \eta_{\Lambda r}}\ \lapmn \varphi \Vert_{(p,q),\R^n} = \Vert \brac{\lapn (1-\eta_{\Lambda r})}\ \lapmn \varphi \Vert_{(p,q),\R^n} \leq C_{p,q}\ \Lambda^{-\gamma} \Vert \varphi \Vert_{(p, \infty)}
\]
for all $\Lambda \geq 1$ and $\varphi \in C_0^\infty(B_r)$, $B_r \subsubset \R^n$.
\end{proposition}
\begin{proofP}{\ref{pr:estvlapn1metalampnvp}}
We have,
\[
 \lapn(1-\eta_{\Lambda r}) = -\lapn \eta_{\Lambda r}.
\]
Let $p \in (1,\infty)$, pick $p_1,p_2 \in (2,\infty)$ such that
\[
 \fracm{p} = \fracm{p_1} + \fracm{p_2},
\]
\[
 \frac{1}{2}+\fracm{p_2} \in \brac{\fracm{p},1}.
\]
We set $p_3 \in (1,p)$ such that
\[
 \fracm{2} + \fracm{p_3} = 1 + \fracm{p_2}.
\]
Then,
\[
\begin{ma}
 \Vert \brac{\lapn \eta_{\Lambda r}}\ \lapmn \varphi \Vert_{(p,q),\R^n}
&\leq& \Vert \lapn \eta_{\Lambda r}\Vert_{(p_1,q),\R^n}\ \Vert \lapmn \varphi \Vert_{(p_2,\infty),\R^n}\\
&\overset{\sref{P}{pr:finswlapsinlp}}{\aleq}& \brac{\Lambda r}^{-\frac{n}{2}+\frac{n}{p_1}}\ \Vert \varphi \Vert_{(p_3,\infty),\R^n}\\
&\aleq& \brac{\Lambda r}^{-\frac{n}{2}+\frac{n}{p_1}}\ r^{\frac{n}{p_3}-\frac{n}{p}}\ \Vert \varphi \Vert_{(p,\infty),\R^n}\\
&=& \Lambda^{-\frac{n}{2}+\frac{n}{p_1}}\ \brac{r^n}^{-\fracm{2}+\fracm{p_1} + \fracm{p_3}-\fracm{p}}\ \Vert \varphi \Vert_{(p,\infty),\R^n}.
\end{ma}
\]
Now,
\[
 -\fracm{2}+\fracm{p_1} + \fracm{p_3}-\fracm{p} =  -\fracm{2}-\fracm{p_2} + \fracm{2} + \fracm{p_2} = 0.
\]
\end{proofP}

\begin{proposition}\label{pr:h1metalapmnvp}
There is $\gamma > 0$ and for any $q \in [1,\infty]$ there is a constant $C_q > 0$ such that for any $\varphi \in C_0^\infty(B_r)$, $\Lambda \geq 1$, $B_r \subsubset \R^n$
\[
 \Vert H((1-\eta_{\Lambda r}),\lapmn \varphi) \Vert_{(2,q),\R^n} = \Vert H((1-\eta_{\Lambda r}),\lapmn \varphi) \Vert_{(2,q),\R^n} \leq C_q\ \Lambda^{-\gamma}\ \Vert \varphi \Vert_{(2,\infty),\R^n}.
\]
\end{proposition}
\begin{proofP}{\ref{pr:h1metalapmnvp}}
We have
\[
 \begin{ma}
  \Vert H((1-\eta_{\Lambda r}),\lapmn \varphi) \Vert_{(2,q),\R^n} &\overset{\sref{P}{pr:lorentzfacts}}{\aleq}& \sup_{\ontop{\psi \in C_0^\infty(\R^n)}{\Vert \psi \Vert_{(2,q')} \leq 1}} \intl_{\R^n} H((1-\eta_{\Lambda r}),\lapmn \varphi)\ \psi.
 \end{ma}
\]
For such a $\psi$ using Lemma~\ref{la:commlowerorder} and the fact that
\[
 H((1-\eta_{\Lambda r}),\lapmn \varphi) = H(\eta_{\Lambda r},\lapmn \varphi) = H(\lapmn \lapn \eta_{\Lambda r},\lapmn \varphi)
\]
we have for certain $s_k \in (0,\frac{n}{2})$ and $t_k \in (0,s_k]$
\[
\begin{ma}
 &&\int H((1-\eta_{\Lambda r}),\lapmn \varphi)\ \psi\\
&\aleq& \sum_{k=1}^L \int \lapms{s_k-t_k} \brac{M_{k,1} \lapms{t_k} \abs {\lapn \eta_{\Lambda r}}\ M_{k,2} \lap^{-\frac{n}{4}+\frac{s_k}{2}} \abs{\varphi}}\ \abs{\psi}\\
&=& \sum_{k=1}^L \int M_{k,1} \lapms{t_k} \abs {\lapn \eta_{\Lambda r}}\ M_{k,2} \lap^{-\frac{n}{4}+\frac{s_k}{2}} \abs{\varphi}\ \lapms{s_k-t_k} \abs{\psi}\\
&\overset{\sref{P}{pr:lapmsabslapsteta}}{\aleq}& \sum_{k=1}^L \brac{\Lambda r}^{t_k-\frac{n}{2}}\ \Vert \lap^{-\frac{n}{4}+\frac{s_k}{2}} \abs{\varphi} \Vert_{(p_1,q),\R^n}\ \Vert \lapms{s_k-t_k} \abs{\psi} \Vert_{(p_2,q'),\R^n}.\\
\end{ma}
\]
Here, (fixing $s := s_k$, $t := t_k$)
\[
 \fracm{p_1} = \frac{1}{2} + \frac{s-t}{n} \in \left. \left [ \frac{1}{2},1 \right . \right ),
\]
and
\[
 \fracm{p_2} = \frac{1}{2} - \frac{s-t}{n} \in \left. \left (0,\frac{1}{2} \right . \right ].
\]
Then,
\[
 \Vert \lapms{s_k-t_k} \abs{\psi} \Vert_{(p_2,q'),\R^n} \aleq \Vert \psi \Vert_{(2,q'),\R^n} \leq 1.
\]
and for $\fracm{p_3} = 1- \frac{t}{n} \in \brac{\fracm{2},1}$
\[
 \Vert \lap^{-\frac{n}{4}+\frac{s_k}{2}} \abs{\varphi} \Vert_{(p_1,q),\R^n} \aleq \Vert \varphi \Vert_{(p_3,q),\R^n}\\
\leq C_{p_2}\ r^{n\brac{\fracm{2}-\frac{t}{n}}}\ \Vert \varphi \Vert_{(2,\infty),\R^n}\\
\aeq{} r^{\frac{n}{2}-t}\ \Vert \varphi \Vert_{(2,\infty),\R^n}.
\]
\end{proofP}

% \begin{lemma}\label{la:estvlapn1metalampnvp}
% For any $\varphi \in C_0^\infty(B_r(0))$
% \[
% \Vert \brac{\lapn \brac{1-\eta_{\Lambda r}}}\ \lapmn \varphi\Vert_{2,\R^n} \aleq  \Lambda^{-\frac{n}{4}} \Vert \varphi \Vert_{(2,\infty),\R^n}
% \]
% \end{lemma}
% \begin{proofL}{\ref{la:estvlapn1metalampnvp}}
% \[
% \begin{ma}
%  &&\Vert \lapn \brac{1-\eta_{\Lambda r}} \Vert_{4,\R^n} \Vert \lapmn \varphi\Vert_{4,\R^n}\\
% &\aleq & \brac{\Lambda r}^{-\frac{n}{2}+\frac{n}{4}}\ \Vert \varphi \Vert_{\frac{4}{3},B_r}\\
% &\aleq & \brac{\Lambda r}^{-\frac{n}{2}+\frac{n}{4}}\ r^{\frac{n}{4}}\ \Vert \varphi \Vert_{(2,\infty),B_r}\\
% &=& \Lambda^{-\frac{n}{4}}\ \Vert \varphi \Vert_{(2,\infty),B_r}.\\
% \end{ma}
% \]
% \end{proofL}

Proposition~\ref{pr:estvlapn1metalampnvp} and Proposition~\ref{pr:h1metalapmnvp} imply in particular
\begin{proposition}\label{pr:estlapn1metalpmnvp}
For all $\varphi \in C_0^\infty(B_r)$, $\Lambda \geq 1$, $q \in [1,\infty]$
\[
 \left \Vert \lapn \brac{\brac{1-\eta_{\Lambda r}}\ \lapmn \varphi} \right \Vert_{(2,q),\R^n} \leq C_q\ \Lambda^{-\gamma}\ \Vert \varphi \Vert_{(2,\infty),\R^n},
\]
and for any $\Lambda \geq 1$,
\[
 \left \Vert \lapn \brac{\eta_{\Lambda r}\ \lapmn \varphi} \right \Vert_{(2,q),\R^n} \leq C_{q}\ \Vert \varphi \Vert_{(2,q),\R^n}.
\]
\end{proposition}
\begin{proofP}{\ref{pr:estlapn1metalpmnvp}}
This follows from
\[
 \lapn \brac{\brac{1-\eta_{\Lambda r}}\ \lapmn \varphi} = H(\brac{1-\eta_{\Lambda r}}, \lapmn \varphi) + \brac{\lapn \brac{1-\eta_{\Lambda r}}}\ \lapmn \varphi + 0,
\]
and
\[
 \lapn \brac{\eta_{\Lambda r}\ \lapmn \varphi} = H(\eta_{\Lambda r}, \lapmn \varphi) + \brac{\lapn \eta_{\Lambda r}}\ \lapmn \varphi + \eta_{\Lambda r} \varphi.
\]
\end{proofP}

\subsection{Localizing Effects of Locally Supported Functions}
\begin{proposition}\label{pr:lapmsvarphiakal}
For $\varphi \in C_0^\infty(A_{r,k})$, for any $l \in \N_0$, if $\abs{l-k} \geq 2$
\[
 \Vert M \lapms{s} \varphi \Vert_{(p_1,q),A_{r,l}} \leq C_{p_1,q}\ 2^{\max\{l,k\} (-n+s) + k \brac{n-\frac{n}{p_1}-s} + l \frac{n}{p_1}}\ \Vert \varphi \Vert_{(p_2,\infty),\R^n},
\]
where $\fracm{p_1} = \fracm{p_2} -\frac{s}{n} \in (0,1)$.
\end{proposition}
\begin{proofP}{\ref{pr:lapmsvarphiakal}}
For some $\psi \in C_0^\infty(A_{r,l})$, $\Vert \psi \Vert_{(p_1',q),\R^n} \leq 1$ we have to estimate
\[
\begin{ma}
\int M \lapms{s}\varphi\ \psi
&\aleq& \int \abs{\cdot}^{-n+s} \abs{\varphi}\ast \abs{\psi}\\
&\aleq& \brac{2^{\max\{l,k\}} r}^{-n+s} \Vert \varphi \Vert_{1,\R^n}\ \Vert \psi \Vert_{1,\R^n}\\
&\aleq& \brac{2^{\max\{l,k\}} r}^{-n+s}\ \brac{2^k r}^{\frac{n}{p_2'}}\ \Vert \varphi \Vert_{(p_2,\infty)}\ \brac{2^l r}^{\frac{n}{p_1}}\\
&\aeq& 2^{\max\{l,k\} (-n+s) + k \brac{n-\frac{n}{p_1}-s} + l \frac{n}{p_1}}\ \Vert \varphi \Vert_{(p_2,\infty)}.
\end{ma}
\]
\end{proofP}

\begin{proposition}\label{pr:lapnvarphibralrk}
For any zero-multiplier operator $M$ and any $p \in (1,\infty)$, $q \in [1,\infty]$ there is a constant $C_{M,p,q}$ such that for any $\varphi \in C_0^\infty(B_r)$, $\Lambda \geq 8$, $k \in \N$ the following holds,
\[
 \Vert M \lapn \varphi \Vert_{(p,q),A_{\Lambda r,k}} \leq C_{M, p,q} \brac{2^k \Lambda}^{-\brac{\frac{3}{2}-\fracm{p}}n} \Vert \lapn \varphi \Vert_{(p,q),\R^n}.
\]
\end{proposition}
\begin{proofP}{\ref{pr:lapnvarphibralrk}}
We have to estimate for some $\psi \in C_0^\infty(A_{\Lambda r,k})$, $\Vert \psi \Vert_{(p',q'),\R^n} \leq 1$
\[
\begin{ma}
 \int \abs{\cdot}^{-\frac{3}{2}n} \abs{\varphi} \ast \abs{\psi}
&\aleq& \brac{2^k \Lambda r}^{-\frac{3}{2}n}\ \Vert \varphi \Vert_{1,\R^n}\ \Vert \psi \Vert_{1,\R^n}\\
&\aleq& \brac{2^k \Lambda r}^{-\frac{3}{2}n}\ r^\brac{\frac{n}{p'}}\ \Vert \varphi \Vert_{(p,q),\R^n}\ \brac{2^k \Lambda r}^{\frac{n}{p}}\ \Vert \psi \Vert_{(p',q'),\R^n}\\
&\overset{\sref{L}{la:poinc}}{\aleq}& \brac{2^k \Lambda r}^{-\frac{3}{2}n}\ r^\brac{\frac{n}{p'}+\frac{n}{2}}\ \Vert \lapn \varphi \Vert_{(p,q),\R^n}\ \brac{2^k \Lambda r}^{\frac{n}{p}}\ \Vert \psi \Vert_{(p',q'),\R^n}\\
&=& \brac{2^k \Lambda}^{-\brac{\frac{3}{2}-\fracm{p}}n}\  \Vert \lapn \varphi \Vert_{(p,q),\R^n}.
\end{ma}
\]
\end{proofP}

\begin{proposition}\label{pr:lapmstlapnvarphibralrk}
Let $\varphi \in C_0^\infty(B_r)$, $t \in (0,\frac{n}{2})$, $\Lambda \geq 8$, and $k \in \N$. Then,
\[
 \Vert M \lapms{t} \abs{\lapn \varphi} \Vert_{(p_1,q),A_{\Lambda r, k}} \leq C_{M,p,q}\ \brac{2^{k} \Lambda}^{-\frac{3}{2}n + \frac{n}{p_2}}\ \Vert \lapn \varphi \Vert_{(p_2,q),\R^n},
\]
and
\[
 \Vert M \lapms{t} \abs{\lapn \varphi} \Vert_{(p_1,q),B_{\Lambda r}} \leq C_{M,p,q}\ \Vert \lapn \varphi \Vert_{(p_2,q),\R^n},
\]
where $\fracm{p_1} = \fracm{p_2} -\frac{t}{n} \in (0,1)$.
\end{proposition}
\begin{proofP}{\ref{pr:lapmstlapnvarphibralrk}}
We have
\[
\begin{ma}
  \Vrac{ M \lapms{t} \abs{\lapn \varphi} }_{(p_1,q),A_{\Lambda r, k}}
&\aleq& \sum_{l=0}^\infty \Vrac{ M \lapms{t} \abs{\eta_{\frac{\Lambda}{2}r}^l \lapn \varphi} }_{(p_1,q),A_{\Lambda r, k}}\\
&\aleq& \sum_{l=k}^{k+2} \Vrac{ M \lapms{t} \abs{\eta_{\frac{\Lambda}{2}r}^l \lapn \varphi} }_{(p_1,q),A_{\Lambda r, k}}\\
&& +\sum_{l=0}^{k-1} \Vrac{ M \lapms{t} \abs{\eta_{\frac{\Lambda}{2}r}^l \lapn \varphi} }_{(p_1,q),A_{\Lambda r, k}}\\
&& +\sum_{l=k+2}^{\infty} \Vrac{ M \lapms{t} \abs{\eta_{\frac{\Lambda}{2}r}^l \lapn \varphi} }_{(p_1,q),A_{\Lambda r, k}}\\
&=:& \sum_{l=k}^{k+2} I_l +\sum_{l=0}^{k-1} II_l +\sum_{l=k+2}^{\infty} III_l.
\end{ma}
\]
As for $I_l$, i.e. $l \aeq k$, $l \geq 1$, we have
\[
\begin{ma}
 \Vert M \lapms{t} \abs{\eta_{\frac{\Lambda}{2}r}^l \lapn \varphi} \Vert_{(p_1,q),A_{\Lambda r, k}}
&\aleq& \Vert \eta_{\frac{\Lambda}{2}r}^l \lapn \varphi \Vert_{(p_2,q),\R^n}\\
&\leq& \Vert \lapn \varphi \Vert_{(p_2,q),A_{\frac{\Lambda}{2} r,l}}\\
&\overset{\sref{P}{pr:lapnvarphibralrk}}{\aleq}& \brac{2^l \Lambda}^{-\brac{\frac{3}{2}-\frac{1}{p_2}}n}\ \Vert \lapn \varphi \Vert_{(p_2,q),A_{\frac{\Lambda}{2} r,l}}\\
&\aeq& \brac{2^k \Lambda}^{-\brac{\frac{3}{2}-\frac{1}{p_2}}n}\ \Vert \lapn \varphi \Vert_{(p_2,q),\R^n}.
\end{ma}
\]
For $II_l$, $III_l$, we have by the different support
\[
\begin{ma}
\Vert M \lapms{t} \abs{\eta_{\frac{\Lambda}{2}r}^l \lapn \varphi} \Vert_{(p_1,q),A_{\Lambda r, k}}
&\aleq& \brac{2^{\max\{k,l\}}\Lambda r}^{-n+t}\ \brac{2^k \Lambda r}^{\frac{n}{p_1}}\ \brac{2^l \Lambda r}^{\frac{n}{p_2'}}\ \Vert \lapn \varphi \Vert_{(p_2,\infty),A_{\frac{\Lambda}{2} r,l}}\\
&\overset{\sref{P}{pr:lapnvarphibralrk}}{\aleq}& \brac{2^{\max\{k,l\}}}^{-n+t}\ \brac{2^k }^{\frac{n}{p_1}}\ \brac{2^l}^{\frac{n}{p_2'}}\ \brac{2^l \Lambda}^{-\brac{\frac{3}{2}-\fracm{p_2}}n} \Vert \lapn \varphi\Vert_{(p_2,\infty)}.
\end{ma}
\]
The claim now follows since
\[
\begin{ma}
 \sum_{l=k}^\infty 2^{l(-n+t) + k \frac{n}{p_1} + l\frac{n}{p_2'} - \frac{3}{2}nl + l\frac{n}{p_2}}
&=&  \sum_{l=k}^\infty 2^{k \frac{n}{p_1} + l(t - \frac{3}{2}n)}\\
&\aeq&  2^{k(- \frac{3}{2}n + \frac{n}{p_2})}\\
\end{ma}
\]
and
\[
 \begin{ma}
 \sum_{l=0}^k 2^{k(-n+t) + k \frac{n}{p_1} + l\frac{n}{p_2'} - \frac{3}{2}nl + l\frac{n}{p_2}}
&=& \sum_{l=0}^k 2^{k(-n+t) + k \frac{n}{p_1} -l\frac{n}{2}}\\
&\aeq& 2^{k(-\frac{3}{2}n + \frac{n}{p_2})},
\end{ma}
\]
\end{proofP}

Finally, we are able to have the following
\begin{lemma}\label{la:estHvpfgL1}
Let $\varphi \in C_0^\infty(B_r)$, $f \in H^{\frac{n}{2}}(\R^n)$, $g \in L^2(\R^n)$. Then for all $\Lambda \geq 50$,
\[
\begin{ma}
 \Vert H(\varphi,f) g \Vert_{1,\R^n} &\leq& C\ \Vert \lapn \varphi \Vert_{2,\R^n}\ \brac{\Vert \lapn f \Vert_{2,B_{2\Lambda^3 r}} + \Lambda^{-\gamma} \Vert \lapn f \Vert_{2,\R^n}}\ \Vert g \Vert_{(2,\infty),B_{2\Lambda r}}\\
&&+ C\ \Lambda^{-\gamma}\ \sum_{k=1}^\infty 2^{-\gamma k}\ \Vert \lapn \varphi \Vert_{2,\R^n}\ \Vert \lapn f \Vert_{2,\R^n}\ \Vert g \Vert_{(2,\infty),A_{k,\Lambda r}}.
\end{ma}
\]
\end{lemma}
\begin{proofL}{\ref{la:estHvpfgL1}}
As always, we have
\[
\begin{ma}
\Vert H(\varphi,f) g \Vert_{1,\R^n}
&\aleq& \Vert H (\varphi,f) \Vert_{(2,1),B_{2\Lambda r}}\ \Vert g \Vert_{(2,\infty),B_{2\Lambda r}} + \sum_{k=1}^\infty \Vert H(\varphi,f) \Vert_{(2,1),A_{k,\Lambda r}}\ \Vert g \Vert_{(2,\infty),A_{k,\Lambda r}}\\
&=:& I\ \Vert g \Vert_{(2,\infty),B_{2\Lambda r}} + \sum_{k=1}^\infty II_k\ \Vert g \Vert_{(2,\infty),A_{k,\Lambda r}}.
\end{ma}
\]
As \underline{for $II_k$}, by Lemma~\ref{la:commlowerorder} we have to estimate terms of the following form for some $\psi \in C_0^\infty(A_{k,\Lambda r})$, $\Vert \psi \Vert_{(2,\infty),\R^n} \leq 1$, $s \in (0,\frac{n}{2})$, $t \in (0,s]$:
\[
\begin{ma}
&&\int M_1 \lapms{s-t} \psi\ M_2 \lapms{t} \abs{\lapn \varphi}\ M_3 \lap^{-\frac{n}{4}+\frac{s}{2}}\abs{\lapn f}\\
&=& \sum_{l=0}^\infty \int M_1 \lapms{s-t} \psi\ \eta^l_{\frac{\Lambda}{4}r}\ M_2\lapms{t} \abs{\lapn \varphi}\ M_3 \lap^{-\frac{n}{4}+\frac{s}{2}}\abs{\lapn f}\\
&=:& II_{k,0} + \sum_{l=2}^\infty \int M_1 \lapms{s-t} \psi\ \eta^l_{\frac{\Lambda}{4}r}\ M_2\lapms{t} \abs{\lapn \varphi}\ M_3 \lap^{-\frac{n}{4}+\frac{s}{2}}\abs{\lapn f}\\

&\aleq& II_{k,0}\\
&&+ \sum_{l=2}^\infty \Vert M_1 \lapms{s-t} \psi \Vert_{(\frac{2n}{n-2(s-t)},\infty),A_{\frac{\Lambda}{4} r,l}}\ \Vert M_2 \lapms{t} \abs{\lapn \varphi} \Vert_{(\frac{2n}{n-2t},2),A_{\frac{\Lambda}{4} r,l}}
\ \Vert M_3 \lap^{-\frac{n}{4}+\frac{s}{2}} \abs{\lapn f} \Vert_{(\frac{n}{s},2),\R^n}\\
&\aleq& II_{k,0} + \sum_{l=2}^\infty  \Vert M_1 \lapms{s-t} \psi \Vert_{(\frac{2n}{n-2(s-t)},\infty),A_{\frac{\Lambda}{4} r,l}}\ 
\Vert M_2 \lapms{t} \abs{\lapn \varphi} \Vert_{(\frac{2n}{n-2t},2),A_{\frac{\Lambda}{4} r,l}} \ \Vert \lapn f \Vert_{2,\R^n}\\
\end{ma}
\]
By Proposition~\ref{pr:lapmsvarphiakal} if $\abs{k-l} \geq 2$,
\[
\Vert M_1 \lapms{s-t} \psi \Vert_{(\frac{2n}{n-2(s-t)},\infty),A_{\frac{\Lambda}{4} r,l}} \aleq
2^{\max\{l,k\} (-n+s-t) + k \frac{n}{2} + l \frac{n-2(s-t)}{2}}\ \Vert \psi \Vert_{(2,\infty),\R^n},
\]
and else
\[
\Vert M_1 \lapms{s-t} \psi \Vert_{(\frac{2n}{n-2(s-t)},\infty),A_{\frac{\Lambda}{4} r,l}} \aleq \Vert \psi \Vert_{(2,\infty),\R^n}.
\]
Moreover, by Proposition~\ref{pr:lapmstlapnvarphibralrk} if $\abs{l} \geq 2$.
\[
 \Vert M_2 \lapms{t} \abs{\lapn \varphi} \Vert_{(\frac{2n}{n-2t},2),A_{\frac{\Lambda}{4} r,l}} \aleq
\brac{2^l \Lambda}^{-n}\ \Vert \lapn \varphi \Vert_{2,\R^n}.
\]
Consequently,
\[
\begin{ma}
&&\int M_1 \lapms{s-t} \psi\ M_2 \lapms{t} \abs{\lapn \varphi}\ M_3 \lap^{-\frac{n}{4}+\frac{s}{2}}\abs{\lapn f}\\
&\aleq& II_{k,0} + \Lambda^{-n} \sum_{l=2}^{k-1}  2^{k (-n+s-t) + k \frac{n}{2} + l \frac{n-2(s-t)}{2}-ln}\ \Vert \lapn \varphi \Vert_{2,\R^n}\ \Vert \lapn f \Vert_{2,\R^n}\\
&& + \Lambda^{-n} \sum_{l=k+1}^\infty 2^{l(-n+s-t) + k \frac{n}{2} + l \frac{n-2(s-t)}{2}-ln}\ \Vert \lapn \varphi \Vert_{2,\R^n}\ \Vert \lapn f \Vert_{2,\R^n}\\
&& + \Lambda^{-n}\ 2^{-kn} \Vert \lapn \varphi \Vert_{2,\R^n}\ \Vert \lapn f \Vert_{2,\R^n}\\
&\aeq& II_{k,0}  + \Lambda^{-n}\ 2^{k (-\frac{n}{2}+s-t)}\ \Vert \lapn \varphi \Vert_{2,\R^n}\ \Vert \lapn f \Vert_{2,\R^n}\\
&& + 2^{-kn}\ \Lambda^{-n}\ \Vert \lapn \varphi \Vert_{2,\R^n}\ \Vert \lapn f \Vert_{2,\R^n}\\
&& + 2^{-kn}\ \Lambda^{-n}\ \Vert \lapn \varphi \Vert_{2,\R^n}\ \Vert \lapn f \Vert_{2,\R^n}\\
&\aeq& II_{k,0} + 2^{-k\gamma}\ \Lambda^{-n}\  \Vert \lapn \varphi \Vert_{2,\R^n}\ \Vert \lapn f \Vert_{2,\R^n}.
\end{ma}
\]
It remains to estimate $II_{k,0}$:
\[
 \begin{ma}
  \abs{II_{k,0}} 
&\aleq& \Vert M_1 \lapms{s-t} \psi \Vert_{(\frac{2n}{n-2(s-t)},\infty),B_{\sqrt{\Lambda}r}}\ \Vert \lapn \varphi \Vert_{2,\R^n}\ \Vert \lapn f \Vert_{2}\\
&& +  \Vert M_1 \lapms{s-t} \psi \Vert_{(\frac{2n}{n-2(s-t)},\infty),B_{2\Lambda r}}\ \Vert M_2 \lapms{t} \abs{\lapn \varphi} \Vert_{(\frac{2n}{n-2t},2),\R^n \backslash B_{\sqrt{\Lambda}r}}\ \Vert \lapn f \Vert_{2}.
 \end{ma}
\]
Again, because $\dist(\supp \psi,B_{\sqrt{\Lambda}r}) \ageq 2^k \Lambda r$
\[
 \Vert M_1 \lapms{s-t} \psi \Vert_{(\frac{2n}{n-2(s-t)},\infty),B_{\sqrt{\Lambda}r}} \aleq{} \brac{2^k \Lambda}^{-\gamma} \Vert \psi \Vert_{(2,\infty),\R^n},
\]
\[
 \Vert M_1 \lapms{s-t} \psi \Vert_{(\frac{2n}{n-2(s-t)},\infty),B_{2\Lambda r}} \aleq{} \brac{2^k}^{-\gamma} \Vert \psi \Vert_{(2,\infty),\R^n},
\]
and Proposition~\ref{pr:lapmstlapnvarphibralrk} implies again (using that Lorentz spaces are normable and thus infinite triangular inequalities hold)
\[
\begin{ma}
  \Vert M_2 \lapms{t} \abs{\lapn \varphi} \Vert_{(\frac{2n}{n-2t},2),\R^n \backslash B_{\sqrt{\Lambda}r}}
&\aleq& \sum_{i = 1}^\infty \brac{2^i \sqrt {\Lambda}}^{-\gamma} \Vert \lapn \varphi \Vert_{2}\\
&\aleq& \sqrt {\Lambda}^{-\gamma} \Vert \lapn \varphi \Vert_{2}.\\
\end{ma}
\]
As \underline{for $I$},
\[
 H(\varphi,f) = H(\varphi,\lapmn \brac{\eta_{\Lambda^3 r}\lapn f}) + H(\varphi,\lapmn \brac{\brac{1-\eta_{\Lambda^3 r}}\lapn f}),
\]
and by the arguments in Remark \ref{rem:lowordest},
\[
 \Vert H(\varphi,\lapmn \brac{\eta_{\Lambda^3 r}\lapn f}) \Vert_{(2,1),\R^n} \aleq \Vert \lapn \varphi \Vert_{2,\R^n}\ \Vert \eta_{\Lambda^3 r} \lapn f \Vert_{2,\R^n}.
\]
It remains to estimate
\[
 \Vert H(\varphi,\lapmn \brac{\brac{1-\eta_{\Lambda^3 r}}\lapn f}) \Vert_{(2,1),B_{2\Lambda r}}.
\]
Again, this is done using Proposition~\ref{la:commlowerorder} and we have to control for some $\psi \in C_0^\infty(B_{2\Lambda r})$, $\Vert \psi \Vert_{(2,\infty),\R^n} \leq 1$,
\[
\begin{ma}
&&\int M_1 \lapms{s-t} \psi\ M_2 \lapms{t} \abs{\lapn \varphi}\ M_3 \lap^{-\frac{n}{4}+\frac{s}{2}}\abs{(1-\eta_{\Lambda^3 r})\lapn f}\\
&=&\int M_1 \lapms{s-t} \psi\ M_2 \lapms{t} \abs{\lapn \varphi}\ \eta_{\Lambda^2 r}\ M_3 \lap^{-\frac{n}{4}+\frac{s}{2}}\abs{(1-\eta_{\Lambda^3 r})\lapn f}\\
&&+ \int \brac{1-\eta_{\Lambda^2 r}}\ M_1 \lapms{s-t} \psi\ M_2 \lapms{t} \abs{\lapn \varphi}\ M_3 \lap^{-\frac{n}{4}+\frac{s}{2}}\abs{(1-\eta_{\Lambda^3 r})\lapn f}\\

% &\aleq& \Vrac{ M_1 \lapms{s-t} \psi}_{(\frac{2n}{n-2(s-t)},\infty),\R^n}\ \Vrac{ M_2 \lapms{t} \abs{\lapn \varphi}\ }_{(\frac{2n}{n-2t},2),\R^n}\ \Vrac{ M_3 \lap^{-\frac{n}{4}+\frac{s}{2}}\abs{(1-\eta_{\Lambda^3 r})\lapn f} }_{(\frac{n}{s},2), B_{\Lambda^2 r}}\\
% 
% && + \Vrac{ M_1 \lapms{s-t} \psi}_{(\frac{2n}{n-2(s-t)},\infty),\R^n \backslash B_{\Lambda^2 r}}\ \Vrac{ M_2 \lapms{t} \abs{\lapn \varphi}\ }_{(\frac{2n}{n-2t},2),\R^n \backslash B_{\Lambda^2 r}}\\
% &&\qquad \cdot \Vrac{ M_3 \lap^{-\frac{n}{4}+\frac{s}{2}}\abs{(1-\eta_{\Lambda^3 r})\lapn f} }_{(\frac{n}{s},2)}\\

&\aleq& \Vrac{ \psi}_{(2,\infty),\R^n}\ \Vrac{ \lapn \varphi }_{2,\R^n}\ \Vrac{ M_3 \lap^{-\frac{n}{4}+\frac{s}{2}}\abs{(1-\eta_{\Lambda^3 r})\lapn f} }_{(\frac{n}{s},2), B_{\Lambda^2 r}}\\

&& + \Vrac{ M_1 \lapms{s-t} \psi}_{(\frac{2n}{n-2(s-t)},\infty),\R^n \backslash B_{\Lambda^2 r}}\ \Vrac{ \lapn \varphi }_{2,\R^n}\ \Vrac{ \lapn f }_{2,\R^n}\\
\end{ma}
\]
As before, one has
\[
 \Vrac{ M_3 \lap^{-\frac{n}{4}+\frac{s}{2}}\abs{(1-\eta_{\Lambda^3 r})\lapn f} }_{(\frac{n}{s},2), B_{\Lambda^2 r}} \aleq  \Lambda^{-3s}\ \Lambda^{-2s}\ \Vrac{ \lapn f }_{2,\R^n},
\]
and
\[
 \Vrac{ M_1 \lapms{s-t} \psi}_{(\frac{2n}{n-2(s-t)},\infty),\R^n \backslash B_{\Lambda^2 r}} \aleq  \Lambda^{-\frac{n}{2}}\ \Vrac{ \psi }_{(2,\infty),\R^n}.
\]
So,
\[
 \abs{I} \aleq \Vrac{ \lapn \varphi }_{2,\R^n}\ \brac{\Vrac{ \eta_{\Lambda^3 r} \lapn f }_{2,\R^n} + \Lambda^{-\gamma} \Vrac{ \lapn f }_{2,\R^n}}.
\]
\end{proofL}

%% file: energy.tex
\section{Picking a good Frame: Improved Control by Energy Minimizing}\label{s:energy}
In this section, we prove that we can replace $\Omega \in L^2(so(N))$ by an $\Omega_P \in L^{2,1}_{loc}$ for an appropriate choice of $P$. 
\subsection{Adaption of H\'elein's Energy Method}
Let $D \subset \R^n$ be a smoothly bounded set. We define the energy functional
\[
E(Q) \equiv E_{D}(Q) := \intl_{\R^n} \abs{Q \lapn (Q^T-I) + Q \Omega Q^T}^2,\quad Q \in H^{\frac{n}{2}}_I(D,SO(N)).
\]
Here, similar to Definition \ref{def:Hs}, we have denoted for the identity matrix $I \in R^{N \times N}$
\[
 H^{\frac{n}{2}}_I(D,SO(N)) := \left \{ Q \in L^2(D):\ Q-I \in H^{\frac{n}{2}}(D,\R^{m\times m}),\ \supp(Q-I) \subset \overline{D}\right \}.
\]
We are going to prove the following two Lemmata which are adaptions of H\'elein's moving frame argument (see \cite{Hel02}, also \cite{Chone95}), and in their spirit similar to \cite[Lemma~2.2, Lemma~2.4]{IchEnergie}.
\begin{lemma}[Existence of a Minimizer]\label{la:en:ex}
Let $\Omega \in L^2(\R^n,\R^{N\times N})$. Then there exists $P \in H_I^{\frac{n}{2}}(\R^n,SO(N))$ such that $E(P) \leq E(Q)$ for any $Q \in H^{\frac{n}{2}}_I(D,SO(N))$. Moreover,
\[
\Vert \lapn (P-I) \Vert_{2,\R^n} \leq 2\ \Vert \Omega \Vert_{2,\R^n}.
\]
\end{lemma}

\begin{lemma}[Euler-Lagrange Equations]\label{la:en:el}
A critical point $P \in H^{\frac{n}{2}}_I(D,SO(N))$ of $E(\cdot)$ satisfies
\[
 \intl so(\Omega_P)\ \lapn \varphi = \intl so(H(\varphi,P-I) P^T \Omega_P) \quad \mbox{for all $\varphi \in C_0^\infty(D)$}.
\]
Here, $\Omega_P = P \lapn (P^T-I) + P \Omega P^T$ and $so(A) := \fracm{2} \brac{A - A^T}$.
\end{lemma}

\begin{proofL}{\ref{la:en:ex}}
Obviously, $Q \equiv I$ is a feasible mapping for $E(\cdot)$. Hence, we can assume the existence of a minimizing sequence $Q_k \in H_I^{\frac{n}{2}}(D,SO(N))$ such that
\[
 E(Q_k) \leq \Vert \Omega \Vert_{2,\R^n}^2.
\]
In particular
\begin{equation}\label{eq:en:uniflnbd}
 \Vert \lapn (Q_k-I) \Vert_{2,\R^n} \leq 2 \Vert \Omega \Vert_{2,\R^n}.
\end{equation}
We denote $R_k = Q_k - I \in H^{\frac{n}{2}}$. The mappings $R_k$ are uniformly bounded in $L^\infty(\R^n)$ because $Q_k \in SO(N)$ a.e. As $D$ is a bounded domain and $\supp R_k \subset \overline{D}$ we have a uniform $L^2(\R^n)$-bound of $R_k$ which together with \eqref{eq:en:uniflnbd} implies a uniform $H^{\frac{n}{2}}(\R^n,\R^{m\times m})$-bound for $R_k$. Consequently, we can choose a subsequence (again denoted with $R_k$) which converges weakly in $H^{\frac{n}{2}}(\R^n,\R^{m\times m})$ to some $R \in H^{\frac{n}{2}}(\R^n,\R^{m\times m})$.\\
Moreover, using the boundedness of $D \subset \R^n$ and Lemma~\ref{la:bs:compact}, up to taking yet again a subsequence, we can assume that $R_k \xrightarrow{k \to \infty} R$ strongly in $L^2(\R^n)$ and pointwise almost everywhere. Pointwise convergence implies that $P := R + I \in SO(N)$ almost everywhere, and thus $P \in H^{\frac{n}{2}}_I (D,SO(N))$.\\
Then
\[
 \begin{ma}
  E(Q_k) &=& \intl \abs{\lapn (Q_k^T-P^T) + \Omega (Q_k^T-P^T) + (\lapn P^T + \Omega P^T)}^2\\
  &=& \intl \abs{I_k + II_k + III}^2,\\
 \end{ma}
\]
where
\[
\begin{split}
&I_k := \lapn (Q_k^T-P^T),\\
&II_k := \Omega (Q_k^T-P^T),\\
&III :=\lapn P^T + \Omega P^T.
\end{split}
\]
We have that $II_k \xrightarrow{k \to \infty} 0$ almost everywhere by the pointwise convergence of $Q_k$. On the other hand, as $Q_k, P$ are bounded in $L^\infty$ and $\Omega \in L^2(\R^n)$, by Lebesgue's dominated convergence theorem $II_k \xrightarrow{k \to \infty} 0$ in $L^2(\R^n)$. This and the weak convergence of $\lapn (Q_k^T - P^T) \rightharpoonup 0$ in $L^2(\R^n)$ imply that the mixed terms
\[
\intl_{\R^n} I_k\ II_k,\quad \intl_{\R^n} I_k\ III,\quad  \intl_{\R^n} II_k\ III \xrightarrow{k \to \infty} 0.
\]
Furthermore,
\[
\intl_{\R^n} \abs{III}^2 = E(P).
\]
Consequently, for $k \to \infty$
\[
E(Q_k) = \Vert I_k \Vert_{2,\R^n}^2 + E(P) + o(1).
\]
Taking the limit $k \to \infty$ on both sides this implies
\[
\inf_Q E(Q) \geq \limsup_{k\to \infty} \Vert I_k \Vert_{2,\R^n}^2 + E(P).
\]
As $E(P) \geq \inf_Q E(Q)$ this implies $E(P) = \inf_Q E(Q)$ and $Q_k - P \xrightarrow{k \to \infty} 0$ in $H^{\frac{n}{2}}(\R^n,\R^{m\times m})$.
\end{proofL}

\begin{proofL}{\ref{la:en:el}}
Let $\varphi \in C_0^\infty(D)$, $\alpha \in so(N)$. We distort $P$ by
\[
Q_\varepsilon := e^{\varepsilon \varphi \alpha} P = P + \varepsilon \varphi\ \alpha\ P + o(\varepsilon) \in H^{\frac{n}{2}}_I(D,SO(N)).
\]
Then,
\[
Q_\varepsilon^T = P^T - \varepsilon \varphi\ P^T\ \alpha + o(\varepsilon),
\]
and 
\[
\begin{ma}
&&\lapn  (Q^T_\varepsilon -I)\\
&=& \lapn  (P^T-I) - \varepsilon \lapn  (\varphi P^T)\ \alpha + o(\varepsilon)\\
&=& \lapn  (P^T-I) - \varepsilon \lapn  (\varphi\ (P^T-I))\ \alpha - \varepsilon \lapn \varphi\ \alpha + o(\varepsilon)\\
&=& \lapn  (P^T-I) - \varepsilon\ \brac{\lapn  \varphi}\ (P^T-I)\ \alpha 
- \varepsilon \varphi\ \lapn (P^T-I)\ \alpha 
- \varepsilon\ H(\varphi, P^T-I)\ \alpha 
- \varepsilon \lapn \varphi\ \alpha + o(\varepsilon)\\
&=& \lapn  (P^T-I) - \varepsilon\ \brac{\lapn  \varphi}\ P^T\ \alpha 
- \varepsilon \varphi\ \lapn (P^T-I)\ \alpha 
- \varepsilon\ H(\varphi, P^T-I)\ \alpha 
+ o(\varepsilon).
\end{ma}
\]
We compute
\begin{equation}\label{eq:en:distqlq}
\begin{ma}
&&Q_\varepsilon\ \lapn (Q_\varepsilon^T-I)\\
&=& \brac{P + \varepsilon \varphi\ \alpha\ P}\ \brac{\lapn (P^T-I) - \varepsilon (\lapn \varphi)\ P^T \ \alpha - \varepsilon \varphi\ \lapn (P^T-I)\ \alpha - \varepsilon\ H(\varphi,P^T-I)\ \alpha} + o(\varepsilon)\\
&=& P\ \lapn (P^T-I) + \varepsilon \varphi \brac{\alpha\ P \lapn (P^T-I) - P \lapn (P^T-I)\ \alpha} - \varepsilon \lapn \varphi\ \alpha - \varepsilon\ P\ H(\varphi,P^T-I)\ \alpha + o(\varepsilon),
\end{ma}
\end{equation}
and
\begin{equation} \label{eq:en:distqoq}
Q_\varepsilon \Omega Q^T_\varepsilon = P \Omega P^T + \varepsilon \varphi\ \brac{\alpha\ P \Omega P^T - P \Omega P^T\ \alpha}+ o(\varepsilon).
\end{equation}
Recall that we denote the term $Q \lapn (Q^T-I) + Q\Omega Q^T$ by $\Omega_Q$. Then we infer from \eqref{eq:en:distqlq} and \eqref{eq:en:distqoq}
\[
\Omega_{Q_\varepsilon} = \Omega_P + \varepsilon \varphi\ \brac{\alpha\ \Omega_P - \Omega_P\ \alpha} - \varepsilon \lapn \varphi\ \alpha - \varepsilon\ P\ H(\varphi,P^T-I)\ \alpha + o(\varepsilon).
\]
In order to compute $\abs{\Omega_{Q_\varepsilon}}^2$ let us denote for $A \in \R^{N\times N}$, $B \in \R^{N\times N}$
\[
 A : B := A_{ij}B_{ij},
\]
thus in particular $\abs{A}^2 = A:A$. Note, that for any matrix $A \in \R^{N\times N}$ and any Matrix $B \in so(N)$,
\[
A:BA = B^TA:A = - BA:A = -A:BA,
\]
hence $A:BA = A:AB = 0$. Consequently,
\begin{equation}\label{eq:en:absoqe}
\abs{\Omega_{Q_\varepsilon}}^2 = \abs{\Omega_P}^2 + o(\varepsilon) - 2 \varepsilon\ \Omega_P: \brac{\lapn \varphi\ \alpha + P\ H(\varphi,P^T-I)\ \alpha}.
\end{equation}
Again we use a fact from Linear Algebra to continue: For any $A \in \R^{N\times N}$ and any $B \in so(N)$ we have
\[
A:B = so(A):B.
\]
Hence, \eqref{eq:en:absoqe} becomes
\[
\begin{ma}
\abs{\Omega_{Q_\varepsilon}}^2 &=& \abs{\Omega_P}^2 + o(\varepsilon) - 2 \varepsilon \brac{so(\Omega_P)\  \lapn \varphi + so(H(\varphi,P-I)\ P^T\ \Omega_P)}: \alpha.
\end{ma}
\]
We integrate this,
\[
\begin{ma}
E(Q_\varepsilon) - E(P)= o(\varepsilon) - 2\varepsilon \intl_{\R^n} so(\Omega_P): \alpha\ \lapn \varphi + so(H(\varphi,P-I)\ P^T\ \Omega_P):\alpha.
\end{ma}
\]
Dividing by $\varepsilon$ and taking the limit $\varepsilon \to 0$ we infer for a critical point $P$
\[
0 = \intl_{\R^n} so(\Omega_P): \alpha\ \lapn \varphi + so(H(\varphi,P-I)\ P^T\ \Omega_P):\alpha.
\]
This holds for every $\alpha \in so(N)$, so component-wise
\[
\intl_{\R^n} so(\Omega_P)\ \lapn \varphi = -\intl_{\R^n} so(H(\varphi,P-I)\ P^T\ \Omega_P).
\]
\end{proofL}

\subsection{Local integrability gain: Proof of Lemma~\ref{la:wlapneqhterm:nonlocalized}}
We want to show, that Euler-Lagrange systems as in Lemma~\ref{la:en:el} imply $L^{2,1}_{loc}$-integrability, because of the $H(\cdot,\cdot)$ appearing on the right-hand side. To this end, we are going to show in this section a localized version of Lemma~\ref{la:wlapneqhterm:nonlocalized}, more precisely we have
\begin{lemma}\label{la:wlapneqhterm}
There exists a constant $C > 0$ such that the following holds: Assume that $f,g,h \in L^2(\R^n)$, $\Lambda > 8$, and that for all $\varphi \in C_0^\infty(B_{\Lambda r})$
\begin{equation}\label{eq:wlapneqhterm:intlflapnveqrhs}
 \intl_{\R^n} f \lapn \varphi = \intl_{\R^n} g\ H(h,\varphi).
\end{equation}
Then,
\[ 
\begin{ma}
 \Vert f \Vert_{(2,1),B_r} &\leq& C\ \brac{\Vert g \Vert_{2,B_{\Lambda^3 r}} + \Lambda^{-\gamma} \Vert g \Vert_{2,\R^n}}\ \Vert \lapn h \Vert_{2,\R^n}\\ 
&&+\  C\ \Vert g \Vert_{2,\R^n}\ \brac{\Vert \lapn h \Vert_{2,B_{\Lambda^3 r}} + \Lambda^{-\gamma} \Vert \lapn h \Vert_{2,\R^n}} + \Lambda^{-\gamma} \Vert f \Vert_{2,\R^n}.
\end{ma}
\]
\end{lemma}

Before giving the proof, let us state several intermediate results.
\begin{proposition}\label{pr:floclapnest}
There is $\gamma > 0$ and $C > 0$ such that for any $B_r \subset \R^n$, $\varphi \in C_0^\infty(B_r)$, $f \in L^2(\R^n)$, and $\Lambda > 8$
\[
 \intl_{\R^n} f \varphi \leq \intl_{\R^n} f\ \lapn (\eta_{\Lambda r} \lapmn \varphi) + C\ \Lambda^{-\gamma}\ \Vert f \Vert_{2,\R^n}\ \Vert \varphi \Vert_{(2,\infty),\R^n}.
\]
\end{proposition}
\begin{proofP}{\ref{pr:floclapnest}}
We have
\[
\begin{ma}
&&\intl_{\R^n} f \lapn \brac{\brac{1-\eta_{\Lambda r}} \lapmn \varphi}\\
&\leq& \Vert f \Vert_{2,\R^n}\ \Vert H(\brac{1-\eta_{\Lambda r}},\lapmn \varphi) \Vert_{2,\R^n} + \Vert f \Vert_{2,\R^n}\ \Vert \lapn \brac{1-\eta_{\Lambda r}} \Vert_{4,\R^n}\ \Vert \lapmn \varphi\Vert_{4,\R^n}\\
&\overset{\sref{P}{pr:finswlapsinlp}}{\aleq}& \Vert f \Vert_{2,\R^n}\ \Vert H(\brac{1-\eta_{\Lambda r}},\lapmn \varphi) \Vert_{2,\R^n} + \Vert f \Vert_{2,\R^n}\ \brac{\Lambda r}^{-\frac{n}{4}}\ \Vert \varphi\Vert_{\frac{4}{3},\R^n}\\
&\aleq& \Vert f \Vert_{2,\R^n}\ \Vert H(\brac{1-\eta_{\Lambda r}},\lapmn \varphi) \Vert_{2,\R^n} + \Vert f \Vert_{2,\R^n}\ \brac{\Lambda r}^{-\frac{n}{4}}\ r^{\frac{n}{4}}\ \Vert \varphi\Vert_{(2,\infty),\R^n}
\end{ma}
\]
The result then follows by Proposition~\ref{pr:h1metalapmnvp}.
\end{proofP}

\begin{proposition}\label{pr:fHgerlapmnvpL1est}
For any $\varphi \in C_0^\infty(B_r)$, $\Lambda > 8$,
\[
\begin{ma}
\Vrac{ g\ H(h,\eta_{r} \lapmn \varphi) }_{1,\R^n} &\aleq& \brac{\Vert g \Vert_{2,B_{\Lambda^2 r}} + \Lambda^{-\gamma} \Vert g \Vert_{2,\R^n}}\ \Vert \lapn h \Vert_{2,\R^n}\ \Vert \varphi \Vert_{(2,\infty),\R^n}\\
&&+\Vert g \Vert_{2,\R^n}\ \brac{\Vert \lapn h \Vert_{2,B_{\Lambda^2 r}} + \Lambda^{-\gamma} \Vert \lapn h \Vert_{2,\R^n}}\ \Vert \varphi \Vert_{(2,\infty),\R^n}.
\end{ma}
\]
\end{proposition}
\begin{proofP}{\ref{pr:fHgerlapmnvpL1est}}
We have
\[
\begin{ma}
 g\ H(h,\eta_{r} \lapmn \varphi) 
&=& \eta_{\Lambda r}g\ H(h,\eta_{r} \lapmn \varphi) + \brac{1-\eta_{\Lambda r}}g\ H(h,\eta_{r} \lapmn \varphi)\\
&=:& I + II.
\end{ma}
\]
We use Remark \ref{rem:lowordest} in order to have
\[
\begin{ma}
 \Vert I \Vert_{1,\R^n} 
&\aleq& \Vert g \Vert_{2,B_{2\Lambda r}}\ \Vert \lapn h \Vert_{2,\R^n}\ \Vert \lapn \brac{\eta_r \lapmn \varphi} \Vert_{(2,\infty),\R^n}.
\end{ma}
\]
By Proposition~\ref{pr:estlapn1metalpmnvp},
\[
 \Vert \lapn \brac{\eta_r \lapmn \varphi} \Vert_{(2,\infty),\R^n} \aleq  \Vert \varphi \Vert_{(2,\infty)},
\]
and thus
\[
 \Vert I \Vert_{1,\R^n} \aleq \Vert g \Vert_{2,B_{2\Lambda r}}\ \Vert \lapn h \Vert_{2,\R^n}.
\]
As \underline{for $II$},
\[
 \Vert II \Vert_{1,\R^n} \aleq{} \Vert g \Vert_{2,\R^n}\ \Vert H(h,\eta_{r} \lapmn \varphi) \Vert_{2,\R^n\backslash B_{\Lambda r}}.
\]
In order to estimate
\[
\Vert H(h,\eta_{r} \lapmn \varphi) \Vert_{2,\R^n\backslash B_{\Lambda r}},
\]
by the lower order estimates in Lemma~\ref{la:commlowerorder} and a usual duality approach we have to estimate for some $\psi \in C_0^\infty(\R^n \backslash B_{\Lambda r})$, $\Vert \psi \Vert_{2,\R^n} \leq 1$, finitely many quantities like the following one, for some $s \in (0,\frac{n}{2})$, $t\in (0,s]$
\[
\begin{ma}
 &&\intl \brac{M_1 \lapms{s-t} \psi}\ \brac{M_2 \lapms{t} \abs{\lapn h}}\ \brac{M_3 \lap^{-\frac{n}{4}+s} \abs{\lapn \brac{\eta_r \lapmn \varphi}}}\\
&=&\intl \eta_{\frac{\Lambda}{2}r}\ \brac{M_1 \lapms{s-t} \psi}\ \brac{M_2 \lapms{t} \abs{\lapn h}} \ \brac{M_3 \lap^{-\frac{n}{4}+s} \abs{\lapn \brac{\eta_r \lapmn \varphi}}}\\
&&+\intl (1-\eta_{\frac{\Lambda}{2}r})\ \brac{M_1 \lapms{s-t} \psi}\ \brac{M_2 \lapms{t} \abs{\lapn h}} \ \brac{M_3 \lap^{-\frac{n}{4}+s} \abs{\lapn \brac{\eta_r \lapmn \varphi}}}\\
&=:& A_1 + A_2.
\end{ma}
\]
\underline{As for $A_1$}, by H\"older's inequality,
\[
 \abs{A_1} \aleq{} \Vert M_1 \lapms{s-t} \psi \Vert_{(\frac{2n}{n-2s+2t},2),\R^n}\ \Vert \eta_{\frac{\Lambda}{2}r} M_2 \lapms{t} \abs{\lapn h} \Vert_{(\frac{2n}{n-2t},2),\R^n}\ \Vert M_3 \lap^{-\frac{n}{4}+s} \abs{\lapn \brac{\eta_r \lapmn \varphi}} \Vert_{(\frac{n}{s},\infty),\R^n}.
\]
By Sobolev inequality Proposition~\ref{pr:sobineq},
\[
 \Vert M_1 \lapms{s-t} \psi \Vert_{(\frac{2n}{n-2s+2t},2),\R^n} \aleq{} \Vert \psi \Vert_{2,\R^n} \leq 1.
\]
By the localized Sobolev inequality, Lemma~\ref{la:locsobineq},
\[
\Vert \eta_{\frac{\Lambda}{2}r} M_2 \lapms{t} \abs{\lapn h} \Vert_{(\frac{2n}{n-2t},2),\R^n} \aleq \Vert \lapn h \Vert_{2,B_{\Lambda^2 r}} + \Lambda^{-\gamma} \Vert \lapn h \Vert_{2,\R^n}.
\]
Finally, by Proposition~\ref{pr:estlapn1metalpmnvp}
\[
 \Vert M_3 \lap^{-\frac{n}{4}+s} \abs{\lapn \brac{\eta_r \lapmn \varphi}} \Vert_{(\frac{n}{s},\infty),\R^n}
\aleq \Vert \lapn \brac{\eta_r \lapmn \varphi} \Vert_{(2,\infty),\R^n} \overset{\sref{P}{pr:estlapn1metalpmnvp}}{\aleq} \Vert \varphi \Vert_{(2,\infty),\R^n}.
\]
Consequently,
\[
 \abs{A_1} \aleq \brac{\Vert \lapn h \Vert_{2,B_{\Lambda^2 r}} + \Lambda^{-\gamma} \Vert \lapn h \Vert_{2,\R^n}}\ \Vert \varphi \Vert_{(2,\infty),\R^n}.
\]
It remains to \underline{estimate $A_2$}, and we have again by H\"older's inequality,
\[
\begin{ma}
 \abs{A_2}
%&\aleq& \Vert M_1 \lapms{s-t} \psi \Vert_{(\frac{2n}{n-2s+2t},2),\R^n}\ \Vert M_2 \lapms{t} \abs{\lapn h} \Vert_{(\frac{2n}{n-2t},2),\R^n}\ \Vert (1-\eta_{\frac{\Lambda}{2}r})\ M \laps{s} \brac{\eta_r \lapmn \varphi} \Vert_{(\frac{n}{s},\infty),\R^n}\\
&\overset{\sref{P}{pr:sobineq}}{\aleq}& \Vert \psi \Vert_{2,\R^n}\ \Vert \lapn h \Vert_{2,\R^n}\ \Vert M_3 \lap^{-\frac{n}{4}+s} \abs{\lapn \brac{\eta_r \lapmn \varphi}} \Vert_{(\frac{n}{s},\infty),\R^n\backslash B_{\frac{\Lambda}{2}r}}\\
&\aleq& \Vert \lapn h \Vert_{2,\R^n}\ \Vert M_3 \lap^{-\frac{n}{4}+s} \abs{\lapn \brac{\eta_r \lapmn \varphi}} \Vert_{(\frac{n}{s},\infty),\R^n\backslash B_{\frac{\Lambda}{2}r}}.
\end{ma}
\]
Next,
\[
\begin{ma}
 &&\Vert M_3 \lap^{-\frac{n}{4}+s} \abs{\lapn \brac{\eta_r \lapmn \varphi}} \Vert_{(\frac{n}{s},\infty),\R^n\backslash B_{\frac{\Lambda}{2}r}}\\
 &\aleq& \Vert M_3 \lap^{-\frac{n}{4}+s} \eta_{\sqrt{\Lambda} r}\abs{\lapn \brac{\eta_r \lapmn \varphi}} \Vert_{(\frac{n}{s},\infty),\R^n\backslash B_{\frac{\Lambda}{2}r}} + \Vert M_3 \lap^{-\frac{n}{4}+s} (1-\eta_{\sqrt{\Lambda} r})\abs{\lapn \brac{\eta_r \lapmn \varphi}} \Vert_{(\frac{n}{s},\infty),\R^n\backslash B_{\frac{\Lambda}{2}r}}\\
&\aleq& \Vert M_3 \lap^{-\frac{n}{4}+s} \eta_{\sqrt{\Lambda} r}\abs{\lapn \brac{\eta_r \lapmn \varphi}} \Vert_{(\frac{n}{s},\infty),\R^n\backslash B_{\frac{\Lambda}{2}r}} 
+ \Vert \lapn \brac{\eta_r \lapmn \varphi} \Vert_{(2,\infty),\R^n\backslash B_{\sqrt{\Lambda}r}}\\
&=:& A_{2,1} + A_{2,2}.
\end{ma}
\]
As for $A_{2,1}$ set $F := \abs{\lapn \brac{\eta_r \lapmn \varphi}}$, then as above by Proposition \ref{pr:estlapn1metalpmnvp}
\begin{equation}\label{eq:ener:Fleqvphi}
 \Vert F \Vert_{(2,\infty),\R^n} \aleq \Vert \varphi \Vert_{(2,\infty),\R^n}.
\end{equation}
We have for some $\phi \in C_0^\infty(\R^n \backslash B_{\frac{\Lambda}{2} r})$, $\Vert \phi \Vert_{\frac{n}{n-s},\R^n} \leq 1$,
\[
 \begin{ma}
  A_{2,1} &\aleq& \int \abs{\cdot}^{-\frac{n}{2}-s}\ \chi_{\abs{\cdot} \geq \Lambda r}\ \abs{\Phi}\ast \abs{\eta_{\sqrt{\Lambda} r}\ F}\\
&\aleq& \Vrac{\abs{\cdot}^{-\frac{n}{2}-s}\ \chi_{\abs{\cdot} \geq \Lambda r}}_{\frac{n}{s},\R^n}\ \Vrac{\Phi}_{\frac{n}{n-s},\R^n}\ \Vrac{\eta_{\sqrt{\Lambda} r}\ F}_{1,\R^n}\\
&\aleq& \brac{\Lambda r}^{-\frac{n}{2}}\ \Vrac{\Phi}_{\brac{\frac{n}{n-s},1},\R^n}\ \brac{\sqrt{\Lambda} r}^{\frac{n}{2}}\ \Vrac{F}_{(2,\infty),\R^n}\\
&\overset{\eqref{eq:ener:Fleqvphi}}{\aleq}& \Lambda^{-\frac{n}{4}}\ \Vrac{\varphi}_{(2,\infty),\R^n}.
\end{ma}
\]
As for $A_{2,2}$,
\[
 A_{2,2} \overset{\sref{P}{pr:lapnvarphibralrk}}{\aleq} \Lambda^{-\gamma}\ \Vert F \Vert_{(2,\infty),\R^n} \overset{\eqref{eq:ener:Fleqvphi}}{\aleq} \Lambda^{-\gamma}\ \Vert \varphi \Vert_{(2,\infty),\R^n}.
\]
\end{proofP}

Now we are able to give the
\begin{proofL}{\ref{la:wlapneqhterm}}
By Proposition~\ref{pr:floclapnest} (using Proposition~\ref{pr:lorentzfacts})
\[
\begin{ma}
\Vert f \Vert_{(2,1),B_r} 
&\overset{\sref{P}{pr:floclapnest}}{\aleq}& \sup_{\ontop{\varphi \in C_0^\infty(B_r)}{\Vert \varphi \Vert_{(2,\infty)} \leq 1}} \intl_{\R^n} f \lapn (\eta_{\Lambda r} \lapmn \varphi) + \Lambda^{-\gamma}\ \Vert f \Vert_{2,\R^n}\\
&\overset{\eqref{eq:wlapneqhterm:intlflapnveqrhs}}{=}& \sup_{\ontop{\varphi \in C_0^\infty(B_r)}{\Vert \varphi \Vert_{(2,\infty)} \leq 1}} \intl_{\R^n} g\ H(h, \eta_{\Lambda r} \lapmn \varphi) + \Lambda^{-\gamma}\ \Vert f \Vert_{2,\R^n}
\end{ma}
\]
The claim then follows from Proposition~\ref{pr:fHgerlapmnvpL1est}.
\end{proofL}

%% file: regproof.tex
\section{Preparations for Dirichlet Growth Theorem}
Let $v \in L^2(\R^n,\R^N)$ be a solution of \eqref{eq:pdelapnvomegav} in $D \subset \R^n$, i.e.
\begin{equation}\label{eq:preg:vpde}
 \intl_{\R^n} v^i\ \lapn \varphi = \intl_{\R^n} \Omega_{il} v^l\ \varphi \quad \mbox{for all $\varphi \in C_0^\infty(D)$}.
\end{equation}
% Let $u \in H^{\frac{n}{2}}(\R^n,\R^m)$ and $D \subset \R^n$ such that for some $\Omega \in L^2(\R^n,so(m))$ (\ToDo vgl. \eqref{eq:pdelapnvomegav})
% \begin{equation}\label{eq:preg:upde}
%  \intl_{\R^n} \lapn u^i\ \lapn \varphi = \intl_{\R^n} \Omega_{il} \lapn u^l\ \varphi \quad \mbox{for all $\varphi \in C_0^\infty(D)$}.
% \end{equation}
%\ToDo: Certain Distortions of $u$ with $L^\infty$-functions should be fine... (\ToDo auch im 2d-Riviere-Fall?)\\
Let $P-I \in H^{\frac{n}{2}}_0(D,\R^{m\times m})$, $P \in SO(m)$ almost everywhere. Set $w := P v$. Then $w \in L^2(\R^n)$ and moreover
\begin{proposition}\label{pr:Pvpisfeasibletf}
$\tilde{\varphi} := P \varphi$ for any $\varphi \in C_0^\infty(D)$ is a feasible test function for \eqref{eq:preg:vpde}.
\end{proposition}
\begin{proofP}{\ref{pr:Pvpisfeasibletf}}
Observe that $P \varphi \in L^2(\R^n)\cap L^\infty(\R^n)$, and
\[
 \lapn \tilde{\varphi} = P \lapn \varphi + \varphi \lapn (P-I) + H(\varphi,P-I) \in L^2(\R^n).
\]
Thus, for any small $\varepsilon > 0$ (depending only on the support of $\varphi$)
\[
\begin{ma}
 &&\int v\ \lapn \tilde{\varphi} - \int \Omega\ v\ \tilde{\varphi}\\
&=& \int v\ \lapn \brac{\tilde{\varphi} - \eta_\varepsilon \ast \tilde{\varphi}} + \int \Omega\ v\  \brac{\eta_\varepsilon \ast \tilde{\varphi}-\tilde{\varphi}}\\
&=:& I_\varepsilon + II_\varepsilon.
\end{ma}
\]
By Lebesgue's dominated convergence one sees that $II \xrightarrow{\varepsilon \to 0} 0$. As for $I$, pointwise almost everywhere, \[\lapn \brac{\tilde{\varphi} \ast \eta_{\varepsilon}} = \lapn \tilde{\varphi} \ast \eta_{\varepsilon}.\] Indeed, for any $\psi \in C_0^\infty(\R^n)$,
\[
\begin{ma}
 \int \lapn \brac{\tilde{\varphi} \ast \eta_{\varepsilon}}\ \psi
&=& \int \brac{\tilde{\varphi} \ast \eta_{\varepsilon}}\ \lapn \psi\\
&=& \int \int \tilde{\varphi}(x-y)\ \eta_{\varepsilon}(y)\ \lapn \psi(x)\ dx\ dy\\
&=& \int \int \brac{\lapn \tilde{\varphi}}(x-y)\ \eta_{\varepsilon}(y)\ \psi(x)\ dx\ dy\\
&=& \int \brac{\lapn \tilde{\varphi}}\ast \eta_{\varepsilon}\ \psi.
\end{ma}
\]
Since moreover, denoting $\mathcal{M}$ the Hardy-Littlewood Maximal function,
\[
 \abs{\lapn \tilde{\varphi} \ast \eta_{\varepsilon}} \leq \mathcal{M}\brac{\lapn \tilde{\varphi}} \in L^2(\R^n),
\]
again, Lebesgue's dominated convergence implies $I \xrightarrow{\varepsilon \to 0} 0$.
\end{proofP}%
Having taken care of this, we start computing
\[
 \intl_{\R^n} w\ \lapn \varphi = \intl_{\R^n} \lapn (\varphi P)\ v + \intl_{\R^n} \brac{-\lapn (\varphi P) + (\lapn \varphi)\ P} \ v.
\]
As for the second part,
\[
\begin{ma}
&& \intl_{\R^n} \brac{-\lapn (\varphi P) + (\lapn \varphi)\ P + \varphi \lapn P - \varphi \lapn P}\  v\\
&=& \intl_{\R^n} \brac{-\lapn (\varphi (P-I)) + (\lapn \varphi)\ (P-I) + \varphi \lapn (P-I) - \varphi \lapn (P-I)} \ v\\
&=& \intl_{\R^n} \brac{-H(\varphi,P-I) - \varphi \lapn (P-I)} \ v.
\end{ma}
\]
Consequently, using Proposition~\ref{pr:Pvpisfeasibletf},
\[
\begin{ma}
\intl_{\R^n} w\ \lapn \varphi &=& \brac{\intl_{\R^n} \lapn (\varphi P_{ji})\ v^i + \intl_{\R^n} \brac{-H(\varphi,P-I)_{ji} - \varphi \lapn (P-I)_{ji}}\  v^i}_j\\
&\overset{\eqref{eq:preg:vpde}}{=}& \brac{\intl_{\R^n} \varphi P_{ji} \Omega_{il} v^l + \intl_{\R^n} \brac{-H(\varphi,P-I)_{ji} - \varphi \lapn (P-I)_{ji}} \ v^i}_j\\
&=& \intl_{\R^n}  \brac{P \Omega - \lapn (P-I)}\ v\ \varphi - \intl_{\R^n} H(\varphi,P-I) \ v\\
&=& \intl_{\R^n}  \brac{P \Omega P^T - \brac{\lapn (P-I)} P^T}\ w\ \varphi - \intl_{\R^n} H(\varphi,P-I) P^T\ w\\
&=& \intl_{\R^n}  \brac{P \Omega P^T - \brac{\lapn (P-I)} P^T}\ w\ \varphi - \intl_{\R^n} H(\varphi,P-I) P^T\ w.
\end{ma}
\]
Now, we have
\[
 \begin{ma}
  &&H(P-I,P^T-I)\\
&=& \lapn \brac{(P-I)(P^T-I)} - \brac{\lapn (P-I)}(P^T-I) - (P-I) \brac{\lapn (P^T-I)}\\
&=& \lapn \brac{2I - P^T - P} - \brac{\lapn P}(P^T-I) - (P-I) \brac{\lapn P^T}\\
&=& -\lapn P^T - \lapn P  - \brac{\lapn P} P^T + \lapn P - P\lapn P^T + \lapn P^T\\
&=& - \brac{\brac{\lapn (P-I)} P^T  + P\lapn (P^T-I)}.
\end{ma}
\]
Then,
\[
 \begin{ma}
\brac{\lapn (P-I)} P^T
&=& \fracm{2} \brac{\lapn (P-I)} P^T - \fracm{2} P \lapn (P^T-I) + \fracm{2} \brac{\lapn (P-I)} P^T + \fracm{2} P \lapn (P^T-I)\\
&=& so\brac{\lapn (P-I) P^T} - \fracm{2} H(P-I,P^T-I).
 \end{ma}
\]
Thus, for any $\varphi \in C_0^\infty(D)$
\begin{equation}\label{eq:pdew}
\begin{ma}
\intl_{\R^n} w\ \lapn \varphi &=& \intl_{\R^n}  so\brac{\Omega_P}\ w\ \varphi +  \intl_{\R^n} \brac{\fracm{2}H(P-I,P^T-I)\ \varphi - H(\varphi,P-I) P^T}\ w.
\end{ma}
\end{equation}

\subsection{Estimates on the left-hand side}
We have
\[
\Vert w \Vert_{(2,\infty),B_r} \overset{\sref{P}{pr:lorentzfacts}}{=} \sup_{\ontop{g \in C_0^\infty(B_r)}{\Vert g \Vert_{(2,1)} \leq 1}} \int w\ g.
\]
For $g \in C_0^\infty(B_r)$, $\lapmn g \in L^{p}$ for any $p \in (2,\infty]$ for $\Lambda > 2$
\[
\lapmn g = \eta_{\Lambda r} \lapmn g + (1-\eta_{\Lambda r})\lapmn g.
\]
Note, that $\lapn$ is well-defined on both parts with estimates by Proposition~\ref{pr:estlapn1metalpmnvp}, so
\[
g = \lapn \brac{\eta_{\Lambda r} \lapmn g} + \lapn \brac{(1-\eta_{\Lambda r})\lapmn g},
\]
and both parts are in $L^2(\R^n)$. Then,
\[
 \int w\ g = \int w\ \lapn \brac{\eta_{\Lambda r} \lapmn g} + \int w\ \lapn \brac{(1-\eta_{\Lambda r}) \lapmn g}
\]
As for the second term, it is 'nicely' controlled:
\begin{lemma}\label{la:cutoffstufflhs}
There are constants $\gamma > 0$, $C > 0$ such that if $\Vert g \Vert_{(2,1)} \leq 1$, $\supp g \subset \overline{B_r}$, $\Lambda \geq 8$,
\[
 \int w\ \lapn \brac{(1-\eta_{\Lambda r}) \lapmn g} \leq C\  \Lambda^{-\gamma} \Vert w \Vert_{(2,\infty), B_{4\Lambda r}} + C\ \Lambda^{-\gamma }\ \sum_{k \in \N} 2^{-\gamma k} \Vert w \Vert_{(2,\infty),A_{4\Lambda r,k}}
\]
\end{lemma}
\begin{proofL}{\ref{la:cutoffstufflhs}}
We have
\[
 \begin{ma}
  &&\int w\ \lapn \brac{(1-\eta_{\Lambda r}) \lapmn g}\\
&=&\int \eta_{2\Lambda r} w\ \lapn \brac{(1-\eta_{\Lambda r}) \lapmn g} + \int (1-\eta_{2\Lambda r})\ w\ \lapn \brac{(1-\eta_{\Lambda r}) \lapmn g}\\
&=:& I + II.
 \end{ma}
\]
As for $I$, by Proposition~\ref{pr:estlapn1metalpmnvp},
\[
 \abs{I} \aleq{} \Vert \eta_{2\Lambda r} w \Vert_{(2,\infty),\R^n}\ \Lambda^{-\gamma}\ \Vert g \Vert_{(2,1),\R^n} \leq \Lambda^{-\gamma}\ \Vert \eta_{2\Lambda r} w \Vert_{(2,\infty),\R^n}.
\]
As for $II$,
\[
\begin{ma}
II &=& \sum_{k=1}^\infty \int \eta_{2\Lambda r}^k\ w\ \lapn \brac{(1-\eta_{\Lambda r}) \lapmn g}\\
&=& \sum_{k=1}^\infty \int \eta_{2\Lambda r}^k\ w\ \brac{\lapn \brac{1-\eta_{\Lambda r}}}\ \lapmn g +\sum_{k=1}^\infty \int \eta_{2\Lambda r}^k\ w\ H(1-\eta_{\Lambda r},\lapmn g)\\
&=& \sum_{k=1}^\infty \int \eta_{2\Lambda r}^k\ w\ \brac{\lapn \brac{1-\eta_{\Lambda r}}}\ \lapmn g +\sum_{k=1}^\infty \int \eta_{2\Lambda r}^k\ w\ H(-\eta_{\Lambda r},\lapmn g)\\
&=:& \sum_{k=1}^\infty \brac{II_{1,k} + II_{2,k}}.
\end{ma}
\]
As for $II_{1,k}$,
\[
\begin{ma}
 II_{1,k} &\aleq{}& \brac{2^k \Lambda r}^{-\frac{n}{2}}\ \Vert g \Vert_{1}\ \Vert \eta_{2\Lambda r}^k\ w\ \brac{\lapn \brac{1-\eta_{\Lambda r}}} \Vert_{1}\\
&\aleq{}& \brac{2^k \Lambda}^{-\frac{n}{2}}\ \Vert g \Vert_{(2,1),\R^n}\ \Vert \eta_{2\Lambda r}^k\ w\Vert_{(2,\infty),\R^n}\ \Vert \lapn (1-\eta_{\Lambda r}) \Vert_{(2,1),\R^n}\\
&\aleq& \brac{2^k \Lambda}^{-\frac{n}{2}} \Vert \eta_{2\Lambda r}^k\ w \Vert_{(2,\infty),\R^n}\ \Vert \lapn \eta_{\Lambda r} \Vert_{(2,1),\R^n}\\
&\overset{\sref{P}{pr:finswlapsinlp}}{\aleq}& \brac{2^k \Lambda}^{-\frac{n}{2}} \Vert \eta_{2\Lambda r}^k\ w \Vert_{(2,\infty),\R^n}.
 \end{ma}
\]
% Now, if $k > 2$ for some $\psi \in C_0^\infty(A_{\Lambda r,k})$, $\Vert \psi \Vert_{(2,\infty),\R^n} \leq 1$, cf. Proposition~\ref{pr:lorentzfacts},
% \[
% \begin{ma}
%  &&\Vert \lapn \eta_{\Lambda r} \Vert_{(2,1),A_{\Lambda r,k}}\\
% &\aleq& \int \psi\ \lapn \eta_{\Lambda r}\\
% &\aeq& \int \abs{\xi}^{-\frac{3}{2}n}\ \chi_{\abs{\xi} > 2^k\Lambda r}\ \abs{\psi} \ast \abs{\eta_{\Lambda r}}(\xi)\ d\xi\\
% &\aleq& \brac{2^k \Lambda r}^{-\frac{3}{2}n}\ \Vert \psi \Vert_{1,\R^n}\ \Vert \eta_{\Lambda r} \Vert_{1,\R^n}\\
% &\aleq& \brac{2^k \Lambda r}^{-\frac{3}{2}n}\ \brac{2^k\Lambda}^{\frac{n}{2}}\ \Vert \psi \Vert_{(2,\infty),\R^n}\ \brac{\Lambda r}^n\\
% &\leq& 2^{-nk}.\\
% \end{ma}
% \]
% If $k = 1,2$ we just estimate
% \[
%  \Vert \lapn \eta_{\Lambda r} \Vert_{(2,1),A_{\Lambda r,k}} \aleq \brac{\Lambda r}^{-\frac{n}{2}}\  \brac{2^k \Lambda r}^{\frac{n}{2}} \overset{k \leq 2}{\aleq} 2^{-nk}.
% \]
In order to estimate \underline{$II_{k,2}$}, we have consider (finitely many) terms of the following form for $s \in (0,\frac{n}{2})$, $t \in (0,s]$
\[
\begin{ma}
 &&\Vert \abs{\eta_{2\Lambda r}^k\ w}\ M_1 \lapms{s-t} \brac{M_2 \lapms{t} \abs{\lapn (1-\eta_{\Lambda r})}\ M_3 \lap^{-\frac{n}{4}+\frac{s}{2}}\abs{g}}\ \Vert_{1,\R^n}\\
&\aleq& \Vert \brac{\eta_{2\Lambda r}^k\ w} \Vert_{(2,\infty),\R^n}\ \Vert M_1 \lapms{s-t} \brac{M_2 \lapms{t} \abs{\lapn \eta_{\Lambda r}}\ M_3 \lap^{-\frac{n}{4}+\frac{s}{2}}\abs{g}}\ \Vert_{(2,1),A_{k,2\Lambda r}}.\\
\end{ma}
\]
For some $\psi \in C_0^\infty(A_{k,2\Lambda r})$, $\Vert \psi \Vert_{(2,\infty)} \leq 1$ we have
\[
\begin{ma}
 &&\Vert M_1 \lapms{s-t} \brac{M_2 \lapms{t} \abs{\lapn (1-\eta_{\Lambda r})}\ M_3 \lap^{-\frac{n}{4}+\frac{s}{2}}\abs{g}}\ \Vert_{(2,1),A_{k,\Lambda r}}\\
&\aleq& \int \brac{M_1 \lapms{s-t}\psi}\ \brac{M_2 \lapms{t} \abs{\lapn (1-\eta_{\Lambda r})}}\ \brac{M_3 \lap^{-\frac{n}{4}+\frac{s}{2}}\abs{g}}\\
&=:& A_0 + \sum_{l=1}^\infty \int (\eta_{\frac{\Lambda}{4} r}^l)\ \brac{M_1 \lapms{s-t}\psi}\ \brac{M_2 \lapms{t} \abs{\lapn (1-\eta_{\Lambda r})}}\ \brac{M_3 \lap^{-\frac{n}{4}+\frac{s}{2}}\abs{g}}\\
&\aleq& \abs{A_0}+ \sum_{l=1}^\infty \Vert M_1 \lapms{s-t} \psi \Vert_{(\frac{2n}{n-2s+2t},\infty),A_{l,\frac{\Lambda}{4} r}}\
\Vert M_2\lapms{t} \abs{\lapn \eta_{\Lambda r}} \Vert_{\frac{2n}{n-2t},A^l_{\Lambda r}}\
\Vert M_3 \lap^{-\frac{n}{4}+\frac{s}{2}} \abs{g}  \Vert_{(\frac{n}{s},1),A^l_{\Lambda r}}\\
&\overset{\sref{P}{pr:lapmstlapnvarphibralrk}}{\aleq}& \abs{A_0}+ \sum_{l=1}^\infty 2^{-ln}\ \Vert M_1 \lapms{s-t} \psi \Vert_{(\frac{2n}{n-2s+2t},\infty),A_{l,\frac{\Lambda}{4} r}}\ \Vert M_3 \lap^{-\frac{n}{4}+\frac{s}{2}} \abs{g}  \Vert_{(\frac{n}{s},1),A^l_{\Lambda r}}\
\end{ma}
\]
As usual, if $l \geq 1$,
\[
 \Vert M_3 \lap^{-\frac{n}{4}+\frac{s}{2}} \abs{g}  \Vert_{(\frac{n}{s},1),A^l_{\Lambda r}} \aleq \brac{2^l \Lambda}^{-\frac{n}{2}}\ \Vert g \Vert_{(2,1),\R^n} \leq \brac{2^l \Lambda}^{-\frac{n}{2}}.
\]
Moreover, by Proposition \ref{pr:lapmsvarphiakal},
\[
\begin{ma}
 &&\sum_{l=1}^\infty 2^{-ln}\ \Vert M_1 \lapms{s-t} \psi \Vert_{(\frac{2n}{n-2s+2t},\infty),A_{l,\frac{\Lambda}{4} r}}\\
&\aleq& 2^{-kn} + \sum_{l=1}^k 2^{k(-n+s-t)+k\brac{n-\frac{n-2s+2t}{2}-s+t}+l\frac{n-2s+2t}{2}-nl} \\
&&+ \sum_{l=k}^\infty 2^{l(-n+s-t)+k\brac{n-\frac{n-2s+2t}{2}-s+t}+l\frac{n-2s+2t}{2}-nl}\\
&\aleq& 2^{-kn} + 2^{k\brac{-\frac{n}{2}+s-t}} \sum_{l=1}^k 2^{l\brac{-\frac{n}{2}-s+t}} \\
&&+ 2^{k\frac{n}{2}} \sum_{l=k}^\infty 2^{-\frac{3}{2}nl}\\
&\aleq& 2^{-\gamma k}.
\end{ma}
\]
%
%
%
% \ToDo If $\abs{l-k+2} \leq 1$, we have
% \[
% \begin{ma}
%  &&\abs{A_l}\\
% &\aleq& \Vert M_1 \lapms{s-t} \psi \Vert_{(\frac{2n}{n-2s+2t},\infty),A_{l,\frac{\Lambda}{4} r}}\
% \Vert \lapms{t} \abs{\lapn \eta_{\Lambda r}} \Vert_{\frac{2n}{n+2t},\R^n}\
% \Vert M_3 \lap^{-\frac{n}{4}+\frac{s}{2}} \abs{g} \Vert_{(\frac{n}{s},1),A_{l,\frac{\Lambda}{4}r}}\\
% &\leq& \Vert M_1 \lapms{s-t} \psi \Vert_{(\frac{2n}{n-2s+2t},\infty),A^l_{\frac{\Lambda}{4}r}}\
% \Vert M_3 \lap^{-\frac{n}{4}+\frac{s}{2}} \abs{g}\Vert_{(\frac{n}{s},1),A^l_{\frac{\Lambda}{4} r}}\\
%
%
% \end{ma}
% \]
% Now for any $l \in \N_0$, if $\abs{k-l} \geq 4$ (\ToDo check!)
% \[
% \begin{ma}
% &&\Vert M_1 \lapms{s-t} \psi \Vert_{(\frac{2n}{n-2s+2t},\infty),A^l_{\Lambda r}}\\
% &\aleq& 2^{\max\{k,l\}(-n+s-t)+k\frac{n}{2} + l (\frac{n}{2}-s+t)}\ \Vert \psi \Vert_{(2,\infty),\R^n}\\
% &\leq& 2^{\max\{k,l\}(-n+s-t)+k\frac{n}{2} + l (\frac{n}{2}-s+t)},
% \end{ma}
% \]
% if $k \approx l$, we use
% \[
% \begin{ma}
%  &&\Vert M_1 \lapms{s-t} \psi \Vert_{(\frac{2n}{n-2s+2t},\infty),A^l_{\Lambda r}}\\
% &\aleq& \Vert \psi \Vert_{(2,\infty),\R^n} \leq 1
% \end{ma}
% \]
% and as usual for any $l \geq 1$,
% \[
% \begin{ma}
% &&\Vert M_3 \lap^{-\frac{n}{4}+\frac{s}{2}} \abs{g}  \Vert_{(\frac{n}{s},1),A^l_{\Lambda r}}\\
% &\aleq& 2^{-l\frac{n}{2}}\ \Lambda^{-\frac{n}{2}}\ \Vert g \Vert_{(2,1)}\\
% &\leq& 2^{-l\frac{n}{2}} \Lambda^{-\frac{n}{2}}.
% \end{ma}
% \]
For $l = 0$ we argue as in the proof of Lemma~\ref{la:estHvpfgL1}:
\[
\begin{ma}
\abs{A_0} &\aleq& \Vert M_1 \lapms{s-t} \psi \Vert_{\brac{\frac{2n}{n-2s+2t},\infty},B_{\sqrt{\Lambda}r}}  +
\Vert M_1 \lapms{s-t} \psi \Vert_{\brac{\frac{2n}{n-2s+2t},\infty},B_{2\Lambda r}}\ \Vert M_3 \lap^{-\frac{n}{4}+\frac{s}{2}} g \Vert_{\brac{\frac{n}{s},1},B_{\sqrt{\Lambda}r}}\\
&\aleq& 2^{-k\gamma}\ \sqrt{\Lambda}^{-\gamma}.
\end{ma}
\]
% Thus, we have
% \[
% \begin{ma}
%  &&\sum_{l=0}^\infty \abs{A_l}\\
% &\aleq& \Lambda^{-\gamma} \sum_{\abs{l-k} \leq 4} 2^{\max\{k,l\}(-n+s-t)+k\frac{n}{2} + l (-s+t)} + \Lambda^{-\gamma} \sum_{l \approx k} 2^{-l\frac{n}{2}}\\
% &\aleq& \Lambda^{-\gamma} \sum_{l = 0}^k 2^{k(-n+s-t)+k\frac{n}{2} + l (-s+t)} + \Lambda^{-\gamma} \sum_{l = k}^\infty 2^{l(-n+s-t)+k\frac{n}{2} + l (-s+t)} +  2^{-k\frac{n}{2}}\\
% &=& 2^{k(-\frac{n}{2}+s-t)}\ \Lambda^{-\gamma}\sum_{l = 0}^k 2^{l (-s+t)} + 2^{k\frac{n}{2}}\ \Lambda^\gamma\ \sum_{l = k}^\infty 2^{-ln} +  2^{-k\frac{n}{2}}\\
% &=& 2^{k(-\frac{n}{2}+s-t)} \Lambda^{-\gamma}\ \sum_{l = -k}^0 2^{(l+k) (-s+t)} + 2^{k\frac{n}{2}}\ \sum_{l = 0}^\infty 2^{-(l+k)n} +  2^{-k\frac{n}{2}}\\
% &\aeq{}& 2^{-k\frac{n}{2}} \Lambda^{-\gamma}.
% \end{ma}
% \]
We conclude that
\[
 \abs{II_{k,2}} \aleq{} \Lambda^{-\gamma}\ 2^{-k\frac{n}{2}} \Vert w \Vert_{(2,\infty),A^k_{\Lambda r}}.
\]
\end{proofL}

\subsection{Estimates on the right-hand side}
\begin{lemma}[Estimates on the right-hand side]\label{la:est:rhs}
Let $w \in L^2(\R^n)$ be a solution to \eqref{eq:pdew}, where $P -I \in H^{\frac{n}{2}}(\R^n,\R^{m\times m})$, $P \in SO(m)$ a.e., which is a minimizer of $E(\cdot)$ defined in Lemma~\ref{la:en:ex}. Then there exists constants $C_w$, $\gamma > 0$, $\Lambda_0 > 0$ and for any $\Lambda > \Lambda_0$ an $R \in (0,1)$ such that for any $\varphi \in C_0^\infty(B_r)$, if $B_{r} \subset D$, $r \in (0,R)$, and $\Vert \lapn \varphi \Vert_{(2,1),\R^n} \leq 1$
\[
\abs{\intl_{\R^n} w\ \lapn \varphi} \leq C_w\ \Lambda^{-\gamma} \Vert w \Vert_{(2,\infty),B_{\Lambda r}} + C_w \Lambda^{-\gamma} \sum_{k=1}^\infty 2^{-k\gamma}\ \Vert w \Vert_{(2,\infty),A_{\Lambda r,k}}.
\]
\end{lemma}
\begin{proofL}{\ref{la:est:rhs}}
Let $\tilde{\gamma}$ be the smallest of the various exponents of $\Lambda^{-1}$ and $2^{-k}$ appearing (see below). Pick first $\Lambda_0 > 0$ such that for all $\Lambda > \Lambda_0$
\[
 \Lambda^{-\frac{\tilde{\gamma}}{4}}\ \Vert \lapn (P-I) \Vert_{2,\R^n} \leq 1 \tag{$\Lambda1$}
\]
\[
 \Lambda^{-\frac{\tilde{\gamma}}{4}}\ \Vert \Omega \Vert_{2,\R^n} \leq 1\tag{$\Lambda2$}
\]
Next, for fixed $\Lambda > \Lambda_0$ pick $R \in (0,1)$ such that for all $B_{r} \subset \R^n$, $r \in (0,R)$
\[
 \Vert H(P-I,P^T - I) \Vert_{(2,1),B_{\Lambda^5 r}} \leq \Lambda^{-\tilde{\gamma}} \tag{R1}
\]
\[
 \Vert \lapn (P - I) \Vert_{2,B_{\Lambda^5 r}} \leq \Lambda^{-\tilde{\gamma}} \tag{R2}
\]
\[
 \Vert \Omega_P \Vert_{2,B_{\Lambda^5 r}} \leq \Lambda^{-\tilde{\gamma}} \tag{R3}
\]
Let $\varphi \in C_0^\infty(B_r)$, $\Vert \lapn \varphi \Vert_{(2,1),B_r} \leq 1$. We denote the three parts to be estimated \[
\begin{ma}
I \quad &:=& \intl_{\R^n}  so\brac{\Omega}_P\ w\ \varphi\\
II &:=& \intl_{\R^n} H(P-I,P^T-I)\ \varphi\ w\\
III &:=& \intl_{\R^n} H(\varphi,P-I) P^T\ w.
\end{ma}
\]
There is a constant $C$ depending only on the dimensions involved, such that
\[
 \abs{II} \overset{(R1)}{\leq} \Lambda^{-\tilde{\gamma}}\ \Vert \varphi \Vert_{\infty,\R^n}\ \Vert w \Vert_{(2,\infty),B_r} \overset{\sref{L}{la:soblinfty}}{\leq}\ C\ \Lambda^{-\tilde{\gamma}}\ \Vert w \Vert_{(2,\infty),B_r}.
\]
By Lemma~\ref{la:estHvpfgL1},
\[
 \begin{ma}
 \abs{III} &\leq& C\ \brac{\Vert \lapn (P-I) \Vert_{2,B_{\Lambda^3 r}} + \Lambda^{-\tilde{\gamma}} \Vert \lapn (P-I) \Vert_{2,\R^n}}\ \Vert w \Vert_{(2,\infty),B_{\Lambda r}}\\
&&+ C\ \Lambda^{-\tilde{\gamma}}\ \sum_{k=1}^\infty 2^{-\tilde{\gamma} k}\ \Vert \lapn (P-I) \Vert_{2,\R^n}\ \Vert w \Vert_{(2,\infty),A_{k,\Lambda r}}\\
&\overset{(R2)}{\leq}& 2C\ \Lambda^{-\frac{\tilde{\gamma}}{2}}\ \Vert w \Vert_{(2,\infty),B_{\Lambda r}}+ C\ \Lambda^{-\frac{\tilde{\gamma}}{2}}\ \sum_{k=1}^\infty 2^{-\tilde{\gamma} k}\ \ \Vert w \Vert_{(2,\infty),A_{k,\Lambda r}}.
\end{ma}
\]
As for $I$, by the choice of $P$ and Lemma~\ref{la:en:el}, we have
\[
 \int so(\Omega_P)\ \lapn \varphi = \frac{1}{2} \int so(H(\varphi,P-I) P^T \Omega_P) \quad \mbox{for all $\varphi \in C_0^\infty(D)$}.
\]
Lemma~\ref{la:wlapneqhterm} implies
\[
\begin{ma}
 \Vert so(\Omega_P) \Vert_{(2,1),B_r} &\aleq&  \brac{\Vert \Omega_P \Vert_{2,B_{\Lambda^3 r}} + \Lambda^{-\tilde{\gamma}} \Vert \Omega_P \Vert_{2,\R^n}}\\
&&+ \Vert \Omega_P \Vert_{2,\R^n}\ \brac{\Vert \lapn (P-I)\Vert_{2,B_{\Lambda^3 r}} + \Lambda^{-\tilde{\gamma}} \Vert \lapn (P-I)\Vert_{2,\R^n}}\\
&& + \Lambda^{-\tilde{\gamma}} \brac{\Vert \lapn P \Vert_{2,\R^n} + \Vert \Omega \Vert_{2,\R^n}}\\
&\leq& C\ \Lambda^{-\frac{\tilde{\gamma}}{4}}.
\end{ma}
\]
Thus,
\[
 \abs{I} \leq C\ \Lambda^{-\frac{\tilde{\gamma}}{4}}\ \Vert \varphi \Vert_{\infty}\ \Vert w \Vert_{(2,\infty),B_r} \leq C\ \Lambda^{-\frac{\tilde{\gamma}}{4}}\ \Vert w \Vert_{(2,\infty),B_r}.
\]
Now we set $\gamma := \frac{\tilde{\gamma}}{4}$, and the claim is proven.
\end{proofL}

\subsection{Controlled Local Behavior: Proof of Lemma~\ref{la:localcontroldPreIt}}
Lemma~\ref{la:cutoffstufflhs} and Lemma~\ref{la:est:rhs} imply\footnote{Here we use also that
\[
\begin{ma}
&&\Lambda^{-\gamma} \Vert w \Vert_{(2,\infty),B_{\Lambda r}} + \ \Lambda^{-\gamma} \sum_{k=1}^\infty 2^{-k\gamma}\ \Vert w \Vert_{(2,\infty),A_{\Lambda r,k}}\\
&\aleq& \Lambda^{-\gamma}\log \Lambda\ \Vert w \Vert_{(2,\infty),B_{\Lambda^2 r}} + \ \Lambda^{-\gamma}\log \Lambda\ \sum_{k=1}^\infty 2^{-k\gamma}\ \Vert w \Vert_{(2,\infty),A_{\Lambda^2 r,k}}
 \end{ma}
\]
and note that $\Lambda^{-\gamma} \log \Lambda \leq \Lambda^{-\frac{\gamma}{2}}$ for sufficiently large $\Lambda$.}:
\begin{lemma}[Estimates on the left-hand side]\label{la:est:lhs}
Let $w \in L^2(\R^n)$ be a solution to \eqref{eq:pdew}, where $P\in H_I^{\frac{n}{2}}(D,\R^{m\times m})$, $P \in SO(m)$ a.e., which is a minimizer of $E(\cdot)$ defined in Lemma~\ref{la:en:ex}. Then there exists $\Lambda_0 > 0$,$\gamma > 0$, $C_w > 0$ such that for any $\Lambda > \Lambda_0$ there is an $R \in (0,1)$  such that if $B_{\Lambda r} \subset D$, $r \in (0,R)$
\[
\Vert w \Vert_{(2,\infty),B_r} \leq C_w\ \Lambda^{-\gamma} \Vert w \Vert_{(2,\infty),B_{\Lambda r}} + C_w\ \Lambda^{-\gamma} \sum_{k=1}^\infty 2^{-k\gamma}\ \Vert w \Vert_{(2,\infty),A_{\Lambda r,k}}.
\]
\end{lemma}
In particular, for $K := \log_2 \Lambda \in \N$ for some $\Lambda > \Lambda_0$, $q := 2^{-\gamma}$, $\varepsilon := \Lambda^{-\gamma}$ the condition \eqref{eq:it:smallness} of the Iteration Lemma~\ref{la:it} is satisfied, where
\[
\Phi(\lambda) := \Vert w \Vert_{(2,\infty),B_{\lambda \Lambda^{-1} R}}
\]
\[
\psi(\lambda) := \Vert w \Vert_{(2,\infty),B_{\lambda \Lambda^{-1} R} \backslash B_{\fracm{2}\lambda \Lambda^{-1} R}}.
\]
Note that also
\[
 \sum_{k =1}^\infty q^k \psi(\lambda) \aleq \brac{\sum_{k =1}^\infty \Vert w \Vert_{2,B_{\lambda \Lambda^{-1} R} \backslash B_{\fracm{2}\lambda \Lambda^{-1} R}}^2}^{\frac{1}{2}} \aleq \Vert w \Vert_{2}.
\]
As a consequence, we have shown Theorem~\ref{th:regthm}; More precisely, we have
\begin{theorem}\label{th:regthm2}
Let $v \in L^2(\R^n)$ be a solution to \eqref{eq:pdelapnvomegav}, i.e.
\[
 \lapn v = \Omega v \quad \mbox{in $D$}.
\]
Then for certain constants $\Lambda > 0$, $R \in (0,1)$, $C > 0$ depending all on $v$ we have that
\[
 \Vert v \Vert_{(2,\infty),B_r} \leq C\ r
\]
whenever $B_{\Lambda r} \subset D$, and $r \in (0,R)$.
\end{theorem}
In particular, we have the following
\begin{corollary}\label{co:locmorreyreg}
Let $v \in L^2(\R^n)$ be a solution to \eqref{eq:pdelapnvomegav}. Then for certain constants $\Lambda > 0$, $R \in (0,1)$, $C > 0$ depending all on $v$ we have that whenever $B_{\Lambda r} \subset D$, and $r \in (0,R)$,
\[
 \sup_{B_t \subset \R^n} t^{-1}\ \Vert \chi_{B_r} v \Vert_{(2,\infty),B_t} \leq C.
\]
\end{corollary}

%% file: boundary.tex
\section{Continuity Estimates: Proof of Theorem~\ref{th:regthmGen}}
% The goal of this section is the proof of Theorem \
% \begin{theorem}\label{th:regthmbd}
% Let $u \in H^{\frac{n}{2}}(\R^n)$, $v := A\ M \lapn u$ for some zero-multiplier operator $M$ and $A \in L^\infty (\R^n) \cap H^{\frac{n}{2}}(\R^n)$ and $D \subset \R^n$ be a bounded domain, $\partial D \in C^\infty$. Assume that there are constants constants $\Lambda > 0$, $R \in (0,1)$, $C > 0$ depending all on $v$ we have that
% \[
%  \Vert v \Vert_{(2,\infty),B_r} \leq C\ r^2
% \]
% whenever $B_{\Lambda r} \subset D$, and $r \in (0,R)$. If $u \in C^0(\R^n \backslash D)$, then $u \in C^0(\R^n) \cap C^{(0,\alpha)}$ for $\alpha = \ToDo$.
% \end{theorem}
In Theorem \ref{th:regthm2} we have shown, that $v$ satisfies local Morrey space estimates in $D$. In this section, we will show how regularity for $\lapmn v$ in the interior and on the boundary follows.

\subsection{In the Interior}
First of all, we note that local behavior of $v$ essentially controls local H\"older continuity of $\lapmn v$. More precisely,
\begin{proposition}\label{pr:locholdcontlocw}
Let $v \in L^2(\R^n)$, then for any $\Lambda \geq 2$, $x,y \in B_r(\overline{x})$
\[
 \abs{\lapmn v(x) - \lapmn v(y)} \aleq{} \frac{\abs{x-y}}{\Lambda r}\ \Vert v \Vert_{2,\infty} + \abs{\lapmn \brac{\chi_{B_{\Lambda r}(\overline{x})} v}(x) - \lapmn\brac{\chi_{B_{\Lambda r}(\overline{x})} v}(y)}.
\]
In particular, if $\lapmn\brac{\chi_{B_{\Lambda r}(\overline{x})} v}$ is H\"older-continuous in $B_{r}(\overline{x})$, so is $\lapmn v(x)$.
\end{proposition}
\begin{proofP}{\ref{pr:locholdcontlocw}}
We have
\[
\begin{ma}
 &&\lapmn v(x) - \lapmn v(y)\\
&\leq& \sum_{k=0}^\infty \intl_{\abs{\xi-\overline{x}} \in \brac{2^k \Lambda r,2^{k+1} \Lambda r}} \abs{\abs{\xi-x}^{-\frac{n}{2}}- \abs{\xi-y}^{-\frac{n}{2}}}\ \abs{v(\xi)}\ d\xi + \abs{\lapmn \brac{\chi_{B_{\Lambda r}(\overline{x})} v}(x) - \lapmn\brac{\chi_{B_{\Lambda r}(\overline{x})}}(y)}\\
&=:& \sum_{k=0}^\infty I_k + \abs{\lapmn \brac{\chi_{B_{\Lambda r}(\overline{x})} v}(x) - \lapmn\brac{\chi_{B_{\Lambda r}(\overline{x})}}(y)}.
\end{ma}
\]
If $\abs{\xi-\overline{x}} \in \brac{2^k \Lambda r,2^{k+1} \Lambda r}$ and $x \in B_r(\overline{x})$,
\[
 \abs{\xi-x} \geq \abs{\xi-x_0}-\abs{x_0-x} \geq 2^k \Lambda r - r \overset{\Lambda \geq 2}{\geq} 2^{k-1} \Lambda r,
\]
so if $x,y \in B_r(\overline{x})$ and $\abs{\xi-\overline{x}} \in \brac{2^k \Lambda r,2^{k+1} \Lambda r}$,
\[
\abs{\abs{\xi-x}^{-\frac{n}{2}}- \abs{\xi-y}^{-\frac{n}{2}}} \overset{\sref{P}{pr:sillyest2}}{\aleq} \brac{2^{k} \Lambda r}^{-\frac{n}{2}-1} \abs{x-y}.
\]
Hence,
\[
 \begin{ma}
I_k &\leq& \brac{2^{k} \Lambda r}^{-1}\ \abs{x-y}\ \Vert v \Vert_{(2,\infty),\R^n},
 \end{ma}
\]
which implies that
\[
 \sum_{k=0}^\infty I_k \leq \frac{\abs{x-y}}{\Lambda r}\ \Vert v \Vert_{(2,\infty),\R^n}.
\]
\end{proofP}

The following is a consequence of \cite[Proposition 3.3.]{Adams75} and the corollary afterwards.
\begin{proposition}\label{pr:ialphamorreyhoelder}
For any $\alpha \in (0,1)$ there exists $\beta > 0$ and a constant $C_\alpha$ such that the following holds. If for some $K > 0$
\begin{equation}\label{eq:ialphamh:morreyv}
 \sup_{x \in \R^n} \sup_{B_t(x)} t^{-\alpha} \Vert v \Vert_{(2,\infty),B_t(x)} \leq K,
\end{equation}
then for any $D \subsubset \R^n$ there exists a constant $C_{D,\alpha}$ such that
\[
 \sup_{x\neq y} \frac{\lapmn v(x) - \lapmn v(y)}{\abs{x-y}^\beta} \leq C_{D,\alpha}\ K.
\]
\end{proposition}
For the convenience of the reader, we sketch (for $n \geq 3$) the
\begin{proofP}{\ref{pr:ialphamorreyhoelder}}
This can be done very similarly to the proof of Proposition~\ref{pr:locholdcontlocw}: It suffices to show that $\lapmn v$ belongs to a Morrey-Campanato space, cf. \cite[Chapter III]{GiaquintaMI83}, that is
\[
\sup_{B_r(x)} r^{-n-\alpha} \intl_{B_r(x)} \abs{\lapmn v-(\lapmn v)_{B_r(x)}} \aleq K.
\]
% By the arguments in \cite[Proposition 3.3.]{Adams75}, there is $s \in (0,\frac{n}{2})$, such that
% \[
%  \sup_{B_r(x)} r^{-n} \int \abs{\lapms{s}v-(\lapms{s}v)_{B_r(x)}} \aleq K.
% \]
% \ToDo (needed?)
We have
\[
\begin{ma}
 &&\intl_{B_r(x)} \abs{\lapmn v-(\lapmn v)_{B_r(x)}}\\
&\aleq& r^{-n} \intl_{B_r(x)}\intl_{B_r(x)} \intl_{\R^n} \abs{\abs{z_1-\xi}^{-\frac{n}{2}}-\abs{z_2-\xi}^{-\frac{n}{2}}}\ \abs{v(\xi)}\ d\xi\ dz_1\ dz_2\\
&\aleq& r^{-n} \intl_{B_r(x)}\intl_{B_r(x)} \intl_{B_{2r}(x)} \abs{\abs{z_1-\xi}^{-\frac{n}{2}}-\abs{z_2-\xi}^{-\frac{n}{2}}}\ \abs{v(\xi)}\ d\xi\ dz_1\ dz_2\\
&&+ \sup_{z_1,z_2 \in B_r(x)}\ r^{n} \sum_{k=1}^\infty \intl_{\abs{\xi-x}\in (2^k,2^{k+1})} \abs{\abs{z_1-\xi}^{-\frac{n}{2}}-\abs{z_2-\xi}^{-\frac{n}{2}}}\ \abs{v(\xi)}\ d\xi\\
&=:& r^{-n} I + \sup_{z_1,z_2 \in B_r(x)}\ r^{n} \sum_{k=1}^\infty II_k.
\end{ma}
\]
As for $II_k$, by virtually the same arguments as in the proof of Proposition~\ref{pr:locholdcontlocw}
\[
 II_k \overset{\eqref{eq:ialphamh:morreyv}}{\leq} C_\alpha\ r^{-1+\alpha}\ \abs{z_1-z_2}\ K \leq C_\alpha\ r^{\alpha}\ K.
\]
As for $I_k$,
\[
 r^{-n} I_k \aleq r^n\ \Vert v \Vert_{(2,\infty),B_4r} \overset{\eqref{eq:ialphamh:morreyv}}{\leq} r^{n+\alpha}\ K.
\]
\end{proofP}
This and Corollary \ref{co:locmorreyreg} imply the interior H\"older continuity of Theorem~\ref{th:regthmGen}.
\subsection{On the Boundary}
We adapt the famous technique by Hildebrandt and Kaul, \cite{HildKaul72}, in order to obtain boundary regularity. A crucial part of this is, that the $BMO$-Norm is small on small sets.
\begin{lemma}[Local BMO-estimate]\label{la:bmohn2loc}
Let $u \in L^2(\R^n)$, $\lapn u \in L^2(\R^n)$. Set
\[
 \mathcal{M}(\lambda) := \sup_{x \in \R^n}\ \sup_{\tilde{\lambda} \in (0,\lambda)} {\tilde{\lambda}}^{-n} \intl_{B_{\tilde{\lambda}}(x)} \abs{u(z)-(u)_{B_{\tilde{\lambda}}(x)}}\ dz.
\]
Here,
\[
(u)_{B_{\tilde{\lambda}}(x)} :=  \abs{B_{\tilde{\lambda}}(x)}^{-1}\ \intl_{B_{\tilde{\lambda}}(x)} u.
\]
Then
\[
 \mathcal{M}(\lambda) \xrightarrow{\lambda \to 0} 0.
\]
\end{lemma}
\begin{proofL}{\ref{la:bmohn2loc}}
If $n = 1$, one just checks that
\[
 \mathcal{M}(\lambda) \aleq \sup_{x \in \R^n}\ \sup_{\tilde{\lambda} \in (0,\lambda)} \brac{\intl_{B_\lambda(x)} \intl_{B_\lambda(x)} \frac{\abs{u(z_1)-u(z_2)}^2}{\abs{z_1-z_2}^{2}}\ dz\ dz_2 }^{\frac{1}{2}},
\]
which tends to zero as $\lambda \to 0$, cf. \cite{NHarmS10Arxiv}. So let from now on $n \geq 2$. We adapt an approach in \cite[Proposition 3.3]{Adams75}. Let $\Lambda \geq 2$ and decompose
\[
\begin{ma}
 u &=& \lapmn \lapn u\\
&=& \lapmn \brac{\chi_{B_{\Lambda \lambda}(x)} \lapn u} + \lapmn \brac{\chi_{\R^n \backslash B_{\Lambda \lambda}(x)} \lapn u}\\
&=:& u_1 + u_2.
\end{ma}
\]
Consequently,
\[
\begin{ma}
&&\intl_{B_\lambda(x)} \abs{u(z_1)-(u)_{B_\lambda(x)}}\ dz_1 \\
&\aleq& 2 \intl_{B_\lambda(x)} \abs{u_1(z)}\ dz + \sum_{k=0}^\infty \lambda^n\ \sup_{z_1,z_2 \in B_\lambda(x)} \intl_{2^k \Lambda \lambda \leq  \abs{\xi-x} \leq 2^{k+1} \Lambda \lambda} \abs{\abs{z_1-\xi}^{-\frac{n}{2}}-\abs{z_2-\xi}^{-\frac{n}{2}}}\ \abs{\lapn u(\xi)}\ d\xi\\
&=:& 2 I + \sum_{k=0}^\infty \lambda^n\ \sup_{z_1,z_2 \in B_\lambda(x)} II_k.
\end{ma}
\]
As for $I$ we have,
\[
\begin{ma}
  \intl_{B_\lambda(x)} \abs{u_1(z)}\ dz
&\aleq& \intl_{B_\lambda(x)} \intl_{\abs{\xi-x} \leq \Lambda \lambda} \abs{z-\xi}^{-\frac{n}{2}}\ \abs{\lapn u(\xi)}\ d\xi\ dz\\
&\aeq& \brac{\Lambda \lambda}^{\frac{n}{2}}\ \intl_{B_{\Lambda \lambda}(x)} \abs{\lapn u(\xi)}\ d\xi\\
&\aleq& \Lambda^{n}\ \lambda^{n}\ \Vert \lapn u \Vert_{2,B_{\Lambda \lambda}(x)}.\\
\end{ma}
\]
As for $II_k$, by Proposition~\ref{pr:sillyest2},
\[
\begin{ma}
II_k &\aleq& \brac{2^k \Lambda \lambda}^{-\frac{n}{2}-1}\ \abs{z_1 - z_2}\ \ \brac{2^k \Lambda \lambda}^{\frac{n}{2}} \Vert \lapn u \Vert_{2,\R^n}\\
&\aleq& 2^{-k}\ \Lambda^{-1}\ \Vert \lapn u \Vert_{2,\R^n}.
\end{ma}
\]
Hence,
\[
\lambda^{-n}\ \intl_{B_\lambda(x)} \abs{u(z_1)-(u)_{B_\lambda(x)}}\ dz_1 \aleq  \Lambda^{n}\ \Vert \lapn u \Vert_{2,B_{\Lambda \lambda}(x)} +  \Lambda^{-1}\ \Vert \lapn u \Vert_{2,\R^n}.
\]
In particular,
\[
 \mathcal{M}(\lambda) \aleq \Lambda^{n}\ \sup_{x \in \R^n}\ \Vert \lapn u \Vert_{2,B_{\Lambda \lambda}(x)} +  \Lambda^{-1}\ \Vert \lapn u \Vert_{2,\R^n}.
\]
For any $\varepsilon > 0$ we choose then $\Lambda > 2$ such that $\Lambda^{-1}\ \Vert \lapn u \Vert_{2,\R^n} \leq \frac{\varepsilon}{2}$. Afterwards we pick $\lambda_0 > 0$ such that
\[
 \Lambda^{\frac{n}{2}}\ \sup_{x \in \R^n}\ \Vert \lapn u \Vert_{2,B_{\Lambda \lambda_0}(x)} \leq \frac{\varepsilon}{2}.
\]
Then, for any $\lambda \in (0,\lambda_0)$,
\[
 \mathcal{M}(\lambda) \aleq \varepsilon.
\]
\end{proofL}

Then, we have the following theorem, which implies the boundary regularity claimed in Theorem~\ref{th:regthmGen}.

\begin{theorem}\label{th:HKfrac}
Let $u \in L^2(\R^n)$, $\lapn u \in L^2(\R^n)$. Assume that $D \subsubset \R^n$, $\partial D \in C^{\infty}$ and that for some $\Lambda \geq 1$, $K \geq 0$
\begin{equation}\label{eq:bd:intcalpha}
\abs{u(x)-u(y)} \leq K\ \brac{\frac{\abs{x-y}^\alpha}{r^\alpha} + \abs{x-y}^\alpha} \quad \mbox{for almost all $x,y \in B_r$ where $B_{\Lambda r} \subset D$}.
\end{equation}
If there exists $\gamma > 0$ such that
\begin{equation}\label{eq:bd:uregonbound}
 u \in C^{0}\brac{B_\gamma(D)\backslash D},
\end{equation}
then $u \in C^{0}\brac{B_\gamma(D)} \cap C^{0,\alpha}(B)$.
\end{theorem}
\begin{proofT}{\ref{th:HKfrac}}
Let $\varepsilon > 0$ be given, and let $\delta \in (0,1)$, $\mu \in (0,\Lambda^{-1})$ to be chosen later depending on $\varepsilon$. Take $\mathcal{M}(\lambda)$ from Lemma~\ref{la:bmohn2loc}. Because of $u \in L^2(\R^n)$, $\lapn u \in L^2(\R^n)$ we have that $\mathcal{M}(\lambda) \xrightarrow{\lambda \to 0} 0$.\\
Since $\partial D$ is a closed, smooth manifold we assume w.l.o.g. that $\gamma > 0$ is small enough so that there exists the nearest-point projection $\Pi \in C^\infty\brac{B_\gamma(\partial D), \partial D}$, cf. \cite{Simon96}. Denote the mirroring function at $\partial D$ by $\psi: B_\gamma(\partial D) \cap \overline{D} \to B_\gamma(D) \backslash D$,
\[
 \psi(x) := 2\Pi(x)-x.
\]
Fix $\overline{x} \in D$, $\overline{z} \in \partial D$, and assume that $\abs{\overline{x}-\overline{z}} \leq \delta < \fracm{2}\gamma$. Set $\overline{y} := \Pi(\overline{x}) \in \partial D$. As $\Pi$ is well defined around $\overline{x}$ this implies as well $\tau := \abs{\overline{x}-\overline{y}} \leq \delta$ and since $\Pi$ is the (unique) nearest point projection into the boundary $\partial D$, we know that $B_{\tau}(\overline{x}) \subset D$. Denote $\sigma := \mu \tau$. Then,
\[
\begin{ma}
 \intl_{B_\sigma(\overline{x})} \abs{u(x)-u(\psi(x))} &\leq& \intl_{B_\sigma(\overline{x})} \abs{u(x)-(u)_{B_{2\tau}(\overline{y})}}\ dx + \intl_{B_\sigma(\overline{x})} \abs{u(\psi(x))-(u)_{B_{2\tau}(\overline{y})}}\ dx\\
&\aleq{}& \tau^n \mathcal{M}(2\tau) + \intl_{B_\sigma(\psi^{-1}(\overline{x}))} \abs{u(x)-(u)_{B_\tau(\overline{y})}}\ dx\\
&\aleq& \tau^n \mathcal{M}(2\tau).
\end{ma}
\]
% Next,
% \[
%  \begin{ma}
%   I_1 &\aleq& \tau^{-n} \intl_{B_\tau(\overline{y})}\ \intl_{B_\sigma(\overline{x})} \abs{u(x)-u(z)}\ dx\\
% &\aleq& \tau^{-\frac{n}{2}}\ \sigma^{\frac{n}{2}}\ \brac{\intl_{B_\tau(\overline{y})}\ \intl_{B_\sigma(\overline{x})} \abs{u(x)-u(z)}^2\ dx}^{\frac{1}{2}}\\
% &\aleq& \tau^{\frac{n}{2}+\frac{s}{2}}\ \mathcal{M}(2\tau).
% \end{ma}
% \]
Here we used the following computation of the behavior of the transformation $x \to \psi(x)$: For $z \in \R^n$
\[
d\psi_x[z] :=\frac{d}{dt}\big \vert_{t=0} \psi(x+tz) = 2d\Pi_x[z]-z.
\]
One can show then (again, cf. \cite{Simon96}) that $d\Pi_x[\cdot]$ is an orthogonal projection of $\R^n$ onto $T_x\partial D$, thus if we take an orthonormal basis $o_1,\ldots,\ o_n \in \R^n$ of $\R^n$ where $o_n \perp T_x\partial D$, we have that
\[
 d\psi_x[o_i] = 2o_i -o_i = o_i \quad \mbox{if $1 \leq i \leq n-1$}
\]
and
\[
 d\psi_x[o_n] = 0-o_n = -o_n \quad \mbox{if $1 \leq i \leq n-1$}.
\]
Thus,
\[
 \intl_{B_\sigma(\overline{x})} \abs{u(x)-u(\psi(x))} \aleq \tau^{n}\ \mathcal{M}(2\tau).
\]
By a contradiction argument this implies that there exists a set of positive measure $E_\sigma \subset B_\sigma(\overline{x})$ such that for any $x' \in E_\sigma$.
\[
 \abs{u(x')-u(\psi(x'))} \aleq \brac{\frac{\tau}{\sigma}}^{n}\ \mathcal{M}(2\tau).
\]
This and \eqref{eq:bd:intcalpha} for $r := \Lambda^{-1}\tau$ imply (using that $B_{\tau}(\overline{x}) \subset D$)
\[
\begin{ma}
 \abs{u(\overline{x})-u(\overline{y})} &\leq&\abs{u(\overline{x})-u(x')} + \abs{u(x')-u(\psi(x'))} + \abs{u(\psi(x'))-u(\overline{y})}\\
&\aleq& \brac{\Lambda \frac{\sigma}{\tau}}^\alpha + \brac{\frac{\tau}{\sigma}}^{n}\ \mathcal{M}(2\tau) + \abs{u(\psi(x'))-u(\overline{y})}.
\end{ma}
\]
Recalling that $\sigma := \mu \tau$ we have shown
\[
\abs{u(\overline{x})-u(\overline{y})} \leq C_{\Lambda,K,\alpha}\ \mu^\alpha \brac{1 + \mu^{\alpha-n} \mathcal{M}(2\tau)} + \abs{u(\psi(x'))-u(\overline{y})}.
\]
Picking $\mu \ll 1$ depending only on $\Lambda$, $K$ and $\alpha$ such that
\[
 C_{\Lambda,K}\ \mu^\alpha \leq \frac{\varepsilon}{4},
\]
and then $\delta \ll 1$ only depending on $\mu$ and $\mathcal{M}$ such that
\[
\mu^{\alpha-n} \mathcal{M}(\tau) \leq \mu^{\alpha-n} \mathcal{M}(\delta) \leq 1,
\]
as well as for any $x' \in B_\frac{\delta}{2}(x) \cap D$
\[
 \abs{u(\psi(x'))-u(\overline{y})} \leq \frac{\varepsilon}{2},
\]
we have shown,
\[
 \abs{u(\overline{x})-u(\overline{y})} \aleq \varepsilon.
\]
As $\overline{z} \in D \cap B_\delta(\overline{y})$, this implies the claim, as $\delta$ was chosen uniformly for $\overline{x}$ and $\overline{z}$.
\end{proofT}

%% file: iteration.tex
\renewcommand{\thesection}{A}
\renewcommand{\thesubsection}{A.\arabic{subsection}}
\section{Iteration Arguments}
We give the proof of the iteration argument which is used in \cite{DR1dMan}, \cite{DndMan}. Note, that an argument as in \cite{DR1dSphere} is not viable for our setting here, cf. Remark \ref{rem:nowidman}.

\begin{lemma}\label{la:it}
Let $q \in (0,1)$, $K \in \N$, $\varepsilon > 0$ and assume (say)
\begin{equation}\label{eq:it:smallness}
 \varepsilon+q^{2K}+2K\ \varepsilon\ q^K + \varepsilon \brac{1 + \fracm{1-q}+ \ 2K\ \fracm{1-q}} \leq \fracm{4}.
\end{equation}
Let moreover $\Phi: (0,\infty) \to (0,\infty)$ be monotone rising, $\psi: (0,\infty) \to (0,\infty)$ such that for all $\lambda \in (0,\infty)$
\begin{equation}\label{eq:it:psiPhi}
 \psi(\lambda) \leq \Phi(\lambda).
\end{equation}
Assume that for all $\lambda \in (0,1]$
\begin{equation}\label{eq:it:start}
 \Phi\brac{2^{-K}\lambda} \leq \varepsilon\ \Phi(\lambda) + \varepsilon \sum_{k=0}^\infty q^k\ \psi\brac{2^k\lambda}.
\end{equation}
If there is $G < \infty$ so that for all $\lambda \in (0,1)$
\[
 \sum_{k=0}^\infty q^k\ \psi(2^k\lambda) \leq G.
\]
Then, for all $\lambda \in (0,1)$,
\[
 \Phi(\lambda) \leq 32\ \lambda \brac{\Phi(\infty) + G}.
\]
\end{lemma}
\begin{proofL}{\ref{la:it}}
Let $l \in \N$, $l \geq 2$. We claim that for all $j \in \{0,\ldots,\left \lfloor \frac{l-1}{2} \right \rfloor \}$
\begin{equation}\label{eq:it:IV}
 \Phi\brac{2^{-Kl}} \leq \brac{\fracm{4}}^j \Phi\brac{2^{-K(l-2j)}} + \brac{\fracm{4}}^j\ \sum_{k=0}^\infty q^k\ \psi\brac{2^{-K(l-2j+1)+k}}.
\end{equation}
This is true for $j = 0$ by the monotonicity of $\Phi$. Assume now \eqref{eq:it:IV} to be true for some $j$, then in order to show \eqref{eq:it:IV} also for $j+1$, we estimate
\[
\begin{ma}
&&\Phi\brac{2^{-Kl}}\\
&\overset{\eqref{eq:it:IV}}{\leq}& \brac{\fracm{4}}^j \Phi\brac{2^{-K(l-2j)}} + \brac{\fracm{4}}^j\ \sum_{k=0}^\infty q^k\ \psi\brac{2^{-K(l-2j+1)+k}}\\ 
&\overset{\eqref{eq:it:start}}{\leq}& \brac{\fracm{4}}^j \brac{\varepsilon\ \Phi\brac{2^{-K(l-2j-1)}} + \varepsilon \sum_{k=0}^\infty q^k \psi\brac{2^{-K(l-2j-1)+k}}} + \brac{\fracm{4}}^j\ \sum_{k=0}^\infty q^k\ \psi\brac{2^{-K(l-2j+1)+k}}\\

&=& \brac{\fracm{4}}^j \varepsilon\ \Phi\brac{2^{-K(l-2j-1)}} + \brac{\fracm{4}}^j \varepsilon \sum_{k=0}^\infty q^k \psi\brac{2^{-K(l-2j-1)+k}}\\
&&+ \brac{\fracm{4}}^j\ \sum_{k=2K}^\infty q^k\ \psi\brac{2^{-K(l-2j-1)+k-2K}} + \brac{\fracm{4}}^j\ \sum_{k=0}^{2K-1} q^k\ \psi\brac{2^{-K(l-2j-1)+k-2K}}\\

&\overset{\eqref{eq:it:psiPhi}}{\leq}& \brac{\fracm{4}}^j \varepsilon\ \Phi\brac{2^{-K(l-2(j+1))}} + \brac{\fracm{4}}^j \varepsilon \sum_{k=0}^\infty q^k \psi\brac{2^{-K(l-2j-1)+k}}\\
&&+ \brac{\fracm{4}}^j\ q^{2K} \sum_{k=0}^\infty q^k\ \psi\brac{2^{-K(l-2j-1)+k}} + \brac{\fracm{4}}^j\ \sum_{k=0}^{2K-1} q^k\ \Phi\brac{2^{-K(l-2j-1)+k-2K}}\\

&=& \brac{\fracm{4}}^j \varepsilon\ \Phi\brac{2^{-K(l-2(j+1))}}  + \brac{\fracm{4}}^j \brac{\varepsilon+q^{2K}} \sum_{k=0}^\infty q^k \psi\brac{2^{-K(l-2(j+1)+1)+k}}\\
&& + \brac{\fracm{4}}^j\ \sum_{k=0}^{2K-1} q^k\ \Phi\brac{2^{-K(l-2j-1)+k-2K}}\\
\end{ma}
\]
\[
\begin{ma}
&\overset{\eqref{eq:it:start}}{\leq}& \brac{\fracm{4}}^j \varepsilon\ \Phi\brac{2^{-K(l-2(j+1))}} + \brac{\fracm{4}}^j \brac{\varepsilon+q^{2K}} \sum_{k=0}^\infty q^k \psi\brac{2^{-K(l-2(j+1)+1)+k}}\\
&&+ \brac{\fracm{4}}^j\ \varepsilon \sum_{k=0}^{2K-1} q^k\ \Phi\brac{2^{-K(l-2(j+1))+k-2K}} + \brac{\fracm{4}}^j\ \varepsilon \sum_{k=0}^{2K-1} \sum_{i=0}^\infty q^{k+i}\ \psi\brac{2^{-K(l-2j-1)+k+i-K}}\\

&\leq& \brac{\fracm{4}}^j \varepsilon\ \Phi\brac{2^{-K(l-2(j+1))}} + \brac{\fracm{4}}^j \brac{\varepsilon+q^{2K}} \sum_{k=0}^\infty q^k \psi\brac{2^{-K(l-2(j+1)+1)+k}}\\
&&+ \brac{\fracm{4}}^j\ \varepsilon\ \brac{\sum_{k=0}^{\infty} q^k}\ \Phi\brac{2^{-K(l-2(j+1))}} + \brac{\fracm{4}}^j\ \varepsilon \ 2K\ \sum_{i=0}^\infty q^{i}\ \psi\brac{2^{-K(l-2j-1)+i-K}}\\

&=& \brac{\fracm{4}}^j \varepsilon \brac{1 + \fracm{1-q}}\ \Phi\brac{2^{-K(l-2(j+1))}} 
+ \brac{\fracm{4}}^j \brac{\varepsilon+q^{2K}} \sum_{k=0}^\infty q^k \psi\brac{2^{-K(l-2(j+1)+1)+k}}\\
&&+ \brac{\fracm{4}}^j\ \varepsilon \ 2K\ \sum_{i=0}^\infty q^{i}\ \psi\brac{2^{-K(l-2j-1)+i-K}}\\

&=& \brac{\fracm{4}}^j \varepsilon \brac{1 + \fracm{1-q}}\ \Phi\brac{2^{-K(l-2(j+1))}}
 + \brac{\fracm{4}}^j \brac{\varepsilon+q^{2K}} \sum_{k=0}^\infty q^k \psi\brac{2^{-K(l-2(j+1)+1)+k}}\\
&&+ \brac{\fracm{4}}^j\ \varepsilon \ 2K\ q^K \sum_{i=0}^\infty q^{i}\ \psi\brac{2^{-K(l-2(j+1)+1)+i}}
+ \brac{\fracm{4}}^j\ \varepsilon \ 2K\ \sum_{i=0}^{K-1} q^{i}\ \psi\brac{2^{-K(l-2(j+1)+1)+i-K}}\\

&\overset{\eqref{eq:it:psiPhi}}{\leq}& \brac{\fracm{4}}^j \varepsilon \brac{1 + \fracm{1-q}}\ \Phi\brac{2^{-K(l-2(j+1))}}
+ \brac{\fracm{4}}^j \brac{\varepsilon+q^{2K}+\varepsilon 2K q^K}\ \sum_{k=0}^\infty q^k \psi\brac{2^{-K(l-2(j+1)+1)+k}}\\
&&+ \brac{\fracm{4}}^j\ \varepsilon \ 2K\ \sum_{i=0}^{K-1} q^{i}\ \Phi\brac{2^{-K(l-2(j+1))+i-2K}}\\

&\leq& \brac{\fracm{4}}^j \varepsilon \brac{1 + \fracm{1-q}}\ \Phi\brac{2^{-K(l-2(j+1))}} + \brac{\fracm{4}}^j \brac{\varepsilon+q^{2K}+\varepsilon 2K q^K}\ \sum_{k=0}^\infty q^k \psi\brac{2^{-K(l-2(j+1)+1)+k}}\\
&&+ \brac{\fracm{4}}^j\ \varepsilon \ 2K\ \fracm{1-q}\ \Phi\brac{2^{-K(l-2(j+1))}}\\
\end{ma}
\]
which finally gives
\[
\begin{ma}
\Phi\brac{2^{-Kl}} &\leq& \brac{\fracm{4}}^j \varepsilon \brac{1 + \fracm{1-q}+ \ 2K\ \fracm{1-q}}\ \Phi\brac{2^{-K(l-2(j+1))}}\\
&&+ \brac{\fracm{4}}^j \brac{\varepsilon+q^{2K}+\varepsilon 2K q^K}\ \sum_{k=0}^\infty q^k \psi\brac{2^{-K(l-2(j+1)+1)+k}}\\
&\overset{\eqref{eq:it:smallness}}{\leq}& \brac{\fracm{4}}^{j+1} \Phi\brac{2^{-K(l-2(j+1))}}\\
&&+ \brac{\fracm{4}}^{j+1} \ \sum_{k=0}^\infty q^k \psi\brac{2^{-K(l-2(j+1)+1)+k}}.
\end{ma}
\]
Consequently, the claim is proven. In particular, if $\lambda \in (2^{-l-1},2^{-l})$ for some $l \in \N$,
\[
 \Phi(\lambda) \leq \brac{\fracm{4}}^{\left \lfloor \frac{l-1}{2} \right \rfloor} \brac{\Phi(\infty) + G} \leq 32\ \lambda \brac{\Phi(\infty) + G}
\]
\end{proofL}